%% file: arxiv_version.tex
\documentclass{article}

\usepackage[round]{natbib}

\usepackage{fullpage}
\usepackage[utf8]{inputenc} 
\usepackage[T1]{fontenc}    
\usepackage{hyperref}       
\usepackage{url}            
\usepackage{booktabs}       
\usepackage{amsfonts}       
\usepackage{nicefrac}       
\usepackage{microtype}      
\usepackage{xcolor}         
\usepackage{amsmath,amssymb}
\usepackage{bm}
\usepackage{lmodern}
\usepackage{comment}
\usepackage{mathrsfs}
\usepackage{graphicx}
\usepackage{dsfont}
\usepackage{subcaption}
\usepackage{enumitem}
\usepackage{xr}
\usepackage{mathtools}
\usepackage{amsthm}
\usepackage[title]{appendix}

\usepackage{extarrows}

\renewcommand{\qed}{\hfill $\blacksquare$}

\newtheorem{theorem}{Theorem}[section]
\newtheorem{proposition}[theorem]{Proposition}
\newtheorem{lemma}[theorem]{Lemma}
\newtheorem{corollary}[theorem]{Corollary}
\newtheorem{assumption}[theorem]{Assumption}
\theoremstyle{remark}
\newtheorem{claim}[theorem]{Claim}
\newtheorem{remark}[theorem]{Remark}

\begin{document}

%

%

\title{On the estimation of  persistence intensity\\ functions and linear representations\\ of persistence diagrams}

\author{Weichen Wu$^1$ \and Jisu Kim$^2$ \and Alessandro Rinaldo$^3$ }
\date{$^1$ Carnegie Mellon University \\ 
$^2$ Seoul National University\\
$^3$ The University of Texas at Austin}

\maketitle

\begin{abstract}
Persistence diagrams are one of the most popular types of data summaries used in Topological Data Analysis.
The prevailing statistical approach to analyzing persistence diagrams is concerned with filtering out topological noise.
In this paper, we adopt a different viewpoint and 
aim at estimating the actual distribution of a random persistence diagram, which captures both topological signal and noise. To that effect, \cite{chazal2018density}  proved that, under general conditions, the expected value of a random persistence diagram is a measure admitting a Lebesgue density, called the persistence intensity function. In this paper, we are concerned with estimating the persistence intensity function and a novel, normalized version of it -- called the persistence density function. We present a class of kernel-based estimators based on an i.i.d. sample of persistence diagrams and derive estimation rates in the supremum norm. As a direct corollary, we obtain uniform consistency rates for estimating linear representations of persistence diagrams, including Betti numbers and persistence surfaces. Interestingly, the persistence density function delivers stronger statistical guarantees. 

\end{abstract}

\vspace{22pt}
\input{intro.tex}

\input{problem.tex}

\input{results.tex}

\input{discussion.tex}

\section*{Acknowledgements}
Alessandro Rinaldo and Weichen Wu were partially supported by NIH Grant R01 NS121913.

\newpage

\bibliography{refs.bib}
\bibliographystyle{plainnat}

\newpage
\begin{appendices}
\input{basics.tex}
\input{support.tex}

\input{preliminary.tex}
\input{proofs.tex}
\input{plots.tex}
\end{appendices}

\end{document}


%

%

\onecolumn
\aistatstitle{On the estimation of persistence intensity functions and linear representations of persistence diagrams: \\
Supplementary Materials}

\section{FORMATTING INSTRUCTIONS}

To prepare a supplementary pdf file, we ask the authors to use \texttt{aistats2024.sty} as a style file and to follow the same formatting instructions as in the main paper.
The only difference is that the supplementary material must be in a \emph{single-column} format.
You can use \texttt{supplement.tex} in our starter pack as a starting point, or append the supplementary content to the main paper and split the final PDF into two separate files.

Note that reviewers are under no obligation to examine your supplementary material.

\section{MISSING PROOFS}

The supplementary materials may contain detailed proofs of the results that are missing in the main paper.

\subsection{Proof of Lemma 3}

\textit{In this section, we present the detailed proof of Lemma 3 and then [ ... ]}

\section{ADDITIONAL EXPERIMENTS}

If you have additional experimental results, you may include them in the supplementary materials.

\subsection{The Effect of Regularization Parameter}

\textit{Our algorithm depends on the regularization parameter $\lambda$. Figure 1 below illustrates the effect of this parameter on the performance of our algorithm. As we can see, [ ... ]}

\vfill

%% file: intro.tex
\emph{Topological Data Analysis} (TDA) is a field at the interface of computational geometry, algebraic topology and data science whose primary objective is to extract topological and geometric features from possibly high-dimensional, noisy and/or incomplete data. See \cite{fredmichelreview} and references therein for a recent review.
The literature on the statistical analysis of TDA summaries has primarily focused on separating topological signatures from the unavoidable topological noise resulting from the data sampling process. In most cases, the primary goal of statistical inference methods for TDA is to isolate points on the sample persistence diagrams that are sufficiently far from the diagonal to be deemed statistically significant, in the sense of expressing underlying topological signals rather than randomness. This paradigm is entirely natural when the target of inference is one unobservable persistence diagram, and the sample persistent diagrams are noisy approximations to it. Towards that goal, practitioners can now deploy a variety of statistical techniques for identifying topological signals and removing topological noise with provable theoretical guarantees. 

On the other hand, empirical evidence has also demonstrated that topological noise is not necessarily unstructured or uninformative and, in fact, may also carry expressive and discriminative power that can be leveraged for various machine-learning tasks.  In some applications,  the distribution of the topological noise itself is of interest;
in cosmology, see e.g., 
\cite{10.1093/mnras/stab2326}. As a result, statistical summaries able to express the properties of both topological signal and topological noise in a unified manner have also been proposed and investigated: e.g., persistence images and linear functional of the persistence diagrams \citep{adams2017persistence}.

Recently, \cite{chazal2018density}  derived sufficient conditions to ensure that the expected persistent measure -- the expected value of the random counting measure corresponding to a noisy persistent diagram -- admits a Lebesgue density, hereafter the {\it persistence intensity function.} 
The significance of this result is multifaceted. First, the persistent intensity function provides an explicit and interpretable representation of the entire distribution of the persistence homology of random filtrations. 
Secondly,  it allows for a straightforward calculation of the expected linear representation of a persistent diagram as a Lebesgue integral. 
Finally, the representation by the persistence intensity function is of functional, as opposed to algebraic, nature and thus can be estimated via well-established theories and methods from the non-parametric statistics functional estimation. Indeed, \cite{chazal2018density}  further analyzed a kernel-based estimator of the persistence intensity function computed using a sample of i.i.d. persistence diagrams and proved its $L_2$ consistency. Similar results were previously established by  \cite{original.intensity}.

In this paper, we derive consistency rates of estimation of the persistence intensity function and of a novel variant called persistence density function in the $\ell_\infty$ norm based on a sample of i.i.d. persistent diagrams. As we argue below in Theorem~\ref{thm:OT-linfty}, controlling the estimation error for the persistence intensity function in the $\ell_\infty$ norm is stronger than controlling the optimal transport measure $\mathsf{OT}_q$ for any $q>0$ and, under mild assumptions, immediately implies uniform control and concentration of any bounded linear representation of the persistence diagram including (persistent) Betti numbers and persistence images. Our analysis and results are different and complement the statistical results of \cite{chazal2018density} and \cite{original.intensity}; in particular, we seek finite sample $\ell_\infty$ estimation guarantees, which are more challenging and demand more sophisticated techniques.

\begin{figure}[h!]
\centering

\begin{subfigure}{0.48\textwidth}
\centering
\includegraphics[width = \linewidth]{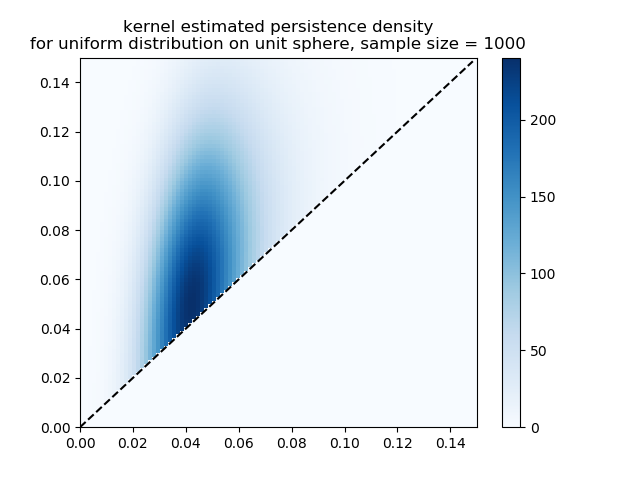}
\end{subfigure}
\hfill 
\begin{subfigure}{0.48\textwidth}
\centering
\includegraphics[width = \linewidth]{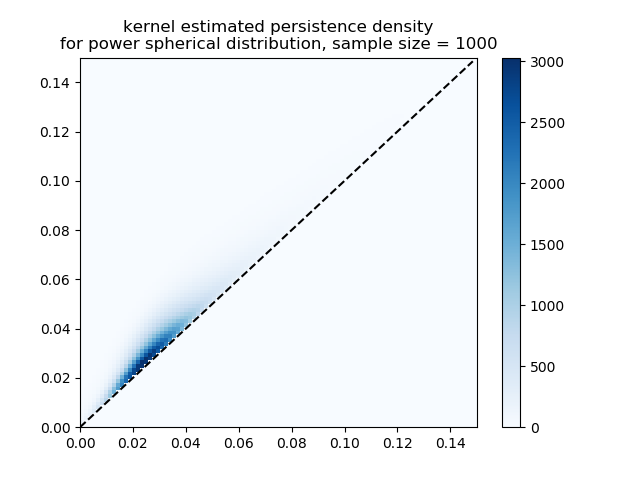}
\end{subfigure}
\caption{Estimated persistence density functions for the uniform distribution (top) and power spherical distribution \cite{decao2020power} (bottom) on the unit circle $\mathbb{S}^1$, based on 1000 diagrams. The parameters for the power spherical distribution are set to $\mu = \frac{\pi}{2}$ and $\kappa = 1$. Each sample contains $1000$ points generated $i.i.d.$ from the distributions on the unit circle corrupted by a $N(0,0.05^2\bm{I}_2)$ additive error.}
\label{fig:sphere}
\end{figure}

We emphasize that the approach and methods considered in this paper are distinct from the prevailing practices in statistical inference for TDA, which focus on extracting topological signals. In contrast, we are interested in capturing the overall randomness of persistence diagrams and describing the topological noise arising from sampling. We aim to develop and explore a statistically grounded approach whose main objective is to describe the {\it distribution} of a random persistence diagram, not any particular realization of it or a target persistence diagram. A second notable point of departure from mainstream TDA is that we assume the availability of an i.i.d. sample of random persistent diagrams, and the accuracy of our rates improve as the number of persistent diagrams increases, not the size of the data used to compute each diagram, which we hold fixed. In both regards, our approach is rather separate from the current TDA paradigm and is not intended as an alternative framework. As an illustrative example, suppose that we are interested in the distribution of persistence diagrams originating from a uniform distribution and a non-uniform distribution on the unit sphere, e.g. the power spherical distribution \citep{decao2020power}. The target topological signature (corresponding to the homology of the unit sphere) is the same in both cases   but the topological noise is different. Figure~\ref{fig:sphere} shows the estimated persistence intensities (for the Vietoris-Rips filtration) based on a sample of 1000 i.i.d. persistence diagrams. The difference in the distributions of the topological noise is apparent and the target of this work. Furthermore, such difference would not be apparent from inspecting individual persistence diagrams; see Figure~\ref{fig:sphere_app} in the supplementary material.


%% file: problem.tex
\section{Background and definitions}

In this section we introduce fundamental concepts from TDA that we will use throughout. We refer the reader to \citep{DBLP:journals/frai/ChazalM21, chazal2018density} for background and extensive references.

\paragraph{Persistence diagrams.}
A persistence diagram is a locally finite multiset of points $D = \{\bm{r}_i = (b_i,d_i) \mid 1 \leq i \leq N(D)\}$ belonging to the set 
\begin{align}\label{eq:defn-Omega-1}
\Omega = \Omega(L) = \{(b,d)\mid 0 < b < d \leq L\} \subset \mathbb{R}^2,
\end{align}
consisting of all the points on the plane in the positive orthant above the identity line and of magnitude no larger than a given constant $L>0$. 
The restriction that the persistence diagrams be contained in a box of side length $L$ is a technical assumption is widely used in the TDA literature; see \cite{divol2021estimation} and the discussion therein. To simplify our notation, we will omit the dependence on $L$, but we will keep track of this parameter in our error bounds.
Some related quantities used throughout are
\begin{align}\label{eq:defn-Omega-2}
&\partial \Omega \coloneqq \{(x,x)\mid 0 \leq x \leq L\}; \quad \quad \overline{\Omega}\coloneqq \Omega \cup \partial \Omega; \nonumber \\ 
&\Omega_\ell: =\left\{\bm{\omega} \in \partial \Omega \big | \|\bm{x} - \partial \Omega\|_2 := \min_{x \in \Omega} \|\bm{\omega} - \bm{x}\|_2 \geq \ell\right\}, \nonumber \\
& \text{for} \quad  \ell \in (0,L/\sqrt{2}).
\end{align}
That is, $\partial \Omega $ is a segment on the diagonal in $\mathbb{R}^2$ and $\Omega_\ell$ consists of all the points in $\Omega$ at a Euclidean distance of $\ell $ or smaller from it.

\paragraph{The expected persistent measure and its normalization.}

A persistence diagram $D = \{\bm{r}_i = (b_i,d_i) \in \Omega \mid 1 \leq i \leq N(D)\}$ can be equivalently represented as a counting measure $\mu$ on $\Omega$  given by
\begin{align*}
A \in \mathcal{B} \mapsto \mu(A) = \sum_{i=1}^{N(D)} \delta_{\bm{r}_i}(A),
\end{align*}
where $\mathcal{B} = \mathcal{B}(\Omega)$ is the class of all Borel subsets of $\Omega$ and $\delta_{\bm{r}}$ denotes the Dirac point mass at $\bm{r} \in \Omega$. We will refer to $\mu$ as the {\it persistence measure} corresponding to $D$ and, with a slight abuse of notation, will treat persistence diagrams as counting measures. If $D$ is a random persistence diagram, then the associated persistence measure is also random. 
We will also study its {\it normalized measure}  $\tilde{\mu}$, which is the persistence measure divided by the total number of points $N(D)$ in the persistence diagram: 
\begin{align*}
A \in \mathcal{B} \mapsto \tilde{\mu}(A) = \frac{1}{N(D)}\sum_{i=1}^{N(D)} \delta_{\bm{r}_i}(A).
\end{align*}
The normalized persistence measure may be more appropriate when the number of points $N(D)$ in the persistence diagram is not of direct interest but their spatial distribution is. This is typically the case when the persistence diagrams at hand contain many points or are obtained from large random filtrations (e.g. the Vietoris-Rips complex built on point clouds), so that the value of $N(D)$ will mostly account for noisy topological fluctuations due to sampling.

We will consider the setting in which the observed persistence diagram $D$ is a random draw from an unknown distribution. Then, the (non-random) measures
\begin{align*}
&A \in \mathcal{B} \mapsto \mathbb{E}[\mu](A) = \mathbb{E}[\mu(A)] \quad \text{and} \\ 
&A \in \mathcal{B} \mapsto \mathbb{E}[\tilde{\mu}](A) = \mathbb{E}[\tilde{\mu}(A)]
\end{align*}
are well defined. We will refer to $\mathbb{E}[\mu]$ and $\mathbb{E}[\tilde{\mu}]$
as the {\it expected persistence measure} and the {\it expected  persistence probability,} respectively.  Neither is a discrete measure (even though persistence measures are discrete by definitions).  Of course, the expected persistence probability $\mathbb{E}[\tilde{\mu}]$ is a probability measure. 

The interpretations of $\mathbb{E}[\mu]$ and of $\mathbb{E}[\tilde{\mu}]$ are straightforward: for any Borel set $A \subset \Omega$, $\mathbb{E}[\mu](A)$ is the expected number of points from the random persistence diagram falling in $A$, while $\mathbb{E}[\tilde{\mu}](A)$ is the probability that a random persistence diagram will intersect $A$. 
Despite their interpretability, the expected persistence measure and probability are not yet standard concepts in the practice and theory of TDA. As a result, they have not been thoroughly investigated.


\paragraph{The persistence intensity and density functions and linear representations.} Recently, \cite{chazal2018density} derived conditions  -- applicable to a wide range to problems -- that ensure that the expected persistence measure $\mathbb{E}[\mu]$ and its normalization $\mathbb{E}[\tilde{\mu}]$ both admit densities with respect to the Lebesgue measure on $\Omega$. Specifically, under fairly mild and general conditions detailed in \cite{chazal2018density} there exist measurable functions $p:\Omega \to \mathbb{R}_{\geq 0}$ and $\tilde{p}:\Omega \to \mathbb{R}_{\geq 0}$, such that for any Borel set $A \subset \Omega$,
\begin{align}\label{eq:intens.dens}
\mathbb{E}[\mu](A) = \int_{A} p(\bm{u}) \mathrm{d}\bm{u}, \quad  \mathbb{E}[\tilde{\mu}](A) = \int_{A}\tilde{p}(\bm{u}) \mathrm{d}\bm{u}.
\end{align}
 In fact, \cite{chazal2018density} provided explicit expressions for $p$ and $\tilde{p}$ (see Section~\ref{app:proof-bound}). Notice that, by construction, $\tilde{p}$ integrates to $1$ over $\Omega$. 
We will refer to the functions $p$ and $\tilde{p}$ as the {\it persistence intensity} and the {\it persistence density} functions, respectively. We remark that the notion of a persistence intensity function was originally put forward by \cite{original.intensity}. 

The persistence intensity and density functions ``operationalize'' the notions of expected persistence measure and expected persistence probability introduced above, allowing us to evaluate, for any set $A \in \mathcal{B}$, $\mathbb{E}[\mu](A)$ and $\mathbb{E}[\tilde{\mu}](A)$ in a straightforward way as Lebesgue integrals. 
The main objective of the paper is to construct estimators of the persistence intensity $p$ and persistence density $\tilde{p}$, respectively, and to provide high probability error bounds with respect to the $L_\infty$ norm. 
As we show below in Theorem~\ref{thm:OT-linfty},$L_\infty$-consistency for the persistence intensity function is a stronger guarantee than consistency in the $\mathsf{OT}_p$ metric, for any $p < \infty$. 

As noted in \cite{chazal2018density}, the persistence intensity and density functions are naturally suited to compute the expected value of linear representations of random persistence diagrams. A linear representation $\Psi$ of the persistence diagram $D = \{\bm{r}_i = (b_i,d_i) \in \Omega \mid 1 \leq i \leq N(D)\}$ with corresponding persistence measure $\mu$ is a summary statistic of $D$ of the form
\begin{align}\label{eq:linear-represent}
\Psi(D) = \sum_{i=1}^{N(D)} f(\bm{r}_i ) = \int_\Omega f(\bm{u}) d \mu(\bm{u}),
\end{align}
for a given measurable function $f$ on $\Omega$. (An analogous definition can be given for the normalized persistence measure $\tilde{\mu}$ instead). Then,
\begin{align}\label{eq:expected-linear-represent}
\mathbb{E}[\Psi(D)] = \int_\Omega f(\bm{u}) d \mathbb{E}[\mu](\bm{u})= \int_\Omega f(\bm{u}) p(\bm{u}) d \bm{u},
\end{align}
where the second identity follows from \eqref{eq:intens.dens}. Linear representations include persistent Betti numbers, persistence surfaces \citep{adams2017persistence}, persistence silhouettes, \cite{Chazal2013StochasticCO} and persistence weighted Gaussian kernels \citep{pmlr-v48-kusano16}.
The persistence surface is an especially popular linear representation. In detail, for a \emph{kernel function} 
$K(\cdot):\mathbb{R}^2 \to \mathbb{R}$
 and any $\bm{x} \in \mathbb{R}^2$, let $K_h(\bm{x}) = \frac{1}{h^2}K(\frac{\bm{x}}{h})$, where $h > 0$ is the bandwidth parameter\footnote{\citep{adams2017persistence} showed empirically that the bandwidth does not have a major influence on the efficiency of the persistence surface.}. The persistence surface of a persistence measure $\mu$ is defined as
\begin{align}\label{eq:defn-pers-surface}
\rho_h(\bm{u}) = \int_{\Omega} f(\bm{\omega}) K_h(\bm{u}-\bm{\omega}) \mathrm{d}\mu(\bm{\omega}),
\end{align}
where $f(\bm{\omega}) \colon \mathbb{R}^2 \rightarrow \mathbb{R} $ is the user-defined \textit{weighting function,} chosen to ensure stability of the representation. Our analysis allows us to immediately obtain consistency rates for the expected persistence surface in $L_\infty$ norm; see Theorem~\ref{thm:surface-var} in the supplementary material. 


\paragraph{Betti and the persistent Betti numbers.}
The \textbf{Betti number} at scale $x \in [0,L]$ is the number of persistent homologies that are in existence at ``time" $x$. Furthermore, the \textbf{persistent Betti number} at a certain point $\bm{x} = (x_1,x_2) \in \Omega$  measures the number of persistent homologies that are born before $x_1$ and die after $x_2$. In our notation,  given a  persistence diagram $D$ and its associated persistence measure $\mu$, for $x \in [0,L]$ and $\bm{x} = (x_1,x_2) \in \Omega$, the corresponding  Betti number and persistent Betti number  are given by
\begin{align*}
\beta_x 
= \mu(B_x) \quad \text{and} \quad \beta_{\bm{x}}  
= \mu(B_{\bm{x}}),
\end{align*}
respectively, where $B_x = [0,x) \times (x,L]$ and $B_{\bm{x}} = [0,x_1) \times (x_2,L]$.
Though Betti numbers are among the most prominent and widely used TDA summaries, relatively little is known about the statistical hardness of estimating their expected values when the sample size is fixed and the number of persistence diagrams increases. Our results will yield error bounds of this type.
We will also consider normalized versions of the Betti numbers defined using the persistence probability $\tilde{\mu}$ of the persistence diagram: 
\begin{align*}
\tilde{\beta}_x = \tilde{\mu}(B_x) \quad \text{and} \quad \tilde{\beta}_{\bm{x}} =  \tilde{\mu}(B_{\bm{x}}).
\end{align*}
Notice that, by definition, $\tilde{\beta}_{\bm{x}} \leq 1$. While their interpretation is not as direct as the Betti numbers computed using persistence diagrams, the expected normalized (persistence) Betti numbers $\mathbb{E}[\tilde{\beta}_{\bm{x}}] = \mathbb{E}[\tilde{\mu}](B_{\bm{x}})$ are informative topological summaries while showing favorable statistical properties (see Corollary~\ref{thm:pers-betti-kernel-var} below).

%% file: results.tex
\section{Results}\label{sec:results}


\subsection{On the OT distance and the $L_{\infty}$ distance between intensity functions} 
We first show that the topology induced by the $L_\infty$ distance between intensity functions is stronger than the one corresponding to the \emph{optimal transport} distance, a natural and very popular metric for persistence diagrams -- and, more generally, locally finite Radon measures such as normalized persistence measures and probabilities; see, in particular, \citep{divol2021estimation}. 
In detail, for two Radon measures $\mu$ and $\nu$ supported on $\overline{\Omega}$, an \emph{admissible transport} from $\mu$ to $\nu$ is defined as a function $\pi:\overline{\Omega} \times \overline{\Omega} \to \mathbb{R}$, such that for any Borel sets $A,B \subset \overline{\Omega}$,
\begin{align*}
\pi(A \times \overline{\Omega}) = \mu(A), \quad \text{and} \quad \pi(\overline{\Omega} \times B) = \nu(B). 
\end{align*}
Let $\mathsf{adm}(\mu,\nu)$ denote all the admissible transports from $\mu$ to $\nu$. For any $q \in \mathbb{R}^+ \cup \{\infty\}$, the $q$-th order Optimal Transport (OT) distance between $\mu$ and $\nu$ is defined as
\begin{align*}
\mathsf{OT}_q(\mu,\nu) = \left(\inf_{\pi \in \mathsf{adm}(\mu,\nu)} \int_{\overline{\Omega} \times \overline{\Omega}} \|\bm{x} - \bm{y}\|_2^q \mathrm{d} \pi(\bm{x},\bm{y})\right)^{\frac{1}{q}}.
\end{align*}
The OT distance is widely used for good reasons: by transporting from and to the diagonal $\partial \Omega$, it captures the distance between two measures that have potentially different total masses, taking advantage of the fact that points on the diagonal have arbitrary multiplicity in persistent diagrams. It also proves to be stable with respect to perturbations of the input to TDA algorithms. It turns out that the $L_{\infty}$ distance between intensity functions provides a tighter control on the difference between two persistent measures. Below, for a real-valued  function on $\Omega$, we let $\| f \|_\infty = \sup_{\bm{x} \in \Omega }|f(\bm{x})|$ be its $L_\infty$ norm.

\begin{theorem}\label{thm:OT-linfty}
 Let $\mu$, $\nu$ be two expected persistent measures on $\Omega$ with intensity functions $p_{\mu}$ and $p_{\nu}$ respectively.  Then
\begin{align*}
\mathsf{OT}_q^q(\mu,\nu) \leq \frac{2}{(q+1)(q+2)} \left(\frac{L}{\sqrt{2}}\right)^{q+2}\|p_{\mu} - p_{\nu}\|_\infty.
\end{align*}

Furthermore, there exist two sequences of expected persistence measures $\{\mu_n\}_{n \in \mathbb{N}}$ and $\{\nu_n\}_{n \in \mathbb{N}}$ with intensity functions $\{p_{\mu_n}\}_{n \in \mathbb{N}}$ and $\{p_{\nu_n}\}_{n \in \mathbb{N}}$ respectively such that, as $n \to \infty$, 
\begin{align*}
\mathsf{OT}_q(\mu_n,\nu_n) \to 0, \quad  \text{while}  \quad \|p_{\mu_n} - p_{\nu_n}\|_\infty \to \infty.
\end{align*}
\end{theorem}

\paragraph{The bottleneck distance} For the case of $q = \infty$, which corresponds to the bottleneck distance when applied to persistence diagrams, there can be no meaningful upper bound of the form of  Theorem~\ref{thm:OT-linfty}: we show in Section~\ref{sec:D1} of the supplementary material that there exist two sequences of measures such that their bottleneck distance converges to a finite number while the $L_{\infty}$ distance between their intensity functions vanishes. Existing contributions in the optimal transport literature \citep{peyre2018comparison,nietert2021smooth} also upper bound the optimal transport distance by a Sobolev-type distance between density functions. Notably, these bounds require, among other things, the measures to have common support and the same total mass, two conditions not assumed in Theorem~\ref{thm:OT-linfty}.

\subsection{Non-parametric estimation of the persistent intensity and density functions}
In this section, we analyze the performance of kernel-based estimators of the persistent intensity function and the persistent density function in the same setting considered by \cite{chazal2018density} and \cite{original.intensity}, where we observe $n$ $i.i.d.$ persistent measures (i.e. diagrams) $\mu_1,\mu_2,\ldots, \mu_n$. 
The proposed procedures are inspired by kernel density estimators for probability densities traditionally used in the non-parametric statistics literature; see, e.g., \cite{gine2021mathematical}. Specifically, we consider the following estimator for $p$ and $\tilde{p}$, respectively:
\begin{subequations}\label{eq:kernel-est-p}
\begin{align}
&\bm{\omega} \in \mathbb{R}^2 \mapsto \hat{p}_h(\bm{\omega}) \coloneqq \frac{1}{n} \sum_{i=1}^n \int_{\Omega} K_{h}(\bm{x}-\bm{\omega})\mathrm{d}\mu_i(\bm{x});
\label{eq:defn-phat}\\
&\bm{\omega} \in \mathbb{R}^2 \mapsto \check{p}_h(\bm{\omega}) =  \frac{1}{n} \sum_{i=1}^n \int_{\Omega} K_{h}(\bm{x}-\bm{\omega})\mathrm{d}\tilde{\mu}_i(\bm{x}),\label{eq:defn-pcheck}
\end{align}
\end{subequations}
where $K(\cdot)$ is a {\it kernel function,} which we assume to satisfy standard conditions used in non-parametric literature, discussed in detail in Section~\ref{sec:assumptions} of the supplementary material.

 \paragraph{Assumptions.}
We require several regularity conditions. Furthermore, we will implicitly assume throughout that both $p$ and $\tilde{p}$ (see \ref{eq:intens.dens}) are well-defined as densities with respect to the Lebesgue measure, though this is not strictly necessary for our main results, Theorems~\ref{thm:intensity-var} and \ref{thm:F-var}.


First, we assume a uniform bound on the $q$-th order total persistence of the $\mu_i$'s,  though not on the total number of points in the persistence diagrams. This can be thought of as a basic moment existence condition that, as elucidated in  \cite{cohen2010lipschitz} and discussed in \cite{divol2019choice} and \cite{divol2021estimation}, is a relatively mild assumption satisfied by a broad variety of data-generating mechanisms. 
\begin{assumption}[Bounded total persistence]\label{as:total-persistence}
There exists a constant $M>0$, such that, for the value of $q$ as in Assumption~\ref{as:bound}, it holds that, almost surely,
\[
\max_{i =1,\ldots,n}\int_{\Omega} \|\bm{\omega} - \partial \Omega\|_2^q \mathrm{d}\mu_i(\omega) < M.
\]
\end{assumption}
 
We will denote with $\mathcal{Z}_{L,M}^q$ the set of persistent measures on $\Omega_{L}$ satisfying Assumption~\ref{as:total-persistence}.

Next, we impose key boundedness conditions on $p$ and $\tilde{p}$, which are needed to apply a concentration inequality for empirical processes that deliver uniform control over the variances of these estimators.

\begin{assumption}[Boundedness]\label{as:bound}
For some $q > 0$, let $\bar{p}(\bm{\omega}) \coloneqq \|\bm{\omega} - \partial \Omega\|_2^q p(\bm{\omega})$. Then, 
\begin{align*}
&\|\bar{p}\|_{\infty} = \sup_{\omega \in \Omega} \|\bm{\omega} - \partial \Omega\|_2^q p(\bm{\omega}) < \infty \quad \text{and} \\ 
&\|\tilde{p}\|_{\infty} = \sup_{\omega \in \Omega} \tilde{p}(\bm{\omega}) < \infty.
\end{align*}
\end{assumption}
We remark that a bound on the  $L_\infty$ norm of the intensity function $p$ is not a realistic assumption because the total mass of the persistence measure may not be uniformly bounded in several common data-generating mechanisms. Indeed, in light of existing results, this condition would be likely violated in many scenarios; see e.g. \cite{divol2019choice}). Thus, we  only require that the weighted intensity function $\bar{p}(\bm{\omega}) = \|\bm{\omega} - \partial \Omega\|_2^q p(\bm{\omega})$ has finite $L_\infty$ norm. Still, it is not a priori clear that Assumption~\ref{as:bound} for $\bar{p}$ itself is realistic; in the supplementary material, we prove that this is indeed the case for the Vietoris-Rips filtration built on i.i.d. samples. On the other hand, assuming that the persistence density  $\tilde{p}$ is uniformly bounded poses no problems. For a formal argument, see Theorems~\ref{thm:bound-intensity} and \ref{thm:bound-density} in the supplementary material. This fact is the primary reason why the persistence probability density function -- unlike the persistence intensity function --  can be estimated uniformly well over the entire set $\Omega$ - see \eqref{thm:intensity-var} below.  We refer readers to Section~\ref{sec:validation} of the supplementary materials for details and a discussion on this subtle but consequential point.

One of the main results of the paper are high probability uniform bounds on the fluctuations of the kernel estimators around their expected values. For a fixed value of the bandwidth $h$, they imply that the estimators $\hat{p}_h$ and $\check{p}_h$ concentrate around their expected value at a parametric rate $1/\sqrt{n}$.

\begin{theorem}\label{thm:intensity-var}
 Suppose that Assumptions~\ref{as:total-persistence} and\ref{as:bound} hold. Then,
 \begin{enumerate}[label=(\alph*)]
 \item there exist positive constants $C_1, C_2$ depending on $M,\|K\|_{\infty}, \|K\|_2, \|\bar{p}\|_{\infty}$ and $q$ such that for any $\delta \in (0,1)$,  it can be guaranteed with probability at least $1-\delta$ that
\begin{align*}
    &\sup_{\bm{\omega}\in \Omega_{2h}} \ell_{\bm{\omega}}^q |\hat{p}_h(\bm{\omega}) - \mathbb{E}\hat{p}_h(\bm{\omega})| \\ 
    &\leq   \max \left\{C_1\frac{1}{nh^2} \log \frac{1}{\delta h^2} , C_2\sqrt{\frac{1}{nh^2}}\sqrt{\log \frac{1}{\delta h^2}}\right\},
\end{align*}
where $\ell_{\omega} \coloneqq \|\bm{\omega} - \partial \Omega \|_2 - h;$ 
\item there exist positive constants $C_1,C_2$ depending on $M,\|K\|_{\infty}, \|K\|_2, \|\tilde{p}\|_{\infty}$ and $q$ such that for any $\delta \in (0,1)$,  it can be guaranteed with probability at least $1-\delta$ that
\begin{align*}
    &\sup_{\bm{\omega}\in \Omega} |\check{p}_h(\bm{\omega}) - \mathbb{E}\check{p}_h(\bm{\omega})|\\ 
    &\leq   \max \left\{C_1\frac{1}{nh^2} \log \frac{1}{\delta h^2} , C_2\sqrt{\frac{1}{nh^2}}\sqrt{\log \frac{1}{\delta h^2}}\right\}.
\end{align*}
\end{enumerate}
\end{theorem}
\paragraph{Remark.} The dependence of the constants on problem related parameters is made explicit in the proofs; see the supplementary material. 

 There is an important difference between the two bounds in Theorem~\ref{thm:intensity-var}: while the variation of $\check{p}_h(\bm{\omega})$ is uniformly bounded everywhere on $\Omega$, the variation of $\hat{p}_h(\bm{\omega})$ is uniformly bounded only when $\bm{\omega}$ is at least $2h$ away from the diagonal $\partial \Omega$, and may increase as $\bm{\omega}$ approaches the diagonal. The difficulty in controlling $\hat{p}_h$ near the diagonal stems from the fact that we only assume the total persistence to be bounded; in other words, the number of points near the diagonal in the sample persistent diagrams can be prohibitively large, since their contribution to the total persistence is negligible. This is expected in noisy settings where the sampling process will result in topological noise consisting of many points in the persistence diagram near the diagonal. We do not know whether this limitation of the estimator $\hat{p}_h$ is intrinsic to the problem or instead an artifact of our proof techniques. Nonetheless, the above result suggests that to achieve uniform control over $\Omega$, relying on density-based rather than intensity-based representations of the persistent measures may be preferable. 

\paragraph{Bias-variance trade-off and minimax lower bound. } 

In order to measure how well $\hat{p}_h$ and $\check{p}_h$ concentrate not just around their expectations but around the target densities $p$ and $\tilde{p}$, respectively, we will need to further control their biases, as a function of the bandwidth $h$. To that effect, we require some degree of smoothness of both $p$ and $\tilde{p}$, as it is standard in non-parametric density estimation. We refer the reader to the appendix for the definition of smooth function spaces. 

\begin{assumption}[Smoothness]\label{as:density}
The persistence intensity function $p$ and persistence probability density function $\tilde{p}$ are  H\"older smooth of the order of $s>0$ with parameters $L_p$ and $L_{\tilde{p}}$, respectively. 
\end{assumption}

Using the above assumption and standard arguments, we obtain that, uniformly over $\bm{\omega} \in \Omega$,  $|\mathbb{E}[\hat{p}_h(\bm{\omega})] - p(\bm{\omega})| $ and $|\mathbb{E}[\check{p}_h(\bm{\omega})] - \tilde{p}(\bm{\omega})|$ are both of order $h^2$. See Theorem \ref{thm:intensity-bias} in the appendix. Next, assuming that the number $n$ of persistent diagram grows unbounded, it follows from Theorems~\ref{thm:intensity-bias} and \ref{thm:intensity-var} that setting the bandwidth to be $h \asymp n^{-\frac{1}{2(s+1)}}$ will optimize the bias-variance trade-off, yielding high-probability estimation errors 
\begin{align*}
&\sup_{\bm{w} \in \Omega_{2h}} \ell_{\omega}^q|\hat{p}_h(\bm{\omega}) - {p}(\omega)| \lesssim O\left(n^{-\frac{s}{2(s+1)}}\right), \quad \text{and}\\ 
&\sup_{\bm{w} \in \Omega} |\check{p}_h(\bm{\omega}) - \tilde{p}(\omega)| \lesssim O\left(n^{-\frac{s}{2(s+1)}}\right).
\end{align*}
In our next result, we show that the above rate is minimax optimal for the persistence density function. For brevity, we here omit a similar result for the persistence intensity function (see Theorem~\ref{thm:weight-intensity-minimax} in the supplementary material). 
\begin{theorem}\label{thm:density-minimax}
Let $\mathscr{F}$ denote the set of functions on $\Omega$ with  Besov norm bounded by $B>0$:
\begin{align*}
\mathscr{F} = \{f: \Omega \to \mathbb{R}, \|f\|_{B_{\infty,\infty}^s} \leq B\}.
\end{align*}
Then,
\begin{align*}
\inf_{\check{p}_n} \sup_{P} \mathbb{E}_{\mu_1,\ldots,\mu_n \overset{\text{i.i.d.}}{\sim} P} \|\check{p}_n - \tilde{p}\|_{\infty} \geq O(n^{-\frac{s}{2(s+1)}}),
\end{align*}
where the infimum is taken over estimator $\check{p}_n$  mapping $\mu_1,\ldots,\mu_n$ to an intensity function in $\mathscr{F}$, the supremum is over the set of all probability distributions on $\mathcal{Z}_{L,M}^q$ and $\tilde{p}$ is the intensity function of $\mathbb{E}_P[\tilde{\mu}]$.
\end{theorem}

\subsection{Kernel-based estimators for linear functionals of the persistent measure}\label{sec:linear}
The kernel estimators \eqref{eq:kernel-est-p} can serve as a basis for estimating bounded linear representations of the expected persistence measure $\mathbb{E}[\mu]$ and its normalized counterpart $\mathbb{E}[\tilde{\mu}]$. Specifically, for $R>0$, let $\mathscr{F}_{2h,R}$ and $\widetilde{\mathscr{F}}_{R}$ denote the set of linear representations of the form
\begin{align*}
\mathscr{F}_{2h,R} = \bigg\{ &\Psi = \int_{\Omega_{2h}} f \mathrm{d} \mathbb{E}[\mu] \bigg | f: \Omega_{2h} \to \mathbb{R}_{\geq 0}, \\ 
&\int_{\Omega_{2h}} \ell_{\bm{\omega}}^{-q} f(\bm{\omega}) \mathrm{d} \bm{\omega} \leq R\bigg\}, \quad \text{and} \\
\widetilde{\mathscr{F}}_{R} = \bigg\{ \Psi = &\int_{\Omega} f \mathrm{d} \mathbb{E}[\mu] \bigg | f: \Omega \to \mathbb{R}_{\geq 0}, \int_{\Omega}  f(\bm{\omega}) \mathrm{d} \bm{\omega} \leq R \bigg\}.
\end{align*}
Then, any linear representations $\Psi \in \mathscr{F}_{2h,R}$ and $\widetilde{\Psi} \in \widetilde{\mathscr{F}}_R$ can be estimated by
\begin{align}\label{eq:defn-F-est}
\hat{\Psi}_{h} = \int_{\Omega_{2h}} f(\bm{\omega}) \hat{p}_h(\bm{\omega}) \mathrm{d}\bm{\omega}  \text{ 
 and  }  \check{\Psi}_{h} = \int_{\Omega} f(\bm{\omega}) \check{p}_h(\bm{\omega}) \mathrm{d}\bm{\omega}.
\end{align}

As a direct corollary of Theorem~\ref{thm:intensity-bias}, we obtain the following uniform high-probability bound on the variance of $\hat{\Psi}_{h}$ and  $\check{\Psi}_{h}$, which, for fixed $h$, yield $1/\sqrt{n}$ rates. In the supplementary material, we also show that, not surprisingly, the biases of both estimators are of order $h^s$ under Assumption ~\ref{as:density}; see Theorem \ref{thm:F-bias} in the supplementary material. 

\begin{theorem}\label{thm:F-var}
Assume that Assumptions~\ref{as:total-persistence} and\ref{as:bound} hold. Then,
\begin{enumerate}[label=(\alph*)]
\item there exist constants $C_1,C_2$ depending on $M,\|K\|_{\infty}, \|K\|_2, \|\bar{p}\|_{\infty}$ and $q$ such that for any $\delta \in (0,1)$,  it can be guaranteed with probability at least $1-\delta$ that
\begin{align*}
&\sup_{\Psi \in \mathscr{F}_{2h,R}} \left|\hat{\Psi}_{h} - \mathbb{E}[\hat{\Psi}_{h}]\right|\\
&\leq R \cdot   \max \left\{C_1\frac{1}{nh^2} \log \frac{1}{\delta h^2} , C_2\sqrt{\frac{1}{nh^2}}\sqrt{\log \frac{1}{\delta h^2}}\right\};
\end{align*}
\item there exist constants $C_1,C_2$  depending on $M,\|K\|_{\infty}, \|K\|_2, \|\tilde{p}\|_{\infty}$ and $q$ such that for any $\delta \in (0,1)$,  it can be guaranteed with probability at least $1-\delta$ that
\begin{align*}
&\sup_{\Psi \in \widetilde{\mathscr{F}}_{R}} \left|\check{\Psi}_{h} - \mathbb{E}[\check{\Psi}_{h}]\right|\\
&\leq R \cdot  \max \left\{C_1\frac{1}{nh^2} \log \frac{1}{\delta h^2} , C_2\sqrt{\frac{1}{nh^2}}\sqrt{\log \frac{1}{\delta h^2}}\right\}.
\end{align*}
\end{enumerate}
\end{theorem}

It is important to highlight the fact that the above bounds hold uniformly over the choice of linear representations under only mild integrability assumptions. 
 We again stress the difference between the two upper bounds: part (a) shows that for a linear functional of the original persistent measure to have controlled variation, we need to be at least $2h$ away from the diagonal $\partial \Omega$, a requirement that is not necessary for linear functionals of the normalized persistent measure, as is shown in part (b). 

We apply our results to the analysis of persistence surfaces and persistence Betti numbers. Due to space limitations, in the main text we focus on the latter and refer the reader to  Theorem~\ref{thm:surface-var} in the supplementary material for novel error rates in estimating persistent surfaces.
For any $\bm{x} \in \Omega$, the persistent Betti number $\beta_{\bm{x}}$ can be estimated in a straightforward way by integrating $\hat{p}_h$ or $\check{p}_h$ over $B_{\bm{x}}$ (which we recall we define to be $B_{\bm{x}} = [0,x_1) \times (x_2,L]$): 
\begin{align}\label{eq:kernel-est-pers-betti}
\hat{\beta}_{\bm{x},h} = \int_{B_{\bm{x}}} \hat{p}_h(\bm{\omega}) \mathrm{d} \bm{\omega} \quad \check{\beta}_{\bm{x},h} = \int_{B_{\bm{x}}} \check{p}_h(\bm{\omega}) \mathrm{d} \bm{\omega}.
\end{align}

 An immediate application of \ref{thm:F-var} yields $1/\sqrt{n}$ high-probability concentration rates.
\begin{corollary}\label{thm:pers-betti-kernel-var}
There exist a constants $C_1,C_2$ depending on $M,\|K\|_{\infty}, \|K\|_2, \|\bar{p}\|_{\infty}, \|\tilde{p}\|_\infty$ and $q >2$ such that, for any $\delta \in (0,1)$,
\begin{enumerate}[label=(\alph*)]
\item for the persistent Betti numbers computed using $\hat{p}_h$,
\begin{align*}
&\sup_{ \bm{x} \in \Omega \colon \ell_{\bm{x}}  >h}\ell_{\bm{x}}^{q-2} \left|\hat{\beta}_{\bm{x},h} - \mathbb{E}[\hat{\beta}_{\bm{x},h}]\right|\\
&\leq \max \left\{C_1\frac{1}{nh^2} \log \frac{1}{\delta h^2} , C_2\sqrt{\frac{1}{nh^2}}\sqrt{\log \frac{1}{\delta h^2}}\right\}
\end{align*}
\item for the persistent Betti numbers computed using $\check{p}_h$,
\begin{align*}
&\sup_{\bm{x} \in \Omega } \left|\check{\beta}_{\bm{x},h} - \mathbb{E}[\check{\beta}_{\bm{x},h}]\right|\\
&\leq \frac{L^2}{4}  \max \left\{C_1\frac{1}{nh^2} \log \frac{1}{\delta h^2} , C_2\sqrt{\frac{1}{nh^2}}\sqrt{\log \frac{1}{\delta h^2}}\right\}
\end{align*}
\end{enumerate}
\end{corollary}


It is also straightforward to see that the biases of both $\hat{\beta}_{\bm{x},h}$ and $\check{\beta}_{\bm{x},h}$  are of order $h^s$, uniformly in $\bm{x}$; see Corollary \ref{thm:pers-betti-kernel-bias} in the supplementary material.


As noted before, the concentration rates of the estimator of the persistence Betti numbers based on the persistence density hold uniformly over $\Omega$, thus suggesting that 
the kernel-based estimator $\hat{p}_h$ \emph{will not be guaranteed to} yield a stable estimation of the Betti number $\beta_x$. As remarked above, this issue arises as the intensity function may not be uniformly bounded near the diagonal. Indeed, in the supplementary material, we describe an alternative proof technique based on an extension of the standard VC inequality and arrive at a very similar rate. On the other hand, this issue does not affect the normalized Betti numbers $\tilde{\beta}_{\bm{x}}$.


An important consequence of the previous result is  a {\it uniform} error bound for the {\it expected normalized Betti curve}
\[
x \in (0,L) \mapsto \mathbb{E}[\tilde{\beta}_x] = \mathbb{E}[\tilde{\mu}](B_x) = \int_{B_x} \tilde{p}(\bm{\omega}) d \bm{\omega},
\]
where  $B_x = [0,x) \times (x,L]$ and $\tilde{\mu}$ is the normalized persistent measure corresponding to a random persistence diagram. In detail, for constants $C_1,C_2>0$ depending on the model parameters, with probability at least $ 1- \delta$, 
\begin{align*}
&\sup_{x \in (0,L)} \left|\check{\beta}_{x,h} - \mathbb{E}[\tilde{\beta}_x]]\right|\\
&\leq \max \left\{C_1\frac{1}{nh^2} \log \frac{1}{\delta h^2} , C_2\sqrt{\frac{1}{nh^2}}\sqrt{\log \frac{1}{\delta h^2}}\right\}  + h^s
\end{align*}
To the best of our knowledge, this is the first result of this kind, as typically one can only establish pointwise and not uniform consistency of Betti numbers.

%% file: discussion.tex
\section{Numerical Illustrations} 
To illustrate our methodology and highlight the differences between the persistence intensity and density functions, we compare the persistence intensity and density functions of 1000 data points drawn from the uniform distribution and the power sphere distribution \cite{decao2020power} on the unit circle $\mathbb{S}^1$. The density functions shown in Figure \ref{fig:sphere} illustrate a clear difference between the structure of topological noise generated by the two distributions. We include plots of the sample points, sample persistence diagrams and kernel-based estimators of persistence intensity functions in Section ~\ref{sec:experiments} of the supplementary material.

We also consider the MNIST handwritten digits dataset and the ORBIT5K dataset.  The ORBIT5K dataset contains independent simulations for the linked twist map, dynamical systems for fluid flow as described in \cite{adams2017persistence}; see also  Appendix G.2 of \cite{NEURIPS2020_b803a925}. In Section~\ref{sec:experiments} of the supplementary material, we show the estimated persistence intensity and density functions computed from persistence diagrams obtained over a varying number of random samples from the ORBIT5K datasets, for different model parameters. The figures confirm our theoretical finding that the values of the persistence density function near the diagonal are not as high (on a relative scale) as those of the persistence intensity function. An analogous conclusion can be reached when inspecting the persistence intensity and density functions for different draws of the MNIST datasets for the digits 4 and 8. We further include plots of the average Betti and normalized Betti curves from the ORBIT5K dataset, along with the curves of the empirical point-wise 5\% and 95\% quantiles. These plots reveal the different scales of the Betti curves and normalized Betti curves, and of their uncertainty.

\section{Discussion}
In this paper, we have taken the first step towards developing a new set of methods and theories for statistical inference for TDA based on i.i.d. samples of persistence diagrams. Our main focus is on the estimation of the persistence intensity function \cite{chazal2018density,original.intensity}, a TDA summary of a functional type that encodes the entire distribution of a random persistence diagram and is naturally suited to handle linear representations. We have analyzed a simple kernel estimator and derived uniform consistency rates that hold under very mild assumptions. We also propose the persistence density function, a novel functional TDA summary that enjoys stronger statistical guarantees. 

A notable advantage of deploying persistence intensity and density functions to quantify the difference between distributions of persistence diagrams compared to more traditional approaches based on optimal transport distances is that our methodology is computationally feasible. Indeed, computing kernel-based estimators of the persistence intensity and density functions is a straightforward task even with very large sample sizes, and so is to evaluate any $L_p$ distanced between them. In contrast, computing optimal transport distances between many persistence diagrams is typically computationally prohibitive. 

There remain various open problems worth pursuing. A natural direction is the study of the topology over the space of normalized persistence measures. For example, based on our results from section \ref{sec:results}, one may expect the normalized persistence measure not to be continuous for the vague topology with respect to the Hausdorff distance. Similarly, it would also be interesting to further investigate the topology induced by convergence of the persistence densities in the $L_\infty$ norm. From the statistical side,  our results guarantee the consistency of the proposed estimators. However, in order to carry out statistical inference, it is necessary to develop more sophisticated procedures that quantify the uncertainty of our estimators. Toward that goal, it would be interesting to develop bootstrap or other resampling-based methods for constructing confidence bands for both the persistence intensity and density functions.

%% file: basics.tex
\section{Notation}
We use boldface small letters like $\bm{u,x,\omega}$ to denote points in $\mathbb{R}^2$ and sub-scripted letters like $x_1,x_2$ to denote their entries. Boldface capital letters like $\bm{X},\bm{Y}$ would be used to denote points on a Riemann manifold. For any positive integer $n$, the symbol $[n]$ refers to the set of all positive integers no larger than $n$. For any set $S$, the symbol $2^S$ represents the power set of $S$, which contains all subsets of $S$ as its elements. The set of all non-negative real numbers would be denoted as $\mathbb{R}_{\geq 0}$. For any function $f$ with domain $\mathcal{A}$, the infinity norm of $f$ is denoted as $\|f\|_{\infty}\coloneqq\sup_{x \in \mathcal{A}} |f(x)|$. 

\section{Background: The persistence diagram}\label{app:basics}

In this section, we give a brief introduction to the persistence diagram. We refer readers to \cite{chazal2018density} for a detailed description. Consider a random point cloud $\bm{X} = (\bm{X}_1,\bm{X}_2,\ldots,\bm{X}_N) \in \mathcal{M}^N$ where $\mathcal{M}$ is a Riemann manifold; and a \emph{filtering function} $\varphi:2^{[N]} \times \mathcal{M}^N \to \mathbb{R}$, which satisfies
\begin{align*}
\varphi(J,\bm{X}) \leq \varphi(J',\bm{X}),\quad \forall J \subset J' \in 2^{[N]}, \bm{X} \in \mathcal{M}^N.
\end{align*}
A simplicial complex given $\bm{X}$ and $\varphi$ at level $\alpha$ is defined as 
\begin{align*}
K_{\alpha}(\bm{X},\varphi) = \{J \subset 2^{[N]} \mid \varphi(J,\bm{X}) \leq \alpha\}.
\end{align*}
Two common examples are the \emph{Cech complex}, where $\varphi(J,\bm{X})$ equals the radius of the circumscribed ball of $\bm{X}[J]$;  and the \emph{Vietoris-Rips complex}, where $\varphi[J,\bm{X}]$ is chosen as the maximum distance between points in $\bm{X}[J]$. 

Throughout the paper, we assume that the filtering function $\varphi$ takes its value in $[0,L]$. For all values $\alpha \in [0,L]$, the sequence of simplicial complexes $\{K_{\alpha}(\bm{X},\varphi)\}_{\alpha \in [0,L]}$ forms a \emph{filtration} denoted as $\mathcal{F}(\bm{X},\varphi)$, where $K_{\alpha}(\bm{X},\varphi) \subseteq K_{\alpha'}(\bm{X},\varphi)$ whenever $\alpha \leq \alpha'$.

\emph{Persistent homology} is a method for computing topological features of a simplicial complex, and can be represented by the \emph{persistence diagram}. In the filtration $\mathcal{F}(\bm{X},\varphi)$, for any persistent homology that begins to appear at level $b$ and disappears at level $d$, we say that the homology is \emph{born} at $b$ and \emph{dies} at $d$. With $\Omega$ defined as in \eqref{eq:defn-Omega-1}, the persistence diagram of the point cloud $\bm{X}$ is a multiset on $\Omega$ that summarizes the birth and death times of all persistent homologies in the filtration $\mathcal{F}(\bm{X},\varphi)$:
\begin{align*}
\mathsf{Dgm}(\bm{X},\varphi) &= \{(b_i,d_i): \text{the } i\text{-th persistent homology in } \mathcal{F}(\bm{X},\varphi) \\
&\text{ that is born at } b_i \text{ and dies at } d_i\}.
\end{align*}

%% file: support.tex
\section{Supportive theoretical results}

\subsection{Validation of Assumption~\ref{as:bound}}\label{sec:validation}
In this part, we provide some common data-generating mechanisms where Assumption~\ref{as:bound} can be validated. 

\begin{theorem}\label{thm:bound-intensity}
Let $q,d$ be two positive integers with $q>d$. Let $\kappa$ be a density on $[0,1]^d$ such that $0 < \inf \kappa \leq \sup \kappa < \infty$. Suppose that $\bm{X}_N$ be either a binomial process with parameters $N$ and $\kappa$ or a Poisson process of intensity $N \kappa$ in the cube $[0,1]^d$. Denote $p(\bm{u})$ as the intensity function for the $k$-dimensional expected persistent measure induced by the Vietoris-Rips filtration. Then when $N$ is sufficiently large, for $\bm{u} \in \Omega$, there exists a polynomial function $\mathsf{poly}(\cdot)$, such that
\begin{equation}
p(\bm{u})\leq\mathsf{poly}(N,d)\sup\kappa,\label{eq:bound-intensity}
\end{equation}
and $\overline{p}(\bm{u})$ can be correspondingly bounded.
\end{theorem}

\begin{theorem}\label{thm:bound-density}
Let $q,d$ be two positive integers with $q>d$. Let $\kappa$ be a density on $[0,1]^{d \times N}$ such that $0 < \inf \kappa < \sup \kappa < \infty$. Suppose that $\bm{X}_1,\bm{X}_2,\ldots,\bm{X}_N \in [0,1]^d$ and that $\bm{X} = (\bm{X}_1,\bm{X}_2,\ldots, \bm{X}_N) \sim \kappa$. Denote $\tilde{p}(\bm{u})$ as the persistence density induced by the Vietoris-Rips filtration of $\bm{X}$. Then there exists a polynomial function $\mathsf(\cdot)$, such that
\begin{equation}
\tilde{p}(\bm{u})\leq\mathsf{poly}(N,d)\sup\kappa.\label{eq:bound-density}
\end{equation}
\end{theorem}

\begin{remark}
The bound in \eqref{eq:bound-intensity} seems to mismatch with Assumption~\ref{as:bound},
since the assumption is on $\bar{p}(\bm{u})=\|\bm{u}-\partial\Omega\|_{2}^{q}p(\bm{u})$
while \eqref{eq:bound-intensity} provides the polynomial bound on
$p(\bm{u})$ directly. So one can imagine that there is no benefit
on considering the assumption on $\bar{p}$ instead of $p$. However,
though not in the formal proof, we believe that $\bar{p}(\bm{u})$
would have a polynomial bound with a slower growth order with respect
to the sample size $N$. This is since as the sample size $N$ grows,
the corresponding persistence diagram tends to have more points that
are close to the diagonal line, as observed in \citep{divol2019choice}.
Hence the term $\|\bm{u}-\partial\Omega\|_{2}^{q}$ can suppress that
effect and $\bar{p}$ can be bounded by a function with a slower growh
order with respect to $N$.
\end{remark}

\begin{remark}
At first glance, comparing \eqref{eq:bound-intensity} and \eqref{eq:bound-density}
seems to indicate that the persistence intensity function $p$ and
the persistence density function $\tilde{p}$ have more or less similar
asymptotic properties with respect to the sample size $N$. However,
this is mainly due to that the polynomial bound $\mathsf{poly}(N,d)$
is comprehensive; for example, $N$ and $N^{10}$ are both in $\mathsf{poly}(N,d)$,
hence the bound $\mathsf{poly}(N,d)$ only guarantees that the function
does not blow up too fast such as an exponential function. And we
believe that, in fact, the growth order is different for $p$
and $\tilde{p}$. This is mainly due to that when the sample size
$N$ is large, the persistence diagram induced by $\bm{X}$ tends
to have more points $N(D)$ in the persistence diagram, and hence
the normalized measure $\tilde{\mu}(A)=\frac{1}{N(D)}\sum_{i=1}^{N(D)}\delta_{\bm{r}_{i}}(A)$
and the corresponding persistence density function $\tilde{p}$ benefits
from the term $\frac{1}{N(D)}$ lowering the growth order with respect
to $N$. In fact, in the proofs of Theorem~\ref{thm:bound-intensity}
and \ref{thm:bound-density} in Section~\ref{app:proof-bound},
the bounds $\mathsf{poly}(N,d)$ are $N^{5}d^{3}$ for \eqref{eq:bound-intensity}
and $N^{4}d^{3}$ for \eqref{eq:bound-density}, although the bounds
need not necessarily equal to the actual asymptotic orders of $p$
and $\tilde{p}$.
\end{remark}

\subsection{Clarification of Assumptions}
\label{sec:assumptions}

In this part, we provide the details in the smoothness assumption of the persistence intensity and density functions, and the regularization assumptions of the kernel function.

\paragraph{H\"older smoothness.} Recall from Assumption~\ref{as:density} that we assume the persistence intensity function $p(\cdot)$ and the persistence density function $\tilde{p}(\cdot)$ are H\"older smooth. A function $f: \Omega \to \mathbb{R}_{\geq 0}$ is H\"older smooth with parameter $L_f$ of oreder $s>0$ if it is $\lfloor s \rfloor$-times continuously differentiable and that for any $\bm{x},\bm{x}' \in \Omega$,
\begin{equation}
\left|f(\bm{x}') - f(\bm{x}) - \sum_{t=1}^{\lfloor s \rfloor} \frac{1}{t!}\sum_{t_1 + t_2 = t,t_{1},t_{2}\geq 0} \frac{\mathrm{d}^t f(\bm{x})}{\mathrm{d}x_1^{t_1}\mathrm{d}x_2^{t_2}} (x_1'-x_1)^{t_1} (x_2'-x_2)^{t_2}\right| \leq L_f \|\bm{x}'-\bm{x}\|_2^s.
\label{eq:density-holder}
\end{equation}

\paragraph{Assumptions regarding the kernel function.} Throughout the paper, we assume the kernel function $K(\cdot)$ satisfies some properties that are commonly used in non-parametric statistics \cite{gine2021mathematical}. Specifically, we make the following assumption.

\begin{assumption}\label{assumption:K}
The kernel function $K:\mathbb{R}^2 \to \mathbb{R}$ satisfies the following conditions:
\begin{enumerate}[label=(\alph*)]
    \item $K(\bm{x}) = 0$ for all $\bm{x}$ with $\|\bm{x}\|_2 > 1$;
    \item $\|K\|_\infty \coloneqq \sup_{\bm{x}} |K(\bm{x})| < \infty$;
    \item $\int_{\mathbb{R}^2} K(\bm{x}) \mathrm{d}\bm{x}=1$;
    \item $\|K\|_2^2 \coloneqq \int_{\mathbb{R}^2} K^2(\bm{x}) \mathrm{d}\bm{x} < \infty$.
    \item There exists a positive integer $s$, such that for all non-negative integers $s_1,s_2$ satisfying $1 \leq s_1 + s_2 < s$, 
    \begin{align*}\int_{\bm{x} \in \mathbb{R}^2} x_1^{s_1}x_2^{s_2}K(\bm{x}) \mathrm{d}\bm{x} = 0.
    \end{align*}
    \item $K$ is $L_K$-Lipchitz with respect to the $\ell_2$ norm on $\mathbb{R}^2$.
\end{enumerate}
\end{assumption}

\subsection{Upper bound for the bias of kernel-based estimators}

In this section, we specify the upper bounds on the  bias of our kernel-based estimators for the persistence intensity/density function, the linear functionals of the persistent measure, and the persistence betti number in Theorems \ref{thm:intensity-bias}, \ref{thm:F-bias} and Corollary \ref{thm:pers-betti-kernel-bias} respectively.

\begin{theorem}\label{thm:intensity-bias}
Under Assumption~\ref{as:density}, there exist  constants $L_{\tilde{p}}, L_p > 0$ such that, for any $\bm{\omega} \in \Omega$,
\begin{align*}
&|\mathbb{E}[\hat{p}_h(\bm{\omega})] - p(\bm{\omega})| \leq L_p h^s\int_{\|\bm{v}\|_2 \leq 1}|K(\bm{v})| \|\bm{v}\|_2^s \mathrm{d}\bm{v}, \quad \text{and} \\ 
&|\mathbb{E}[\check{p}_h(\bm{\omega})] - \tilde{p}(\bm{\omega})| \leq L_{\tilde{p}} h^s\int_{\|\bm{v}\|_2 \leq 1}|K(\bm{v})| \|\bm{v}\|_2^s \mathrm{d}\bm{v}.
\end{align*}
\end{theorem}

The following theorems provide uniform bounds on the bias and variation of these kernel-based estimators.
\begin{theorem}\label{thm:F-bias}
Under Assumption~\ref{as:density},  there exist  constants $L_{\tilde{p}}, L_p > 0$ such that 
\begin{align*}
&\sup_{\Psi \in \mathscr{F}_{2h,R}} \left|\mathbb{E}[\hat{\Psi}_{h}] - \Psi \right|\leq L_{p}\left(\frac{L}{\sqrt{2}}\right)^{q}h^{s} \int_{\|\bm{v}\|_2 \leq 1} |K(\bm{v})|\|\bm{v}\|_2^2 \mathrm{d} \bm{v}; \quad \text{and} \\ 
&\sup_{\Psi \in \widetilde{\mathscr{F}}_{R}} \left|\mathbb{E}[\check{\Psi}_{h}] - \widetilde{\Psi} \right| \leq L_{\tilde{p}} h^s R \int_{\|\bm{v}\|_2 \leq 1} |K(\bm{v})|\|\bm{v}\|_2^2 \mathrm{d} \bm{v}.
\end{align*}
\end{theorem}


It is also straightforward to see that "the biases of both $\hat{\beta}_{\bm{x},h}$ and $\check{\beta}_{\bm{x},h}$  are of order $h^s$", uniformly in $\bm{x}$; see Corollary \ref{thm:pers-betti-kernel-bias} below.


\begin{corollary}\label{thm:pers-betti-kernel-bias}
Under Assumption~\ref{as:density}, it holds that
\begin{align*}
&\sup_{\bm{x} \in \Omega} \left|\mathbb{E}[\hat{\beta}_{\bm{x},h}] - \beta_{\bm{x}}\right| \leq L_p h^s \frac{L^2}{4}\int_{\|\bm{v}\|_2 \leq 1} K(\bm{v})\|\bm{v}\|_2^2 \mathrm{d} \bm{v}, \quad \text{and} \\ 
&\sup_{\bm{x} \in \Omega} \left|\mathbb{E}[\check{\beta}_{\bm{x},h}] - \tilde{\beta}_{\bm{x}}\right| \leq L_{\tilde{p}} h^s \frac{L^2}{4}\int_{\|\bm{v}\|_2 \leq 1} K(\bm{v})\|\bm{v}\|_2^2 \mathrm{d} \bm{v}
\end{align*}
\end{corollary}


\subsection{Minimax lower bound for estimating the persistence intensity function}

Below we provide a minimax lower bound on the $L_\infty$ estimation error of the persistence intensity function by levering well-known minimax arguments for estimating a smooth probability density function based on an i.i.d. sample; see \cite{gine2021mathematical} for details, as well for the definition of Besov norms.

\begin{theorem}\label{thm:weight-intensity-minimax}
Let $\mathscr{F}$ denote the set of functions on $\Omega$ with  Besov norm bounded by $B>0$:
\begin{align*}
\mathscr{F} = \{f: \Omega \to \mathbb{R}, \|f\|_{B_{\infty,\infty}^s} \leq B\}.
\end{align*}
Then,
\begin{align*}
\inf_{\hat{p}_n} \sup_{P} \mathbb{E}_{\mu_1,\ldots,\mu_n \overset{\text{i.i.d.}}{\sim} P} \sup_{\bm{\omega} \in \Omega}  
 \|\bm{\omega} - \partial \Omega\|_2^q |\hat{p}_n(\bm{\omega}) - p(\bm{\omega})| \geq O(n^{-\frac{s}{2(s+1)}}),
\end{align*}
where the infimum is taken over estimator $\hat{p}_n$  mapping $\mu_1,\ldots,\mu_n$ to an intensity function in $\mathscr{F}$, the supremum is over the set of all probability distribution on $\mathcal{Z}_{L,M}^q$ and $p$ is the intensity function of $\mathbb{E}_P[\mu]$.


\end{theorem}
\subsection{Estimating the persistence surface}\label{appendix:persistence.surface}
For estimating the persistence surface in \eqref{eq:defn-pers-surface}, we directly generate the persistence surface from the empirical averaged persistence measure $\bar{\mu}_n$ given by
\begin{align*}
A \in \mathcal{B} \mapsto \bar{\mu}_n(A) = \frac{1}{n}\sum_{i=1}^n \mu_i(A).
\end{align*}
Since $\bar{\mu}_n$ is unbiased for $\mathbb{E}[\mu]$ and $\rho$ is a linear transformation, $\rho_h(\bar{\mu}_n)$ is also unbiased for $\rho_h(\mathbb{E}[\mu])$. The following theorem bounds its variation.
\begin{theorem}\label{thm:surface-var}
With the choice of the weight function
\begin{align*}
f(\bm{\omega}) = \|\bm{\omega} - \partial \Omega\|_2^q,
\end{align*}
when Assumptions~\ref{as:bound}(a) and \ref{as:total-persistence} hold true, there exists a constant $C$ depending on $L,M,L_K, \|K\|_{\infty}$ and $\|\bar{p}\|_{\infty}$, such that for any $\delta \in (0,1)$, it can be guaranteed with probability at least $1-\delta$ that
\begin{align*}
\|\rho_h(\bar{\mu}_n) - \rho_h(\mathbb{E}[\mu])\|_{\infty} \leq C \max\left\{ \frac{1}{nh^2} \log \frac{1}{\delta h^2}, \sqrt{\frac{1}{nh^2}} \sqrt{\log \frac{1}{\delta h^2}}\right\}.
\end{align*}
\end{theorem}


\subsection{Estimating the persistent betti number by the empirical averaged persistence measure}
As an alternative to the kernel-based estimator for the persistent betti number in \eqref{eq:kernel-est-pers-betti}, we can directly use the empirical persistent betti number as the estimator:
\begin{align*}
\bar{\beta}_{\bm{x}} = \bar{\mu}_n(B_{\bm{x}}).
\end{align*}

Since $\bar{\mu}_n$ is an unbiased estimator for $\mathbb{E}[\mu]$, $\bar{\beta}_{\bm{x}}$ is an unbiased estimator for $\beta_{\bm{x}}$. As for the variation of the estimator, we provide the following theorem.
\begin{theorem}\label{thm:pers-betti-var}
Under Assumptions~\ref{as:density}, \ref{as:bound}(a) and \ref{as:total-persistence}, for any $\delta \in (0,1)$, there exists a universal constant $C$ such that with probability at least $1-\delta$ , it can be guaranteed that
\begin{align*}
& \sup_{\bm{x} \in \Omega_{\ell}} |\bar{\beta}_{\bm{x}} - \beta_{\bm{x}}| \leq C \Bigg( \frac{M\ell^{-q}}{n}\left(2\log (M\ell^{-q}n+1) + \log \frac{1}{\delta} \right)\\
&\qquad\qquad+\sqrt{\min\left\{ \frac{M^{2}\ell^{-2q}}{n},\frac{\sqrt{2}ML\ell^{1-2q}\left\Vert \bar{p}\right\Vert _{\infty}}{(q-1)_{+}n}\right\}} \left(\sqrt{2\log (M\ell^{-q}n+1)} + \sqrt{\log \frac{1}{\delta}} \right)\Bigg),
\end{align*}
where $(q-1)_{+}=\max\{q-1,0\}$.
\end{theorem}

%% file: preliminary.tex
\section{Preliminary facts}

In this section we present and prove various auxiliary results that are needed in the proofs of the main theorems.

\subsection{Preliminary facts for the proof of Theorem~\ref{thm:bound-intensity}}
Bounding the weighted intensity function as in Theorem~\ref{thm:bound-intensity} requires a detailed exploration of the persistent diagram for the Vietoris-Rips filtration. Throughout this section, we will consider the filtering function corresponding to the Vietoris-Rips filtration
\begin{align*}
\varphi[J](\bm{X}) = \min_{i,j \in J,i \neq j} \|\bm{X}_i - \bm{X}_j\|_2.
\end{align*}

Firstly, we state a form of the \textbf{area formula}  given by \citep{morgan2016geometric}, which would be useful for a change of variable in deriving the intensity function for the expected persistence measure.

\begin{theorem}\label{thm:area-formula}
Denote $\mathscr{L}^M$ as the $M$-dimensional Lebesgue measure and $\mathscr{H}^M$ as the $M$-dimensional Hausdorff measure. Consider a Lipchitz function $f:\mathbb{R}^M \to \mathbb{R}^N$ for $M \leq N$. If $h:\mathbb{R}^M \to \mathbb{R}$ is an $\mathscr{L}^M$-integrable function, then
\begin{align*}
\int_{\mathbb{R}^M} h(\bm{X})J_{\bm{X}} f(\bm{X}) \mathrm{d} \mathscr{L}^M(\bm{X}) = \int_{\mathbb{R}^N} \sum_{\bm{X} \in f^{-1}\{\bm{Y}\}} h(\bm{X}) \mathrm{d}\mathscr{H}^M \bm{Y},
\end{align*} 
where $J_{\bm{X}}f(\bm{X})$ is the Jacobian determinant of the function $f$:
\begin{align*}
J_{\bm{X}}f(\bm{X}) = \sqrt{\mathsf{det} \left(
\left(\frac{\mathrm{d}f}{\mathrm{d}\bm{X}}\right)^\top \left(\frac{\mathrm{d}f}{\mathrm{d}\bm{X}}\right)\right)}.
\end{align*}
\end{theorem}

Theorem~\ref{thm:area-formula} directly implies the following corollary, the proof of which would be omitted.

\begin{corollary}\label{cor:area-formula}
Let $\psi: \mathbb{R}^M \to \mathbb{R}^N$ be a Lipchitz bijection with $M \leq N$, and $\kappa: \mathbb{R}^N \to \mathbb{R}$ be a function which satisfies that $h \coloneqq \kappa \circ \psi$ is $\mathscr{L}^M$-integrable. Then
\[
\int_{\mathbb{R}^M} \kappa \circ \psi (\bm{X}) J_{\bm{X}} \psi(\bm{X}) \mathrm{d} \mathscr{L}^{M} (\bm{X}) = \int_{\mathbb{R}^N} \kappa(\bm{Y}) \mathrm{d} \mathscr{H}^{M} (\bm{Y}).
\]
\end{corollary}

The following proposition considers two kinds of partitions of the unit cube $[0,1]^{d \times N}$,  with each part satisfying some desired properties. 
\begin{proposition}\label{prop:partition-W}
There exists a set $S$ with cardinality $\text{card}(S) = 4d^2$, such that for any $J_1,J_2 \subset [N]$ that satisfies $J_1 \neq J_2, |J_1| = |J_2| = 2$,  bearing a zero-measured set, $[0,1]^{d \times n}$ can be partitioned as
\[
[0,1]^{d \times n} = \bigcup_{s \in S} W_{J_1,J_2}^s,
\]
such that within each part $W_{J_1,J_2}^s$, there exists a diffeomorphism $\Psi_{J_1,J_2}^s:W_{J_1,J_2}^s \to  \mathbb{R}^2 \times [0,1]^{nd-2}$, such that:
\begin{enumerate}
\item For every $\bm{X} \in W_{J_1,J_2}^s$, $\Psi_{J_1,J_2}^s(\bm{X})_1 = \varphi[J_1](\bm{X})$ and $\Psi_{J_1,J_2}^s(\bm{X})_2 = \varphi[J_2](\bm{X})$;
\item The Jacobian determinant $J_{\bm{X}} \Psi_{J_1,J_2}^s(\bm{X}) \geq \frac{1}{d}$.
\end{enumerate}
\end{proposition}
\emph{Proof}: Let $S = [d]^2 \times \{-1,+1\}^2$, then it is easy to see that $|S| = 4d^2$. For any $J_1,J_2 \subset [n]$ with $J_1 \neq J_2$ and $|J_1| = |J_2| = 2$, let denote $J_1 = \{i_1,j_1\}$, $J_2 = \{i_2,j_2\}$ with $j_2  = \max\{j \in J_2: j \notin J_1\}$. For any $s = (k_1,k_2,s_1,s_2) \in S$, let
\begin{align*}
W_{J_1,J_2}^s = \{X:&\{k_1\} = \text{argmax}_k |X_{i_1}^k - X_{j_1}^k|, s_1(X_{j_1}^k-X_{i_1}^k) > 0,\\
&\{k_2\} = \text{argmax}_k |X_{i_2}^k - X_{j_2}^k|, s_2(X_{j_2}^k-X_{i_2}^k) > 0.\}
\end{align*}

Notice here that $\{k_1\} = \text{argmax}_k |X_{i_1}^k - X_{j_1}^k|$ means $k_1$ is the \emph{only} index for $|X_{i_1}^k - X_{j_1}^k|$ to reach its maximum. 

We begin by proving that $\{W_{J_1,J_2}^s\}_{s \in S}$ forms a partition of $[0,1]^{d \times n}$ bearing a zero-measured set. Firstly, for $s, s' \in S$ with $s \neq s'$, it is easy to see that  $W_{J_1,J_2}^s$ and $W_{J_1,J_2}^{s'}$ are disjoint. Secondly, if 
\[
\bm{X} \in [0,1]^{d \times n} - \bigcup_{s \in S} W_{J_1,J_2}^s,
\]
then by definition, there exists $k,k' \in [d]$, such that $k \neq k'$ and that either
\[
|X_{j_1}^k - X_{i_1}^k| = |X_{j_1}^{k'} - X_{i_1}^{k'}|
\]
or 
\[
|X_{j_2}^k - X_{i_2}^k| = |X_{j_2}^{k'} - X_{i_2}^{k'}|.
\]

Notice that for any $k,k' \in [d]$ with $k \neq k'$, the set
\begin{align*}
&\left\{\bm{X}:|X_{j_1}^k - X_{i_1}^k| = |X_{j_1}^{k'} - X_{i_1}^{k'}|\right\}\\
&= \left\{\bm{X}:X_{j_1}^k - X_{i_1}^k = |X_{j_1}^{k'} - X_{i_1}^{k'}|\right\} \cup \left\{\bm{X}:X_{j_1}^k - X_{i_1}^k = -|X_{j_1}^{k'} - X_{i_1}^{k'}|\right\},
\end{align*}
where the sets 
\begin{align*}
&\left\{\bm{X} \in [0,1]^{d \times n}:X_{j_1}^k - X_{i_1}^k = |X_{j_1}^{k'} - X_{i_1}^{k'}|\right\} \quad \text{and  }\\
&\left\{\bm{X} \in [0,1]^{d \times n}:X_{j_1}^k - X_{i_1}^k = -|X_{j_1}^{k'} - X_{i_1}^{k'}|\right\} 
\end{align*}
 are a subsets of $(nd-1)$ dimensional linear manifolds in $[0,1]^{d \times n}$, and are therefore zero-measured in $\mathscr{L}^{nd}$. Similarly, we can prove that the set $[0,1]^{d \times n} - \bigcup_{s \in S} W_{J_1,J_2}^s$ is the union of a finite number of subsets of $(nd-1)$ dimensional linear manifolds in $[0,1]^{d \times n}$. Consequently, 
\[
\bigcup_{s \in S} W_{J_1,J_2}^s
\]
is a partition of $[0,1]^{ d \times n}$ bearing a zero-measured set.

Furthermore, define $\Psi_{J_1,J_2}^s$ as
\[
\Psi_{J_1,J_2}^s(\bm{X}) = \left(\varphi[J_1](\bm{X}),\varphi[J_2](\bm{X}),\{X_j^k\}_{\substack{1 \leq j \leq n \\ 1 \leq k \leq d\\ (j,k) \neq (j_1,k_1) \\ (j,k) \neq (j_2,k_2)}}\right), \quad \forall \bm{X} \in W_{J_1,J_2}^s.
\]
Then we can firstly notice that
\begin{align*}
&X_{j_1}^{k_1} = s_1 \sqrt{u_1^2-\sum_{k \neq k_1} \left(X_{j_1}^k\right)^2} + X_{i_1}^{k_1} \quad \text{and}\\
&X_{j_2}^{k_2} = s_2 \sqrt{u_2^2-\sum_{k \neq k_2} \left(X_{j_2}^k\right)^2} + X_{i_2}^{k_2},
\end{align*}
for $u_1 = \varphi[J_1](X)$ and $u_2 = \varphi[J_2](X)$. This validates $\Psi_{J_1,J_2}^s$ as a diffeomorphism. The proof now boils down to bounding the Jacobian of $\Psi_{J_1,J_2}^s$. Towards this end, notice that the partial derivative of $\varphi$ is bounded by
\begin{align*}
\left|\frac{\partial \varphi[J_1](\bm{X})}{\partial X_{j_1}^{k_1}}\right| &= \left|\frac{\partial }{\partial X_{j_1}^{k_1}} \sqrt{\sum_{k=1}^d (X_{i_1}^k-X_{j_1}^k)^2}\right|\\
&= \left|\frac{X_{j_1}^{k_1} - X_{i_1}^{k_1}}{\sqrt{\sum_{k=1}^d (X_{i_1}^k-X_{j_1}^k)^2}}\right|\\
&\geq \frac{1}{\sqrt{d}},
\end{align*}
where in the last line we applied the fact that
\[
\left|X_{j_1}^{k_1} - X_{i_1}^{k_1}\right| = \max_{1 \leq k \leq d} \left|X_{j_1}^{k} - X_{i_1}^{k}\right| \geq \sqrt{\frac{1}{d} \sum_{k=1}^d (X_{i_1}^k-X_{j_1}^k)^2}.
\]
Similarly,
\begin{align*}
\left|\frac{\partial \varphi[J_2](\bm{X})}{\partial X_{j_2}^{k_2}}\right| &= \left|\frac{\partial }{\partial X_{j_2}^{k_2}} \sqrt{\sum_{k=1}^d (X_{i_2}^k-X_{j_2}^k)^2}\right| \geq \frac{1}{\sqrt{d}}.
\end{align*}

Furthermore, since $j_2 \notin J_1$, it is easy to see that
\[
\frac{\partial \varphi[J_1](\bm{X})}{\partial X_{j_2}^{k_2}} = 0.
\]
Therefore, the Jacobian determinant of $\Psi_{J_1,J_2}^s$ is bounded by
\begin{align*}
J_{\bm{X}}\Psi_{J_1,J_2}^s(\bm{X}) &= \left|\text{det}\left(\frac{\mathrm{d}\Psi_{J_1,J_2}^s(\bm{X})}{\mathrm{d}\bm{X}}\right)\right|\\
&= \left|\text{det}\left(\begin{pmatrix}
\mathbf{I}_{nd-2} & \mathbf{0}_{(nd-2)\times 1} & \mathbf{0}_{(nd-2)\times 1}\\
\mathbf{0}_{1 \times (nd-2)} & \frac{\partial \varphi[J_1](\bm{X})}{\partial X_{j_1}^{k_1}} & \frac{\partial \varphi[J_1](\bm{X})}{\partial X_{j_2}^{k_2}}\\
\mathbf{0}_{1 \times (nd-2)} & \frac{\partial \varphi[J_2](\bm{X})}{\partial X_{j_1}^{k_1}} & \frac{\partial \varphi[J_2](\bm{X})}{\partial X_{j_2}^{k_2}}
\end{pmatrix}\right)\right|\\
&= \left|\frac{\partial \varphi[J_1](\bm{X})}{\partial X_{j_1}^{k_1}} \cdot \frac{\partial \varphi[J_2](\bm{X})}{\partial X_{j_2}^{k_2}}\right| \geq \frac{1}{d}.
\end{align*}

This completes the proof.
\qed

The following is important for representing of the persistence intensity function $p$ and the persistence density function $\tilde{p}$.

\begin{proposition}\label{prop:partition-V}
Bearing a zero-measured set, $[0,1]^{d \times n}$ can be partitioned as
\[
[0,1]^{d \times n} = \bigcup_{r=1}^{R} V_r,
\]
such that

\begin{enumerate}
\item For every $\bm{X}, \bm{X}' \in V_r$, $J_1,J_2 \subset [n]$ with $|J_1| = |J_2| = 2$, it is guaranteed that $\varphi[J_1](\bm{X}) \neq \varphi[J_2](\bm{X})$; furthermore, if $\varphi[J_1](\bm{X}) < \varphi[J_2](\bm{X})$, then $\varphi[J_1](\bm{X}') < \varphi[J_2](\bm{X}')$;
\item For every $\bm{X}, \bm{X} \in V_r$, $J_1,J_2,J_3,J_4 \subset [n]$ with $|J_1| = |J_2| = |J_3| = |J_4|=2$, it is guaranteed that $\varphi[J_1](\bm{X}) - \varphi[J_2](\bm{X}) \neq \varphi[J_3](\bm{X}) - \varphi[J_4](\bm{X})$; furthermore, if $\varphi[J_1](\bm{X}) - \varphi[J_2](\bm{X}) > \varphi[J_3](\bm{X}) - \varphi[J_4](\bm{X}) > 0$, then $\varphi[J_1](\bm{X}
') - \varphi[J_2](\bm{X}') > \varphi[J_3](\bm{X}') - \varphi[J_4](\bm{X}) > 0$.
\item For every  $r \in [R]$ and $\bm{X} \in V_r$, there are $N_r$ points in $\mathsf{Dgm}(\bm{X},\varphi)$; furthermore, all these points can be ordered by their orthogonal distance to the diagonal, and the order is fixed for all $\bm{X} \in V_r$.
\end{enumerate}

Furthermore, the expected persistence measure $\mathbb{E}[\mu]$ and its normalized counterpart $\mathbb{E}[\tilde{\mu}]$ can be characterized such that for any Borel set $B \subset \Omega$, 
\begin{align*}
&\mathbb{E}[\mu](B) = \sum_{r=1}^R \sum_{i=1}^{N_r}\int_{x \in \Phi^{-1}[J_{ir}^1,J_{ir}^2](B) \cap V_r} \kappa(\bm{X}) \mathrm{d} \bm{X} \quad \text{and} \\ 
&\mathbb{E}[\tilde{\mu}](B) = \sum_{r=1}^R \frac{1}{N_r}\sum_{i=1}^{N_r}\int_{x \in \Phi^{-1}[J_{ir}^1,J_{ir}^2](B) \cap V_r} \kappa(\bm{X}) \mathrm{d} \bm{X}
\end{align*},
in which
\[
\Phi[J_1,J_2](\bm{X}) = (\varphi[J_1](\bm{X}),\varphi[J_2](\bm{X})),
\]
and $J_{ir}^1,J_{ir}^2$ are the simplicial complexes corresponding to the birth and death of the $i$-th persistence homology for all $\bm{X} \in V_r$.

\end{proposition}

\emph{Proof:} For simplicity, we only give a sketch of the proof for this proposition. A weaker version of this proposition is proved in \citep{chazal2018density}, where the second property of the partition is not required. Therefore, the partition we aim to construct here is a refinement of the partition given in \citep{chazal2018density}. In order to see that the second condition can be reached, we firstly prove that the set
\begin{align*}
A = \big\{\bm{X} \in [0,1]^{d \times n}: &\exists J_1,J_2,J_3,J_4 \subset [n], \text{ s.t.}\\
&|J_1| = |J_2| = |J_3| = |J_4| = 2,\\
& J_1 \neq J_2, J_3 \neq J_4, (J_1,J_2) \neq (J_3,J_4),\\
& \varphi[J_1](\bm{X}) -\varphi[J_2](\bm{X}) = \varphi[J_3](\bm{X}) - \varphi[J_4](\bm{X}) \big\}
\end{align*}

is zero-measured. For this step, the technique in proving Lemma 4.1 in \citep{chazal2018density} can be applied to prove that $A$ does not contain any open set, and all its points are singular.

We can further define
\[
\mathcal{F}_n^2 = \{(J_1,J_2):J_1,J_2 \subset [n],|J_1| = |J_2| = 2, J_1 \neq J_2\}.
\]

Since $A$ is zero-measured, we can only consider the set $[0,1]^{d \times n} - A$, on which
\[
\left\{\Delta \varphi[J_1,J_2](\bm{X})\coloneqq\varphi[J_1](\bm{X}) - \varphi[J_2](\bm{X}) \right\}_{(J_1,J_2) \in \mathcal{F}_n^2}
\]

must take different values for different $(J_1,J_2) \in \mathcal{F}_n^2$. Denote these values as $r_1 < r_2 < ... < r_L$, and let $E_{\ell}(\bm{X})$ denote the element $(J_1,J_2 ) \subset \mathcal{F}_n^2$ such that $\Delta \varphi[J_1,J_2](
\bm{X}) = r_\ell$. The sets $E_1(\bm{X}),E_2(\bm{X}),...,E_L(\bm{X})$ then form a partition of $\mathcal{F}_n^2$. With similar techniques as Lemma 4.2 in \citep{chazal2018density}, we can prove that the map $\bm{X} \mapsto \mathcal{A}^2(\bm{X})$ is locally constant almost surely everywhere. This essentially completes the proof.

\qed


The following lemma is a direct application of Proposition~4.6 in \cite{divol2019choice}, and guarantees that the number of points in the persistence diagram $\mathsf{Dgm}(\bm{X},\varphi)$ that are far enough from the diagonal is upper bounded in terms of the expectation.
\begin{lemma}
\label{lem:persistence-diagram-number-bound}
Let $\kappa$ be a probability density function on $[0,1]^d$ that satisfies $0 < \inf \kappa < \sup \kappa < \infty$. Denote $\mathbb{X}_n$ as a binomial process with parameters $n$ and $\kappa$ or a Poisson process with parameter $n\kappa$ on $[0,1]^d$. In the $k$th dimensional persistence diagram of the Vietoris-Rips filtration of $\mathbb{X}_n$, let $N_{\ell}$ be the number of points with persistence of at least $\ell$. Then there are some universal constant $C$ that the expectation of $N_{\ell}$ is upper bounded as
\[
\mathbb{E}\left[N_{\ell}\right]\leq Cn\exp\left(-Cn\ell^{d}\right),
\]
where $C$ is a constant depends only on $k$.
\end{lemma}

\emph{Proof:} 
Let $\mu$ be the persistence measure corresponding to the $k$-th
dimensional persistence diagram of the Vietoris-Rips filtration of
$\mathbb{X}_{n}$. From Proposition~4.6 in \cite{divol2019choice},

\[
P\left(\mu(\mathbb{R}\times[\ell,\infty))>t\right)\leq c_{1}\exp\left(-c_{2}\left(n\ell^{d}+(\frac{t}{n})^{1/(k+1)}\right)\right).
\]
And hence the expectation of $\mu(\mathbb{R}\times[\ell,\infty))$
is bounded as 
\begin{align*}
\mathbb{E}\left[\mu(\mathbb{R}\times[\ell,\infty))\right] & \leq\int_{0}^{\infty}c_{1}\exp\left(-c_{2}\left(n\ell^{d}+(\frac{t}{n})^{1/(k+1)}\right)\right)dt\\
 & =c_{1}\exp\left(-c_{2}(n\ell^{d})\right)\int_{0}^{\infty}\exp\left(-c_{2}(\frac{t}{n})^{1/(k+1)}\right)dt\\
 & =c_{1}\exp\left(-c_{2}(n\ell^{d})\right)\int_{0}^{\infty}(k+1)nu^{k}\exp\left(-c_{2}u\right)du\\
 & =Cn\exp\left(-Cn\ell^{d}\right),
\end{align*}
for some constant $C$ that depends on $k$. Now, $\mathbb{R}\times[\ell,\infty)$
contains all the homological features whose persistence is at least
$\ell$, so 
\[
N_{\ell}\leq\mu(\mathbb{R}\times[\ell,\infty)).
\]
And hence 
\[
\mathbb{E}\left[N_{\ell}\right]\leq Cn\exp\left(-Cn\ell^{d}\right).
\]
\qed

\subsection{Uniform tail bounds}

In this section, we provide some uniform tail bound theorems that are important for bounding the variation of estimators. We will omit the proofs of these theorems in the paper.

\paragraph{The Talagrand's inequality.} The following form of the Talagrand's inequality was shown in \citep{steinwart2008support}. 
\begin{theorem}\label{thm:Talagrand} 
Let $(\mathcal{Z},\mathscr{F},P)$ be a probability space and $(T,d)$ be a separable metric space.  Consider a function class $\mathcal{G} = \{g_t:t \in T\} \in L_0(\mathcal{Z})$, such that the function $t \mapsto g_t(z)$ is continuous in $t$ for all $z \in \mathcal{Z}$.  Furthermore, suppose that there exists a constant $B>0,\sigma^2 > 0$ such that for all $g \in \mathcal{G}$, $\mathbb{E}[g] = 0,\mathbb{E}[g^2] \leq \sigma^2,||g||_\infty \leq B$. Let $Z_1,Z_2,...,Z_n \sim \text{ i.i.d. } P$, and define
\[
G = \sup_{g \in \mathcal{G}} \left|\frac{1}{n} \sum_{i=1}^n g(Z_i)\right| .
\]
Then for any $\delta \in (0,1)$, with probability of at least $1-\delta$,
\begin{equation}
G \leq 4 \mathbb{E}[G] + \sqrt{\frac{2\sigma^2}{n} \log \frac{1}{\delta}} +  \frac{B}{n} \log \frac{1}{\delta}.
\end{equation}
\end{theorem}

Theorem~\ref{thm:Talagrand} implies that the expectation of $G$ is an important factor in bounding $G$. The following theorem gives and upper bound of $\mathbb{E}[G]$ by the covering number of $\mathcal{G}$.

\begin{theorem}\label{thm:expectation}
Under the same conditions as in Theorem~\ref{thm:Talagrand}, if for any $\eta \in (0,B)$, there exists $A>0$, $\nu>0$ such that for any probability measure $Q$ on $\mathcal{Z}$, the covering number
\[
\mathscr{N}(\mathcal{G},L_2(Q),\eta) \leq \left(\frac{AB}{\eta}\right)^\nu,
\]
then there exists a constant $C$ such that
\[
\mathbb{E}[G] \leq C\left(\frac{\nu B}{n} \log \left(\frac{AB}{\sigma}\right) + \sqrt{\frac{\nu \sigma^2}{n} \log\left(\frac{AB}{\sigma}\right) }  \right).
\]
\end{theorem}

\paragraph{Tail bound by polynomial discrimination.}
As an alternative to the Talagrand's inequality, the following theorem bounds $G$ with high probability when the function class $\mathcal{G}$ has \emph{polynomial discrimination}. The proof applies the Bernstein's inequality and a straightforward union bound argument.
\begin{theorem}\label{thm:discrimination}
Under the same conditions as in Theorem~\ref{thm:Talagrand}, define 
\begin{align}
\mathcal{G}(\bm{Z}_1^n) = \{(g(Z_1),g(Z_2),...,g(Z_n)): g \in \mathcal{G}\}.
\end{align}
If the cardinality of the set $\mathcal{G}(\bm{Z}_1^n)$ is bounded by
\begin{align}
\text{Card}(\mathcal{G}(\bm{Z}_1^n)) \leq (An+1)^\nu
\end{align}
for some $\nu>0$, then there exists a universal constant $C$ such that with probability at least $1-\delta$,
\begin{align}
G \leq C\left(\sqrt{\frac{\sigma^2}{n}}\left(\sqrt{\nu \log (An+1)} + \sqrt{\log \frac{1}{\delta}}\right) + \frac{B}{n}\left(\nu \log (An+1) + \log \frac{1}{\delta}\right)\right)
\end{align}
\end{theorem}

The following lemma shows that for persistent measures with bounded total persistence, the total mass of the set away from the diagonal $\partial \Omega$ is upper bounded.

\begin{lemma}\label{lemma:q}
Let $\Omega_{\ell}$ denote the set of points in $\Omega$ that are at least $\ell$ away from the diagonal:
\[
\Omega_{\ell} = \{\bm{\omega} \in \Omega: \|\bm{\omega} - \partial \Omega\|_2 \geq \ell\}.
\]
Then for a persistent measure $\mu$, if $\text{Pers}_q(\mu) \leq M$, then $\mu(\Omega_{\ell}) \leq M\ell^{-q}$.
\end{lemma}

The following theorem shown in \cite{divol2021estimation} provides a standard lower bound for the minimax rate of estimating a probability density function using independent samples. This is useful for deducting the minimax rate for estimating the (weighted) intensity functions.

\begin{theorem}\label{thm:pdf-minimax}
Let $\mathscr{F}$ denote the set of probability density functions on $[0,1]^2$ with Bounded Besov norm:
\[
\mathscr{F}= \{f: [0,1]^2 \to \mathbb{R}, \int_{[0,1]^2} f(x) \mathrm{d}x = 1, ||f||_{\infty,\infty}^r \leq B\}. 
\]
Then for any estimator (measurable function) 
\[
\hat{f}_n: ([0,1]^2)^n \to \mathscr{F},
\]
there exists $f \in \mathscr{F}$, such that if $X_1,X_2,...,X_n \sim \text{ i.i.d. } f$, then
\[
\mathbb{E}\|\hat{f}_n(X_1,X_2,...,X_n) - f\|_{\infty} \geq O\left(n^{-\frac{r}{2r+2}}\right).
\]
\end{theorem}

%% file: proofs.tex
\section{Proof of theorems and supportive propositions}

\subsection{Proof of Theorem~\ref{thm:OT-linfty}}\label{sec:D1}
In order to prove Theorem~\ref{thm:OT-linfty}, we firstly show the following supportive lemma.
\begin{lemma}\label{lemma:int-q}
Let $\Omega$ and $\partial \Omega$ be defined as in \eqref{eq:defn-Omega-1} and \eqref{eq:defn-Omega-2}. Then for any $q> 0$,
\begin{align*}
\int_{\Omega} \|\bm{x}-\partial \Omega\|_2^q \mathrm{d}\bm{x} = \frac{2}{(q+1)(q+2)}\left(\frac{L}{\sqrt{2}}\right)^{q+2}.
\end{align*}
\end{lemma}
\emph{Proof of Lemma \ref{lemma:int-q}}: Take the coordinate transformation
\begin{align*}
\begin{cases}
y_1 = \frac{x_2 - x_1}{\sqrt{2}} = \|\bm{x}-\partial \Omega\|_2;\\
y_2 = \frac{x_2 + x_1}{\sqrt{2}}.
\end{cases}
\end{align*}
Then it can be easily verified that the determinant of the Jacobian matrix between $\bm{x}$ and $\bm{y}$ coordinates is 1, and that the $\ell_1$ ball $\Omega$ can be represented using $\bm{y}$ coordinates by
\[
\Omega = \left\{(y_1,y_2): 0 < y_1 \leq \frac{L}{\sqrt{2}},  y_1 \leq y_2 \leq \sqrt{2}L - y_1 \right\}.
\]
Therefore,
\begin{align*}
\int_{\Omega} \|\bm{x}-\partial \Omega\|_2^q \mathrm{d}\bm{x} &= \int_{0}^{\frac{L}{\sqrt{2}}} \left(\int_{y_1}^{\sqrt{2}L-y_1} \mathrm{d}y_2\right)  y_1^q \mathrm{d}y_1\\
&= \int_{0}^{\frac{L}{\sqrt{2}}} (\sqrt{2}L-2y_1) y_1^q \mathrm{d}y_1\\
&= \frac{2}{(q+1)(q+2)}\left(\frac{L}{\sqrt{2}}\right)^{q+2}.
\end{align*}

With this lemma, we can now prove Theorem~\ref{thm:OT-linfty}.

\emph{Proof of Theorem~\ref{thm:OT-linfty}}: The main idea of bounding the OT distance is to construct an admissible transport between $\mu$ and $\nu$, and then control the cost of this transport. We will separate the proof into three steps accordingly.

\paragraph{Step 1: Construct an admissible transport from $\mu$ to $\nu$.} Define $\hat{\pi}$ as a measure on $\overline{\Omega} \times \overline{\Omega}$ such that for any Borel sets $A,B \subset \overline{\Omega}$,

\begin{equation}\label{eq:define-pi-hat}
\begin{split}
\hat{\pi}(A \times B) &= \int_{A\cap B \cap \Omega} \min \{p_{\mu}(\bm{x}),p_{\nu}(\bm{x})\} \mathrm{d}\bm{x} + \\
& \int_{A \cap \mathsf{Proj}_{\partial \Omega}^{-1}(B) \cap \Omega} [p_{\mu}(\bm{x}) - p_{\nu}(\bm{x})]^+ \mathrm{d}\bm{x} + \int_{B \cap \mathsf{Proj}_{\partial \Omega}^{-1}(A) \cap \Omega} [p_{\nu}(\bm{x}) - p_{\mu}(\bm{x})]^+ \mathrm{d}\bm{x}.
\end{split}
\end{equation}
Here, for any set $A \subset \overline{\Omega}$,
\begin{align*}
\mathsf{Proj}_{\partial \Omega}^{-1}(A) = \{\bm{\omega} \in \Omega: \mathsf{Proj}_{\partial \Omega}(\omega) \in A\}.
\end{align*}
Intuitively, $\hat{\pi}$ represents such a transport: at each point $\bm{x} \in \Omega$, if $p_{\mu}(\bm{x}) > p_{\nu}(\bm{x})$, then we transport the mass of $p_{\nu}$ from $\bm{x}$ to $\bm{x}$, and the remaining mass from $\bm{x}$ to its projection onto $\partial \Omega$; if $p_{\nu}(\bm{x}) > p_{\mu}(\bm{x})$, then the opposite is done.

Firstly, we prove that this is an admissible transport between $\mu$ and $\nu$. Notice that for any Borel set $A \subset \Omega$, $A \cap \overline{\Omega} \cap \Omega = A$, $A \cap \mathsf{Proj}_{\partial \Omega}^{-1}(\overline{\Omega}) \cap \Omega = A$ and $\mathsf{Proj}_{\partial \Omega}^{-1}(A) = \emptyset$. Therefore, by taking $B  = \overline{\Omega}$ in \eqref{eq:define-pi-hat}, we get
\begin{align*}
\hat{\pi}(A \times \overline{\Omega}) &= \int_{A} \min \{p_{\mu}(\bm{x}),p_{\nu}(\bm{x})\} \mathrm{d}\bm{x} + \int_{A}  [p_{\mu}(\bm{x}) - p_{\nu}(\bm{x})]^+ \mathrm{d}\bm{x} + 0 \\
&= \int_{A} \left\{\min \{p_{\mu}(\bm{x}),p_{\nu}(\bm{x})\} + [p_{\mu}(\bm{x}) - p_{\nu}(\bm{x})]^+\right\}\mathrm{d}\bm{x}  \\
&= \int_{A} p_{\mu}(\bm{x}) \mathrm{d}\bm{x} = \mu(A).
\end{align*}
Similarly, we can prove that $\hat{\pi}(\overline{\Omega} \times B) = \nu(B)$ for any Borel set $B \subset \Omega$. Therefore, $\hat{\pi}$ is an admissible transport between $\mu$ and $\nu$. 

\paragraph{Step 2: Present $\mathrm{d} \hat{\pi}$.} In order to calculate the transport cost of $\hat{\pi}$, we firstly need to present $\mathrm{d} \hat{\pi}$. For this, we would make use of \emph{pushforward measures}. Define $\imath:\bar{\Omega}\to\bar{\Omega}\times\bar{\Omega}$ by
$\imath(\bm{x})=(\bm{x},\bm{x})$, and let $\jmath:\bar{\Omega}\times\bar{\Omega}\to\bar{\Omega}$
be satisfying $\jmath\circ\imath=id$. Furthermore, let
$\imath_{*}(\lambda_{\Omega})$ be the pushforward measure on $\bar{\Omega}\times\bar{\Omega}$ generated by $\imath$.  Then for any Borel sets ${A,B} \subset \overline{\Omega}$, one has $\imath^{-1}(A\times B)=A\cap B$,
and the first term in \eqref{eq:define-pi-hat} can be presented as
\begin{align*}
 & \int_{A\cap B\cap\Omega}\min\left\{ p_{\mu}(\bm{x}),p_{\nu}(\bm{x})\right\} d\bm{x}\\
 & =\int_{\imath^{-1}(A\times B)}\min\left\{ (p_{\mu}\circ\jmath)(\imath(\bm{x})),(p_{\nu}\circ\jmath)(\imath(\bm{x}))\right\} d\lambda_{\Omega}(\bm{x})\\
 & =\int_{A\times B}\min\left\{ (p_{\mu}\circ\jmath)(\bm{x},\bm{y}),(p_{\nu}\circ\jmath)(\bm{x},\bm{y})\right\} d\imath_{*}(\lambda_{\Omega})(\bm{x},\bm{y}).
\end{align*}

For the second term in \eqref{eq:define-pi-hat}, we can similarly, define $\imath^{(1)}:\bar{\Omega}\to\bar{\Omega}\times\bar{\Omega}$
by $\imath^{(1)}(\bm{x})=(\bm{x},{\rm Proj}_{\partial\Omega}(\bm{x}))$, let $\jmath^{(1)}:\bar{\Omega}\times\bar{\Omega}\to\bar{\Omega}$
be satisfying $\jmath^{(1)}\circ\imath^{(1)}=id$, and consider the
pushforward measure $\imath_{*}^{(1)}(\lambda_{\Omega})$. Then $(\imath^{(1)})^{-1}(A\times B)=A\cap{\rm Proj}_{\partial\Omega}^{-1}(B)$,
and 
\begin{align*}
 & \int_{A\cap{\rm Proj}_{\partial\Omega}^{-1}(B)\cap\Omega}\left[p_{\mu}(\bm{x})-p_{\nu}(\bm{x})\right]^{+}d\bm{x}\\
 & =\int_{(\imath^{(1)})^{-1}(A\times B)}\left[(p_{\mu}\circ\jmath^{(1)})(\imath^{(1)}(\bm{x}))-(p_{\nu}\circ\jmath^{(1)})(\imath^{(1)}(\bm{x}))\right]^{+}d\lambda_{\Omega}(\bm{x})\\
 & =\int_{A\times B}\left[(p_{\mu}\circ\jmath^{(1)})(\bm{x},\bm{y})-(p_{\nu}\circ\jmath^{(1)})(\bm{x},\bm{y})\right]^{+}d\imath_{*}^{(1)}(\lambda_{\Omega})(\bm{x},\bm{y}).
\end{align*}
For the third term in \eqref{eq:define-pi-hat}, we can similarly define $\imath^{(2)}:\bar{\Omega}\to\bar{\Omega}\times\bar{\Omega}$ by $\imath^{(2)}(\bm{x})=({\rm Proj}_{\partial\Omega}(\bm{x}),\bm{x})$, let $\jmath^{(2)}:\bar{\Omega}\times\bar{\Omega}\to\bar{\Omega}$ be satisfying $\jmath^{(2)}\circ\imath^{(2)}=id$, and consider a pushforward measure $\imath_{*}^{(2)}(\lambda_{\Omega})$. Then $(\imath^{(2)})^{-1}(A\times B)={\rm Proj}_{\partial\Omega}^{-1}(A)\cap B$,
and 
\begin{align*}
 & \int_{{\rm Proj}_{\partial\Omega}^{-1}(A)\cap B\cap\Omega}\left[p_{\mu}(\bm{x})-p_{\nu}(\bm{x})\right]^{+}d\bm{x}\\
 & =\int_{(\imath^{(2)})^{-1}(A\times B)}\left[(p_{\mu}\circ\jmath^{(2)})(\imath^{(2)}(\bm{x}))-(p_{\nu}\circ\jmath^{(2)})(\imath^{(2)}(\bm{x}))\right]^{+}d\lambda_{\Omega}(\bm{x})\\
 & =\int_{A\times B}\left[(p_{\mu}\circ\jmath^{(2)})(\bm{x},\bm{y})-(p_{\nu}\circ\jmath^{(2)})(\bm{x},\bm{y})\right]^{+}d\imath_{*}^{(1)}(\lambda_{\Omega})(\bm{x},\bm{y}).
\end{align*}
Combining these results, we can obtain the following presentation of $\mathrm{d}\hat{\pi}$:
\begin{align*}
\mathrm{d}\hat{\pi} & =\min\left\{ (p_{\mu}\circ\jmath)(\bm{x},\bm{y}),(p_{\nu}\circ\jmath)(\bm{x},\bm{y})\right\} d\imath_{*}(\lambda_{\Omega})\\
& \qquad +\left[(p_{\mu}\circ\jmath^{(1)})(\bm{x},\bm{y})-(p_{\nu}\circ\jmath^{(1)})(\bm{x},\bm{y})\right]^{+}d\imath_{*}^{(1)}(\lambda_{\Omega})\\
 & \qquad+\left[(p_{\mu}\circ\jmath^{(2)})(\bm{x},\bm{y})-(p_{\nu}\circ\jmath^{(2)})(\bm{x},\bm{y})\right]^{+}d\imath_{*}^{(2)}(\lambda_{\Omega}).
\end{align*}

\paragraph{Step 3: Calculate the transportation cost of $\hat{\pi}$}. Based on our presentation of $\mathrm{d}{\hat{\pi}}$, the $q$-th order transportation cost of $\hat{\pi}$ is, by definition:

\begin{align}\label{eq:cost-pi-hat}
C_q^q(\hat{\pi}) &= \int_{\overline{\Omega} \times \overline{\Omega}} \|\bm{x}-\bm{y}\|_2^q \mathrm{d}\hat{\pi}(\bm{x},\bm{y}) \nonumber \\
&= \int_{\overline{\Omega} \times \overline{\Omega}} \|\bm{x}-\bm{y}\|_2^q  \min\left\{ (p_{\nu}\circ\jmath)(\bm{x},\bm{y}),(p_{\mu}\circ\jmath)(\bm{x},\bm{y})\right\} \mathrm{d}\imath_{*}(\lambda_{\Omega}) \nonumber \\
&+ \int_{\overline{\Omega} \times \overline{\Omega}} \|\bm{x}-\bm{y}\|_2^q \left[(p_{\mu}\circ\jmath^{(1)})(\bm{x},\bm{y})-(p_{\nu}\circ\jmath^{(1)})(\bm{x},\bm{y})\right]^{+}\mathrm{d}\imath_{*}^{(1)}(\lambda_{\Omega}) \nonumber \\
&+ \int_{\overline{\Omega} \times \overline{\Omega}} \|\bm{x}-\bm{y}\|_2^q \left[(p_{\mu}\circ\jmath^{(2)})(\bm{x},\bm{y})-(p_{\nu}\circ\jmath^{(2)})(\bm{x},\bm{y})\right]^{+}\mathrm{d}\imath_{*}^{(2)}(\lambda_{\Omega}).
\end{align}

We now explore the three terms in \eqref{eq:cost-pi-hat}. First of all, since $\imath_{*}(\lambda_{\Omega})$ is a pushforward measure generated by the function $\imath(\bm{x}) = (\bm{x},\bm{x})$, it is easy to see that
\[
\imath_{*}(\lambda_{\Omega})(\{(\bm{x},\bm{y}) \in \Omega \times \Omega:\bm{x} \neq \bm{y}\}) = 0.
\]
Therefore, the first term in \eqref{eq:cost-pi-hat} is simply
\begin{align*}
&\int_{\overline{\Omega} \times \overline{\Omega}} \|\bm{x}-\bm{y}\|_2^q  \min\left\{ (p_{\nu}\circ\jmath)(\bm{x},\bm{y}),(p_{\mu}\circ\jmath)(\bm{x},\bm{y})\right\} d\imath_{*}(\lambda_{\Omega})\\
&=\int_{(\bm{x},\bm{y}) \in \overline{\Omega} \times \overline{\Omega}, \bm{x} = \bm{y}} \|\bm{x}-\bm{y}\|_2^q  \min\left\{ (p_{\nu}\circ\jmath)(\bm{x},\bm{y}),(p_{\mu}\circ\jmath)(\bm{x},\bm{y})\right\} d\imath_{*}(\lambda_{\Omega}) \\
&= \int_{\bm{x} \in \overline{\Omega}} \|\bm{x}-\bm{x}\|_2^q \min\{p_{\mu}(\bm{x}),p_{\nu}(\bm{x})\} \mathrm{d}\bm{x}= 0.
\end{align*}
As for the second term, notice that $\imath^{(1)}_{*}(\lambda_{\Omega})$ is a pushforward measure generated by the function $\imath^{(1)}(\bm{x}) = (\bm{x},\mathsf{Proj}_{\partial \Omega}(\bm{x}))$. Therefore by definition,
\[
\imath^{(1)}_{*}(\lambda_{\Omega})(\{(\bm{x},\bm{y}) \in \Omega \times \Omega: \bm{y} \neq \mathsf{Proj}_{\partial \Omega}(\bm{x})\}) = 0.
\]

Hence, the second term in \eqref{eq:cost-pi-hat} is equal to
\begin{align*}
&\int_{\overline{\Omega} \times \overline{\Omega}} \|\bm{x}-\bm{y}\|_2^q \left[(p_{\mu}\circ\jmath^{(1)})(\bm{x},\bm{y})-(p_{\nu}\circ\jmath^{(1)})(\bm{x},\bm{y})\right]^{+}\mathrm{d}\imath_{*}^{(1)}(\lambda_{\Omega})\\
&=\int_{(\bm{x},\bm{y}) \in \overline{\Omega} \times \overline{\Omega}, \bm{y} = \mathsf{Proj}_{\partial \Omega}(\bm{x})} \|\bm{x}-\bm{y}\|_2^q \\
&\qquad \times \left[(p_{\mu}\circ\jmath^{(1)})(\bm{x},\mathsf{Proj}_{\partial \Omega}(\bm{x}))-(p_{\nu}\circ\jmath^{(1)})(\bm{x},\mathsf{Proj}_{\partial \Omega}(\bm{x}))\right]^{+}\mathrm{d}\imath_{*}^{(1)}(\lambda_{\Omega})\\
&= \int_{\bm{x} \in \overline{\Omega}} \|\bm{x}-\mathsf{Proj}_{\partial \Omega}(\bm{x})\|_2^q \left[(p_{\mu}\circ\jmath^{(1)}\circ \imath^{(1)})(\bm{x}) - (p_{\nu}\circ\jmath^{(1)}\circ \imath^{(1)})(\bm{x})\right] \mathrm{d}\bm{x}\\
&= \int_{\Omega}  \|\bm{x}-\partial \Omega\|_2^q \left[p_{\mu}(\bm{x}) - p_{\nu}(\bm{x})\right]^+ \mathrm{d}\bm{x}.
\end{align*}

Similarly, we can obtain that the third term of \eqref{eq:cost-pi-hat} is equal to
\begin{align*}
&\int_{\overline{\Omega} \times \overline{\Omega}} \|\bm{x}-\bm{y}\|_2^q \left[(p_{\mu}\circ\jmath^{(2)})(\bm{x},\bm{y})-(p_{\nu}\circ\jmath^{(2)})(\bm{x},\bm{y})\right]^{+}\mathrm{d}\imath_{*}^{(2)}(\lambda_{\Omega})\\
&= \int_{\Omega} [p_{\nu}(\bm{x}) - p_{\mu}(\bm{x})]^+ \|\bm{x}-\partial \Omega\|_2^q \mathrm{d}\bm{x}.
\end{align*}

Combining these results, we obtain
\begin{align*}
C_q^q(\hat{\pi})& = \int_{\Omega} [p_{\mu}(\bm{x}) - p_{\nu}(\bm{x})]^+ \|\bm{x}-\partial \Omega\|_2^q \mathrm{d}\bm{x} + \int_{\Omega} [p_{\nu}(\bm{x}) - p_{\mu}(\bm{x})]^+ \|\bm{x}-\partial \Omega\|_2^q \mathrm{d}\bm{x}\\
&= \int_{\Omega}  |p_{\mu}(\bm{x})-p_{\nu}(\bm{x})| \|\bm{x}-\partial \Omega\|_2^q \mathrm{d}\bm{x}\\
&\leq \|p_{\mu} - p_{\nu}\|_{\infty} \int_{\Omega}  \|\bm{x}-\partial \Omega\|_2^q \mathrm{d}\bm{x} = \frac{2}{(q+1)(q+2)}\left(\frac{L}{\sqrt{2}}\right)^{q+2} \|p_{\mu} - p_{\nu}\|_{\infty} .
\end{align*}
Notice that the last equality uses Lemma \ref{lemma:int-q}. 

Finally, since $\hat{\pi}$ is an admissible transport from $\mu$ to $\nu$, the optimal transport distance between $\mu$ and $\nu$, $\mathsf{OT}_q(\mu,\nu)$, should be at most $C_q(\hat{\pi})$. The bound in Theorem~\ref{thm:OT-linfty} follows naturally.


\paragraph{Example of converging OT distance while intensity functions diverge.} Consider the following sequences of intensity functions
\begin{align*}
&p_{\mu_n} = \frac{4^n}{L^2} \mathds{1}\left\{\|\bm{x}-\bm{u}_n\|_1 < \frac{\sqrt{2}L}{2^{n+1}}\right\}\\
&p_{\nu_n} = \frac{4^n}{L^2} \mathds{1}\left\{\|\bm{x}-\bm{d}_n\|_1 < \frac{\sqrt{2}L}{2^{n+1}}\right\},
\end{align*}
in which 
\begin{align*}
&\bm{u}_n = \left(\frac{\sqrt{2}L}{4},\frac{\sqrt{2}L}{4} + \frac{\sqrt{2}L}{2^{n+1}}\right)\\
&\bm{d}_n = \left(\frac{\sqrt{2}L}{4} - \frac{\sqrt{2}L}{2^{n+1}},\frac{\sqrt{2}L}{4}\right).
\end{align*}
Essentially, $\mu_n$ and $\nu_n$ are uniform distributions on two adjacent $\ell_1$ balls. It is easy to verify that the total mass of both $\mu_n$ and $\nu_n$ is 1, and the optimal transport distance between $\mu_n$ and $\nu_n$ is upper bounded by 
\begin{align*}
\mathsf{OT}_q(\mu_n,\nu_n) \leq \frac{L}{2^n} \to 0;
\end{align*}
on the other hand, the $\ell_\infty$ distance between the intensity functions clearly diverges as $n \to \infty$:
\begin{align*}
\|p_{\mu_n} - p_{\nu_n}\|_\infty \geq |p_{\mu_n}(\bm{u}_n)-p_{\nu_n}(\bm{u}_n)| = \frac{4^n}{L^2} \to \infty.
\end{align*}
\qed

\emph{A remark on the bottleneck distance.} We argue that there can be no meaningful upper bound for the bottleneck distance $\mathsf{OT}_{\infty}$ by the $\ell_{\infty}$ distance between the intensity or density functions. Consider the following example:  define $T_h$ as an upper-left triangle in $\Omega$:
\[
T_h \coloneqq \left\{\bm{\omega} \in \Omega \mid \|\bm{\omega} - \partial \Omega\|_2 \geq \frac{L-h}{\sqrt{2}}\right\},
\]
and $T_h'$ as a triangle tangent to the diagonal:
\[
T_{h}'\coloneqq\left\{ \bm{\omega}\in\Omega\mid\left\Vert \bm{\omega}-\left(\frac{L}{2},\frac{L}{2}\right)\right\Vert _{\infty}\leq\frac{h}{2}\right\} .
\]
We define $\mu_h$ as the uniform distribution on $T_h$, so that
\begin{align*}
p_{\mu_h}(\bm{\omega}) = \frac{2}{h^2} \mathds{1}\{\bm{\omega} \in T_h\};
\end{align*}
on the other hand $\nu$ is very similar to $\mu$ but has a small part of its mass on $T_h'$:
\begin{align*}
p_{\nu_h}(\bm{\omega}) = \left(\frac{2}{h^2}-h\right) \mathds{1}\{\bm{\omega} \in T_h\} + h \mathds{1}\{\bm{\omega} \in T_h'\}.
\end{align*}
As $h \to 0$, it is easy to verify that $\|p_{\mu_h} - p_{\nu_h}\|_{\infty} = h\to 0$, while $\mathsf{OT}(\mu_h,\nu_h) \to L/\sqrt{2}.$ This is because although the densities for $\mu$ and $\nu$ becomes very close, there is always a small part of the mass of $\mu$ in $T_h$ that has to be transported to $T_h'$; since the bottleneck distance only considers the \emph{maximum} transport cost, it would converge to the limiting distance between $T_h$ and $T_h'$, which is $L/\sqrt{2}$. It is easy to generalize this example to the case where $p_{\mu_h}$ and $p_{\nu_h}$ are smooth.

\subsection{Proof of Theorem~\ref{thm:intensity-bias}}

Both theorems are classic results on the bias of kernel estimators and are proved by the smoothness of the target functions as supposed by Assumption~\ref{as:density}. We here provides the proof of Theorem~\ref{thm:intensity-bias} (a), and part (b) can be proved in a completely similar fashion.

We firstly clarify the specific smoothness condition proposed by Assumption~\ref{as:density}. It guarantees 
Hence, we can represent the bias of $\mathbb{E}[\hat{p}_h(\bm{\omega})]$ as an integral. Since $\bar{\mu}_n$ is an unbiased estimator for $\mathbb{E}[\mu]$,
\begin{align*}
\mathbb{E}[\hat{p}_h(\bm{\omega})] - p(\bm{\omega}) &= \mathbb{E}\left[\int_{\bm{x}} \frac{1}{h^2} K\left(\frac{\bm{x}-\bm{\omega}}{h}\right) \mathrm{d} \bar{\mu}_n \right] - p(\omega)\\ 
&= \int_{\bm{x}} \frac{1}{h^2} K\left(\frac{\bm{x}-\bm{\omega}}{h}\right) \mathrm{d} \mathbb{E}[\bar{\mu}_n]  - p(\omega)\\ 
&= \int_{\bm{x}} \frac{1}{h^2} K\left(\frac{\bm{x}-\bm{\omega}}{h}\right) p(\bm{x}) \mathrm{d}\bm{x} - p(\omega)\\
&= \int_{\bm{x}} \frac{1}{h^2} K\left(\frac{\bm{x}-\bm{\omega}}{h}\right) [p(\bm{x}) -p(\bm{\omega})]\mathrm{d}\bm{x},
\end{align*}
where in the last line we applied the property that the kernel function $K(\cdot)$ integrals to $1$. We can then apply the smoothness of $p(\cdot)$ as in \eqref{eq:density-holder} and obtain that
\begin{align*}
 &\left|\mathbb{E}[\hat{p}_h(\bm{\omega})] - p(\bm{\omega})\right|\\ 
 &\leq \left|\int_{\bm{x}} \frac{1}{h^2} K\left(\frac{\bm{x}-\bm{\omega}}{h}\right) \sum_{t=1}^{s-1} \frac{1}{t!}\sum_{t_1 + t_2 = t} \frac{\mathrm{d}^t p(\bm{\omega})}{\mathrm{d}\omega_1^{t_1}\mathrm{d}\omega_2^{t_2}} (x_1-\omega_1)^{t_1} (x_2-\omega_2)^{t_2}\mathrm{d}\bm{x}\right|\\
 &+ \int_{\bm{x}} \frac{1}{h^2} \left| K\left(\frac{\bm{x}-\bm{\omega}}{h}\right) \right|L_p \|\bm{x} - \bm{\omega}\|_2^s \mathrm{d}\bm{x}
\end{align*}

By taking a change of variable $\bm{v} = \frac{\bm{x} - \bm{\omega}}{h}$ , the first term can be represented as
\begin{align*}
\sum_{t=1}^{s-1} \frac{1}{t!}\sum_{t_1 + t_2 = t} \frac{\mathrm{d}^t p(\bm{\omega})}{\mathrm{d}\omega_1^{t_1}\mathrm{d}\omega_2^{t_2}} \int_{\|\bm{v}\|_2 \leq 1} K(\bm{v}) h^t v_1^{t_1}v_2^{t_2} \mathrm{d}\bm{v}.
\end{align*}
The zero-moment condition of the kernel function in Assumption~\ref{assumption:K} guarantees that this term equals to 0. Hence,
\begin{align*}
\left|\mathbb{E}[\hat{p}_h(\bm{\omega})] - p(\bm{\omega})\right| &\leq \int_{\bm{x}} \frac{1}{h^2}\left| K\left(\frac{\bm{x}-\bm{\omega}}{h}\right) \right|L_p \|\bm{x} - \bm{\omega}\|_2^s \mathrm{d} \bm{x} \\ 
&\xlongequal{\bm{v} = (\bm{x} - \bm{\omega})/h} L_p h^s \int_{\|\bm{v}\|_2 \leq 1} |K(\bm{v})| \|\bm{v}\|_2^s \mathrm{d}\bm{v}. 
\end{align*}

\subsection{Proof of Theorem~\ref{thm:intensity-var} (a)}
\paragraph{A useful claim.} The following claim can be applied for easing calculation in Theorem~\ref{thm:intensity-var}.
\begin{claim}
\label{claim:ineq-jensen}
For $q\in\mathbb{R}$ and $x\in[0,1]$, 
\[
1-x^{q}\leq(q\vee1)(1-x),
\]
where $q\vee1=\max\{q,1\}$.
\end{claim}
\paragraph{Proof of Claim~\ref{claim:ineq-jensen}.}
If $q\geq1$ or $q\leq0$, let $f(x)=1-x^{q}$. Then $f'(x)=-qx^{q-1}$
and $f''(x)=-q(q-1)x^{q-2}$, so $f''(x)\leq0$ for $x\in[0,1]$ and
$f$ is concave on $[0,1]$. Then by Jensen's inequality, 
\[
1-x^{q}=f(x)\leq f(1)+f'(1)(x-1)=q(1-x).
\]
If $q\in[0,1]$, then $x^{q}\geq x$ implies 
\[
1-x^{q}\leq1-x.
\]
Hence combining these gives 
\[
1-x^{q}\leq(q\vee1)(1-x).
\]
\qed

This proof applies the Talagrand's inequality. For this purpose, we firstly define an auxiliary family of functions, and then verify the conditions in Theorems~\ref{thm:Talagrand} and \ref{thm:expectation} .

\paragraph{Defining an auxiliary function class.} Let $\mu_1,\mu_2,....,\mu_n$ be i.i.d. random measures in $\mathcal{Z}_{L,M}^q$, $\ell_{\bm{\omega}} = \|\bm{\omega} - \partial \Omega\|_2 - h$ and $g_{\bm{\omega}}$ be defined as
\begin{equation}\label{eq:define-g}
g_{\bm{\omega}}(\mu) = \ell_{\bm{\omega}}^q \left(\int_{\Omega} \frac{1}{h^2} K\left(\frac{\bm{x}-\bm{\omega}}{h}\right)\mathrm{d}\mu - \int_{\Omega} \frac{1}{h^2} K\left(\frac{\bm{x}-\bm{\omega}}{h}\right)\mathrm{d}\mathbb{E}[\mu]\right),
\end{equation}
and $K$ satisfy Assumption~\ref{assumption:K}. Take $\mathcal{Z} = \mathcal{Z}_{L,M}^q$, $(T,d) = (\Omega_{2h},\|\cdot\|_2)$, and for all $\mu \in \mathcal{Z}_{L,M}^q$, define $\mathcal{G} = \{g_{\bm{\omega}}:\bm{\omega} \in \Omega_{2h}\}$. By definition, $g_{\bm{\omega}}(\mu)$ has zero mean and the variation of the kernel estimator $\hat{p}_h(\cdot)$ can be represented by
\begin{align*}
\sup_{\bm{\omega} \in \Omega_{2h}} \ell_{\bm{\omega}}^q |\hat{p}_h(\bm{\omega})-\mathbb{E}[\hat{p}_h(\bm{\omega})]| = \sup_{\bm{\omega} \in \Omega_{2h}} \left|\frac{1}{n}\sum_{i=1}^n g_{\bm{\omega}}(\mu)\right|.
\end{align*}
Hence, in order to apply the Talagrand's inequality, we need to bound $\|g_{\bm{\omega}}(\mu)\|_{\infty}$, $\mathbb{E}[g_{\bm{\omega}}(\mu)^2]$ and the covering number of $\mathcal{G}$. We provide these upper bound accordingly in the following paragraphs.
\paragraph{Bounding$\|g_{\bm{\omega}}(\mu)\|_{\infty}$ and  $\mathbb{E}[g_{\bm{\omega}}(\mu)^2]$. } Notice that since $K$ vanishes outside the unit circle of $\mathbb{R}^2$, for any $\bm{x} \notin \Omega_{\ell_{\bm{\omega}}}$, we have $ \left|\left|\frac{\bm{x}-\bm{\omega}}{h}\right|\right|_2 > 1$ and therefore $K\left(\frac{\bm{x}-\bm{\omega}}{h}\right) = 0$. Hence, for all $\bm{\omega} \in \Omega_{2h}$,
\begin{align}\label{eq:intensity-infty-bound}
|g_{\bm{\omega}}(\mu)| &= \ell_{\bm{\omega}}^q \left|\int_{\Omega} \frac{1}{h^2} K\left(\frac{\bm{x}-\bm{\omega}}{h}\right)\mathrm{d}\mu - \int_{\Omega} \frac{1}{h^2} K\left(\frac{\bm{x}-\bm{\omega}}{h}\right)\mathrm{d}\mathbb{E}[\mu]\right| \nonumber \\
&\leq \ell_{\bm{\omega}}^q \max\left\{\left|\int_{\Omega} \frac{1}{h^2} K\left(\frac{\bm{x}-\bm{\omega}}{h}\right)\mathrm{d}\mu\right| , \left|\int_{\Omega} \frac{1}{h^2} K\left(\frac{\bm{x}-\bm{\omega}}{h}\right)\mathrm{d}\mathbb{E}[\mu]\right| \right\}\nonumber \\
&= \ell_{\bm{\omega}}^q \max \left\{\left|\int_{\Omega_{\ell_{\bm{\omega}}}} \frac{1}{h^2} K\left(\frac{\bm{x}-\bm{\omega}}{h}\right)\mathrm{d}\mu\right| , \left|\int_{\Omega_{\ell_{\bm{\omega}}}} \frac{1}{h^2} K\left(\frac{\bm{x}-\bm{\omega}}{h}\right)\mathrm{d}\mathbb{E}[\mu]\right|\right\}\nonumber \\
&\leq \ell_{\bm{\omega}}^q \frac{\|K\|_{\infty}}{h^2} \max \left\{\left(\mu(\Omega_{\ell_{\bm{\omega}}}),  \mathbb{E}[\mu](\Omega_{\ell_{\bm{\omega}}})\right)\right\}\nonumber \\
&\leq \ell_{\bm{\omega}}^q \frac{\|K\|_{\infty}M}{h^2 \ell_{\bm{\omega}}^q} = \frac{\|K\|_{\infty}M}{h^2}
\end{align}
where in the last inequality we used Lemma \ref{lemma:q}. On the other hand, the variance of $g_{\bm{\omega}}$ is bounded by
\begin{align}\label{eq:intensity-var-bound}
\mathbb{E}[g_{\bm{\omega}}(\mu)^2] &= \ell_{\bm{\omega}}^{2q} \mathbb{E}\left|\int \frac{1}{h^2} K\left(\frac{\bm{x}-\bm{\omega}}{h}\right)\mathrm{d}\mu - \int \frac{1}{h^2} K\left(\frac{\bm{x}-\bm{\omega}}{h}\right)\mathrm{d}\mathbb{E}[\mu]\right|^2 \nonumber \\
&\leq \ell_{\bm{\omega}}^{2q} \mathbb{E}\left|\int_{\Omega_{\ell_{\bm{\omega}}}} \frac{1}{h^2} K\left(\frac{\bm{x}-\bm{\omega}}{h}\right)\mathrm{d}\mu\right|^2 \nonumber \\
&\leq \ell_{\bm{\omega}}^{2q}\mathbb{E} \left\{\mu(\Omega_{\ell_{\bm{\omega}}}) \cdot \int_{\Omega_{\ell}} \frac{1}{h^4} K^2\left(\frac{\bm{x}-\bm{\omega}}{h}\right)\mathrm{d}\mu\right\}\nonumber \\
&= \ell_{\bm{\omega}}^{2q} \mu(\Omega_{\ell}) \int_{\Omega_{\ell_{\bm{\omega}}}} \frac{1}{h^4} K^2\left(\frac{\bm{x}-\bm{\omega}}{h}\right)\mathrm{d}\mathbb{E}[\mu]\\
&\leq \ell_{\bm{\omega}}^{2q} \cdot \frac{M}{\ell_{\bm{\omega}}^q} \int_{\|\bm{x}-\bm{\omega}\|_2 \leq h} \frac{1}{h^4} K^2\left(\frac{\bm{x}-\bm{\omega}}{h}\right) p(\bm{x}) \mathrm{d}\bm{x} \nonumber \\
&\xlongequal{\bm{v} = (\bm{x}-\bm{\omega})/h} \ell_{\bm{\omega}}^{q} M \int_{\|\bm{v}\|_2 \leq 1} \frac{1}{h^2} K^2(\bm{v})p(\bm{\omega}+\bm{v}h) \mathrm{d}\bm{v} \nonumber \\
&\leq  \ell_{\bm{\omega}}^{q} M  \frac{1}{h^2}\frac{\|\bar{p}\|_{\infty}}{\ell_{\bm{\omega}}^q} \int_{\|\bm{v}\|_2 \leq 1}  K^2(\bm{v})\mathrm{d}\bm{v} = \frac{M\|\bar{p}\|_{\infty}\|K\|_2^2}{ h^2}.
\end{align}
\paragraph{Bounding the covering number of $\mathcal{G}$.} For any probability measure $Q$ on $\mathcal{Z}_{L,M}^q$ and any  $\eta \in(0,\frac{\|K\|_{\infty}M}{h^2})$, we aim to bound the covering number of $\mathcal{G}$ with respect to $L_2(Q)$ distance. This requires relating the $L_2(Q)$ distance in $\mathcal{G}$ and the $\ell_2$ distance in $\mathbb{R}^2$. Specifically, for any $\bm{\omega},\bm{\omega}' \in \Omega_{2h}$ and $\mu \in \mathcal{Z}_{L,M}^q$, we can assume without loss of generality that $\ell_{\bm{\omega}} \leq \ell_{\bm{\omega}'}$. In this case, we firstly observe that 
\begin{align}\label{eq:cover-1}
&\left|\ell_{\bm{\omega}}^q \int K\left(\frac{\bm{x}-\bm{\omega}}{h}\right)\mathrm{d} \mu- \ell_{\bm{\omega}'}^q \int K\left(\frac{\bm{x}-\bm{\omega}'}{h}\right)\mathrm{d}\mu \right|\nonumber \\
&\leq \left|\int \ell_{\bm{\omega}}^q \left[K\left(\frac{\bm{x}-\bm{\omega}}{h}\right) - K\left(\frac{\bm{x}-\bm{\omega}'}{h}\right)\right]\mathrm{d}\mu\right| + \left| \int (\ell_{\bm{\omega}}^q - \ell_{\bm{\omega}'}^q) K\left(\frac{\bm{x}-\bm{\omega}'}{h}\right) \mathrm{d}\mu\right|\nonumber \\
&\leq \ell_{\bm{\omega}}^q \int_{\Omega_{\ell_{\bm{\omega}}}}\frac{L_k}{h}\|\bm{\omega} - \bm{\omega}'\|_2  \mathrm{d}\mu + \int_{\Omega_{\ell_{\bm{\omega}'}}} (\ell_{\bm{\omega}'}^q - \ell_{\bm{\omega}}^q) \|K\|_{\infty} \mathrm{d}\mu \nonumber \\
&\leq \ell_{\bm{\omega}}^q \frac{L_k}{h} \|\bm{\omega} - \bm{\omega}'\|_2 \mu(\Omega_{\ell_{\bm{\omega}}}) + \|K\|_{\infty} (\ell_{\bm{\omega}'}^q - \ell_{\bm{\omega}}^q) \mu(\Omega_{\ell_{\bm{\omega}'}})\nonumber \\
&\leq \frac{ML_k}{h}\|\bm{\omega} - \bm{\omega}'\|_2 + M\|K\|_{\infty} \left[1-\left(\frac{\ell_{\bm{\omega}}}{\ell_{\bm{\omega}'}}\right)^q\right].
\end{align}

Since $\ell_{\bm{\omega}} \geq \ell_{\bm{\omega}'} - \|\bm{\omega}- \bm{\omega}'\|_2$, the last term of \eqref{eq:cover-1} can be bounded by using Claim~\ref{claim:ineq-jensen}
and $\ell_{\bm{\omega}}\geq\ell_{\bm{\omega}'}-\|\bm{\omega}-\bm{\omega}'\|_{2}$
as 
\begin{align}
1-\left(\frac{\ell_{\bm{\omega}}}{\ell_{\bm{\omega}}'}\right)^{q} & \leq(q\vee1)\left(1-\frac{\ell_{\bm{\omega}}}{\ell_{\bm{\omega}}'}\right) \nonumber\\
 & \leq\frac{q\vee1}{\ell_{\bm{\omega}}'}\|\bm{\omega}-\bm{\omega}'\|_{2} \nonumber\\
 & \leq\frac{q\vee1}{h}\|\bm{\omega}-\bm{\omega}'\|_{2}. \label{eq:cover-lwdiff-bound}
\end{align}
Notice that in the last line, we applied the fact that since $\bm{\omega}'\in\Omega_{2h}$,
$\ell_{\bm{\omega}'}=\|\bm{\omega}-\partial\Omega\|_{2}-h\geq h$.

From now on, we use $q'$ to denote $q \vee 1$ for simplicity.  Equations \eqref{eq:cover-1} and \eqref{eq:cover-lwdiff-bound} imply that 

\begin{align*}
\left|\ell_{\bm{\omega}}^q \int K\left(\frac{\bm{x}-\bm{\omega}}{h}\right)\mathrm{d} \mu- \ell_{\bm{\omega}'}^q \int K\left(\frac{\bm{x}-\bm{\omega}'}{h}\right)\mathrm{d}\mu \right| \leq \frac{M(L_k + q'\|K\|_{\infty})}{h} \|\bm{\omega} - \bm{\omega}'\|_2.
\end{align*}
Therefore, the difference between $g_{\bm{\omega}}(\mu)$ and $g_{\bm{\omega}'}(\mu)$ can be bounded by 
\begin{align*}
|g_{\bm{\omega}}(\mu)-g_{\bm{\omega}'}(\mu)|&\leq \left|\ell_{\bm{\omega}}^q \int \frac{1}{h^2}K\left(\frac{\bm{x}-\bm{\omega}}{h}\right)\mathrm{d} \mu- \ell_{\bm{\omega}'}^q \int \frac{1}{h^2}K\left(\frac{\bm{x}-\bm{\omega}'}{h}\right)\mathrm{d}\mu \right|\\
&+ \left|\ell_{\bm{\omega}}^q \int \frac{1}{h^2}K\left(\frac{\bm{x}-\bm{\omega}}{h}\right)\mathrm{d} \mathbb{E}[\mu]- \ell_{\bm{\omega}'}^q \int \frac{1}{h^2}K\left(\frac{\bm{x}-\bm{\omega}'}{h}\right)\mathrm{d}\mathbb{E}[\mu] \right|\\
&\leq \frac{2M(L_k + q'\|K\|_{\infty})}{h^3} \|\bm{\omega} - \bm{\omega}'\|_2.
\end{align*}
In this way, we have related the distance between $g_{\bm{\omega}}$ and $g_{\bm{\omega}'}$ to the distance between $\bm{\omega}$ and $\bm{\omega}'$.
Now, for any $\eta \in (0,\frac{\|K\|_{\infty}M}{h^2})$, we can set $\epsilon = \frac{\eta h^3}{2M(L_K+q'\|K\|_{\infty})}$. It is easy to verify that
\begin{align*}
\epsilon < \frac{h^3}{2M(L_K + q'\|K\|_{\infty})} \frac{\|K\|_{\infty}M}{h^2} = \frac{\|K\|_{\infty}}{2(L_K+ q'\|K\|_{\infty})} h < h.
\end{align*}
Hence, we can construct a $\epsilon$-covering of $\Omega_{2h}$ in the $\ell_2$ distance, denoted as $S$. It is easy to show that the covering number
\[
\mathscr{N}(\Omega_{2h},\|\cdot\|_2,\epsilon) \leq \frac{2L^2}{\epsilon^2}.
\]
By definition, for any $\bm{\omega} \in \Omega_{2h}$, there exists $\bm{\omega}' \in S$, such that $\|\bm{\omega} - \bm{\omega}'\|_2 \leq \epsilon < h < \ell_{\bm{\omega}'}$. Therefore, for any measure $Q$ on $\mathcal{Z}_{L,M}^q$,
\begin{align*}
\|g_{\bm{\omega}}(\mu)- g_{\bm{\omega}'}(\mu)\|_{L_2(Q)} &\leq \sup_{\mu \in \mathcal{Z}_{L,M}^q} |g_{\bm{\omega}}(\mu) - g_{\bm{\omega}'}(\mu)|\\
&\leq \frac{2 M (L_K + q'\|K\|_{\infty})}{h^3} \|\bm{\omega} - \bm{\omega}'\|_2 \leq \frac{2 M (L_K + q'\|K\|_{\infty})}{h^3} \epsilon  = \eta.
\end{align*}
In conclusion,
\begin{align}\label{eq:intensity-cover-bound}
\mathscr{N}(\mathcal{G},L_2(Q),\eta)&\leq \mathcal{N}\left(\Omega_{2h},\|\cdot\|_2, \frac{\eta h^3 }{2M(L_K + q'\|K\|_{\infty})}\right) \nonumber \\
&< \left(\frac{4LM(L_K + q'\|K\|_{\infty})}{\eta h^3 }\right)^2.
\end{align}

\paragraph{Completing the proof.} With $\|g_{\bm{\omega}}(\mu)\|_{\infty}$, $\mathbb{E}[g_{\bm{\omega}}(\mu)^2]$ and the covering number of $\mathcal{G}$ bounded as in \eqref{eq:intensity-infty-bound}, \eqref{eq:intensity-var-bound} and \eqref{eq:intensity-cover-bound}, we can apply Theorems~\ref{thm:Talagrand} and \ref{thm:expectation} with
\begin{align*}
\begin{cases}
AB = \frac{4LM(L_K + q'\|K\|_{\infty})}{h^3 };\\
B = \frac{\|K\|_{\infty} M}{h^2 };\\
\sigma^2 = \frac{M\|\bar{p}\|_{\infty}}{ h^2 }\|K\|_2^2;\\
\nu = 2.
\end{cases}
\end{align*}
This gives us the conclusion that with probability at least $1-\delta$, 
\begin{align*}
\sup_{\bm{\omega} \in \Omega_{2h}} \left|\frac{1}{n} \sum_{i=1}^n g_{\bm{\omega}}(\mu) \right| &\lesssim \frac{2\|K\|_{\infty}M}{nh^2 } \log \left(\frac{4L(L_K+q'\|K\|_{\infty})}{\delta h^2 \|K\|_2} \sqrt{\frac{M}{\|\bar{p}\|_{\infty}}}\right) +\\
&\sqrt{\frac{2M\|\bar{p}\|_{\infty}}{n}} \frac{\|K\|_2}{h} \sqrt{\log \left(\frac{4L(L_K+q'\|K\|_{\infty})}{\delta h^2 \|K\|_2} \sqrt{\frac{M}{\|\bar{p}\|_{\infty}}}\right) }.
\end{align*}
\qed

\subsection{Proof of Theorem \ref{thm:intensity-var}(b)}
Part (b) of Theorem \ref{thm:intensity-var} can be proved in a similar, though slightly easier, fashion to part (a). We therefore provide a sketch of the proof and omit the details.
\paragraph{Defining an auxiliary function class.} For every $\tilde{\mu}$ and $\bm{\omega} \in \Omega$, define
\begin{align*}
g_{\bm{\omega}}(\tilde{\mu}) = \int_{\Omega} \frac{1}{h^2} K\left(\frac{\bm{x}-\bm{\omega}}{h}\right)\mathrm{d}\tilde\mu - \int_{\Omega} \frac{1}{h^2} K\left(\frac{\bm{x}-\bm{\omega}}{h}\right)\mathrm{d}\mathbb{E}[\tilde\mu],
\end{align*}
and let $\mathcal{G} = \{g_{\bm{\omega}}: \bm{\omega} \in \Omega\}$. It is easy to verify that $\mathbb{E}[g] \equiv 0$ for all $\bm{\omega} \in \Omega$, and that
\begin{align*}
\|\check{p}_h(\bm\omega) - \tilde{p}(\bm{\omega})\| = \sup_{g \in \mathcal{G}} \left|\frac{1}{n} \sum_{i=1}^n g(\mu_i)\right|.
\end{align*}
\paragraph{Bounding $\|g\|_{\infty}$ and $\mathbb{E}[g^2]$.} Since $\tilde{\mu}$ and $\mathbb{E}[\tilde{\mu}]$ are normalized measures with a total mass of $1$, $\|g\|_{\infty}$ can be bounded by
\begin{align*}
\|g\|_{\infty} \leq \frac{\|K\|_{\infty}}{h^2};
\end{align*}
in the mean time, Assumption \ref{as:bound} (b) guarantees that $\mathbb{E}[g_{\bm{\omega}}(\tilde\mu)^2]$ can be bounded by
\begin{align*}
\mathbb{E}[g_{\bm{\omega}}(\tilde{\mu})^2] \leq \frac{\|\tilde{p}\|_{\infty} \|K\|_2^2}{h^2}.
\end{align*}
\paragraph{Bounding the covering number of $\mathcal{G}$.} We again apply the Lipchitz property of the kernel function $K(\cdot)$ to conclude that for any $\bm{\omega,\omega}' \in \Omega$,
\begin{align*}
|g_{\bm{\omega}}(\tilde{\mu}) - g_{\bm{\omega'}}(\tilde{\mu})| \leq \frac{2L_K}{h^3} \|\bm{\omega- \omega'}\|_2.
\end{align*}
Hence, using a similar reasoning to the proof of part (a), we can bound the covering number of $\mathcal{G}$ by 
\begin{align*}
\mathscr{N}(\mathcal{G},L^2(Q),\eta) < \left(\frac{4LL_K}{\eta h^3}\right)^2.
\end{align*}
\paragraph{Completing the proof.} Theorem \ref{thm:intensity-var} (b) is a direct corollary of Theorems \ref{thm:Talagrand} and \ref{thm:expectation} with the following choice of parameters:
\begin{align*}
\begin{cases}
AB = \frac{4LL_K }{h^3 };\\
B = \frac{\|K\|_{\infty} }{h^2 };\\
\sigma^2 = \frac{\|\tilde{p}\|_{\infty}}{ h^2 }\|K\|_2^2;\\
\nu = 2.
\end{cases}
\end{align*}

\subsection{Proof of Theorems ~\ref{thm:density-minimax} and ~\ref{thm:weight-intensity-minimax}}
In this section, we provide the proof of Theorem \ref{thm:weight-intensity-minimax}, which gives a minimax lower bound for estimating the weighted persistence intensity function. Theorem \ref{thm:density-minimax}, which gives the minimax lower bound for estimating the persistence density function, can be proved in a similar while simpler fashion, so we omit its proof for brevity.

The main idea of this proof is to build a connection of weighted intensity function $\bar{p}(\cdot)$ and a probability density function.  First of all, we can observe the conclusion of Theorem \ref{thm:pdf-minimax} holds true also when the support for the density function is $\Omega$ instead of $[0,1]^2$. Now, notice that for any $\bm{x} \in \Omega$, we can define the following measure:
\begin{align}\label{eq:mu-x}
\mu_{\bm{x}} = M \delta_{\bm{x}} ||\bm{x}-\partial \Omega||_2^{-q}.
\end{align}
It is easy to verify that $\mathsf{Pers}_q(\mu_{\bm{x}}) = M$, so $\mu_{\bm{x}} \in \mathcal{Z}_{L,M}^q$ . Therefore, for any estimator $\hat{p}_n: (\mathcal{Z}_{L,M}^q)^n \to \mathcal{F}$, we can construct the following estimator $\hat{f}_n$:
\begin{align*}
\hat{f}_n(\bm{x}_1,\bm{x}_2,...,\bm{x}_n) = \hat{p}_n(\mu_{\bm{x}_1},\mu_{\bm{x}_2},...,\mu_{\bm{x}_n}).
\end{align*}
Theorem \ref{thm:pdf-minimax} states that there exists a probability density function $f:\Omega \to \mathbb{R}$ with $||f||_{\infty,\infty}^r \leq B$ such that when $X_1,X_2,...,X_n \sim \text{ i.i.d. } f$, 
\[
\mathbb{E} ||\hat{f}_n(X_1,X_2,...,X_n) - f||_\infty \geq O\left(n^{-\frac{r}{2r+2}}\right).
\]
We can apply the probability density function $f$ to construct a probability measure on $\mathcal{Z}_{L,M}^q$.
First, define a map $\Phi:\Omega\to\mathcal{Z}_{L,M}^{q}$ by $\Phi(\bm{x})=\mu_{\bm{x}}$
in \eqref{eq:mu-x}. Impose a measure structure on $\mathcal{Z}_{L,M}^{q}$
by pushforwarding the measure structure on $\Omega$, i.e. $\mathcal{Y}\subset\mathcal{Z}_{L,M}^{q}$
is measurable if and only if $\Phi^{-1}(\mathcal{Y})$ is measurable
in $\Omega$. Define a probability measure $P$ on $\mathcal{Z}_{L,M}^{q}$
as a pushforward measure, i.e., for any measurable set $\mathcal{Y}\subset\mathcal{Z}_{L,M}^{q}$,
\[
P(\mathcal{Y})=\int_{\Phi^{-1}(\mathcal{Y})}f(\bm{x})\mathrm{d}\bm{x}.
\]
Then from the change of variables, 
\[
\int_{\mathcal{Y}}g(\mu)dP(\mu)=\int_{\Phi^{-1}(\mathcal{Y})}g(\Phi(\bm{x}))f(\bm{x})\mathrm{d}\bm{x}.
\]
Now, the intensity for $P$ can be represented as follows:
let $p(\cdot)$ be the intensity function for $\mathbb{E}[\mu]$ when $\mu \sim P$, then for all $u \in \Omega,$
\begin{align}\label{eq:p-f}
\bar{p}(\bm{u}):= \|\bm{u}-\partial \Omega\|_2^q p(\bm{u}) =M f(\bm{u}).
\end{align}
To see this fact, consider any Borel set $\mathcal{A} \subset \Omega$. 
 By definition, the expected measure $\mathbb{E}[\mu]$ satisfies
\begin{align*}
\mathbb{E}[\mu](\mathcal{A}) & =\mathbb{E}[\mu(\mathcal{A})]=\int_{\mathcal{Z}_{L,M}^{q}}\mu(\mathcal{A})dP(\mu)\\
 & =\int_{\Phi^{-1}(\mathcal{Z}_{L,M}^{q})}\Phi(\bm{x})(\mathcal{A})f(\bm{x})\mathrm{d}\bm{x}\\
 & =\int_{\Omega}\mu_{\bm{x}}(\mathcal{A})f(\bm{x})\mathrm{d}\bm{x}\\
 & =\int_{\Omega}M||\bm{x}-\partial\Omega||_{2}^{-q}\mathbf{1}\{\bm{x}\in\mathcal{A}\}f(\bm{x})\mathrm{d}\bm{x}\\
 & =\int_{\mathcal{A}}M||\bm{x}-\partial\Omega||_{2}^{-q}f(\bm{x})\mathrm{d}\bm{x}.
\end{align*}

Since $\mathcal{A}$ can be any Borel set, we get $p(\bm{u}) = M||\bm{u}-\partial \Omega||_2^{-q}$ by definition, and Equation \eqref{eq:p-f} follows naturally. Since the $\ell_{\infty}$ difference between $\hat{f}_n $and $f$ is lower bounded, we can obtain
\begin{align*}
\mathbb{E}_P \sup_{\bm{\omega} \in \Omega} \|\bm{\omega} - \partial \Omega\|_2^q |\hat{p}_n(\bm{\omega})-p(\bm{\omega})| = M \mathbb{E}_f \|\hat{f}_n - f\|_\infty \geq O\left(n^{-\frac{r}{2r+2}}\right).
\end{align*}
\qed

\subsection{Proof of Theorems and Corollaries regarding linear representations of the persistence measure}
The theoretical results regarding linear representations of the persistence measure in Section \ref{sec:linear} are rather direct applications of the theoretical results on estimating the persistence intensity and density functions. We therefore combine their proofs in this section.

\paragraph{Proof of Theorem \ref{thm:F-bias}.}

First, consider the bias of $\hat{\Psi}$. For any $\Psi\in\mathscr{F}_{2h,R}$,
the bias of $\hat{\Psi}$ is upper bounded by 
\begin{align*}
\left|\mathbb{E}[\hat{\Psi}]-\Psi\right| & =\left|\int_{\bm{\omega}\in\Omega_{2h}}f(\bm{\omega})(\mathbb{E}[\hat{p}_{h}(\bm{\omega})]-p(\bm{\omega}))\mathrm{d}\bm{\omega}\right|\\
 & \leq\int_{\bm{\omega}\in\Omega_{2h}}f(\bm{\omega})\left|\mathbb{E}[\hat{p}_{h}(\bm{\omega})]-p(\bm{\omega})\right|\mathrm{d}\bm{\omega}\\
 & \leq\sup_{\bm{\omega}\in\Omega_{2h}}\left|\mathbb{E}[\hat{p}(\bm{\omega})]-p(\bm{\omega})\right|\int_{\bm{\omega}\in\Omega_{2h}}f(\bm{\omega})\mathrm{d}\omega.
\end{align*}
Then under Assumption~\ref{as:density}, Theorem~\ref{thm:intensity-bias}
gives an upper bound as 
\[
\left|\mathbb{E}[\hat{\Psi}]-\Psi\right|\leq L_{p}h^{s}\int_{\bm{\omega}\in\Omega_{2h}}f(\bm{\omega})\mathrm{d}\omega\int_{\|\bm{v}\|_{2}\leq1}|K(\bm{v})|\left\Vert \bm{v}\right\Vert _{2}^{2}\mathrm{d}v.
\]
And then, the definition of $\mathscr{F}_{2h,R}$ implies 
\[
\int_{\bm{\omega}\in\Omega_{2h}}f(\bm{\omega})\mathrm{d}\omega\leq\left(\frac{L}{\sqrt{2}}\right)^{q}\int_{\bm{\omega}\in\Omega_{2h}}\ell_{{\bf \omega}}^{-q}f(\bm{\omega})\mathrm{d}\omega\leq\left(\frac{L}{\sqrt{2}}\right)^{q}R,
\]
Hence it gives a further upper bound for the bias $\left|\mathbb{E}[\hat{\Psi}]-\Psi\right|$ as
\[
\left|\mathbb{E}[\hat{\Psi}]-\Psi\right|\leq L_{p}\left(\frac{L}{\sqrt{2}}\right)^{q}h^{s}R\int_{\|\bm{v}\|_{2}\leq1}|K(\bm{v})|\left\Vert \bm{v}\right\Vert _{2}^{2}\mathrm{d}v.
\]

Second, consider the bias of $\check{\Psi}$. For any $\Psi\in\widetilde{\mathscr{F}}_{R}$,
the bias of $\check{\Psi}$ is upper bounded by

\begin{align*}
\left|\mathbb{E}[\check{\Psi}]-\Psi\right| & =\left|\int_{\bm{\omega}\in\Omega}f(\bm{\omega})(\mathbb{E}[\check{p}_{h}(\bm{\omega})]-\tilde{p}(\bm{\omega}))\mathrm{d}\bm{\omega}\right|\\
 & \leq\int_{\bm{\omega}\in\Omega}f(\bm{\omega})\left|\mathbb{E}[\check{p}_{h}(\bm{\omega})]-\tilde{p}(\bm{\omega})\right|\mathrm{d}\bm{\omega}\\
 & \leq\sup_{\bm{\omega}\in\Omega}\left|\mathbb{E}[\check{p}_{h}(\bm{\omega})]-\tilde{p}(\bm{\omega})\right|\int_{\bm{\omega}\in\Omega}f(\bm{\omega})\mathrm{d}\omega.
\end{align*}
Then under Assumption~\ref{as:density}, Theorem~\ref{thm:intensity-bias}
and the definition of $\widetilde{\mathscr{F}}_{R}$ give the upper
bound as
\[
\left|\mathbb{E}[\check{\Psi}]-\Psi\right|\leq L_{\tilde{p}}h^{s}R\int_{\|\bm{v}\|_{2}\leq1}|K(\bm{v})|\left\Vert \bm{v}\right\Vert _{2}^{2}\mathrm{d}v.
\]

\paragraph{Proof of Theorem \ref{thm:F-var}.} The upper bound for the variation of $\hat{\Psi}$ is a direct corollary of Theorem~\ref{thm:intensity-var} (a) and the fact that

\begin{align*}
\sup_{\Psi\in\mathscr{F}_{2h,R}}\left|\hat{\Psi}-\mathbb{E}[\hat{\Psi}]\right| & =\sup_{\Psi\in\mathscr{F}_{2h,R}}\left|\int_{\bm{\omega}\in\Omega_{2h}}f(\bm{\omega})[\hat{p}_{h}(\bm{\omega})-\mathbb{E}[\hat{p}_{h}](\bm{\omega})]\mathrm{d}\bm{\omega}\right|\\
 & \leq\int_{\bm{\omega}\in\Omega_{2h}}\ell_{\bm{\omega}}^{-q}f(\bm{\omega})\mathrm{d}\bm{\omega}\cdot\sup_{\bm{\omega}\in\Omega_{2h}}\ell_{\bm{\omega}}^{q}\left|\hat{p}_{h}(\bm{\omega})-\mathbb{E}[\hat{p}_{h}](\bm{\omega})\right|\\
 & \leq R\cdot\sup_{\bm{\omega}\in\Omega_{2h}}\ell_{\bm{\omega}}^{q}\left|\hat{p}_{h}(\bm{\omega})-\mathbb{E}[\hat{p}_{h}](\bm{\omega})\right|;
\end{align*}

The upper bound for the variation of $\check{\Psi}$ follows from Theorem~\ref{thm:intensity-var} (b) and a similar relation:
\[
\sup_{\tilde{\Psi} \in \mathscr{F}_R} \left|\check{\Psi} - \mathbb{E}[\check{\Psi}]\right| \leq R \cdot \sup_{\bm{\omega} \in \Omega} |\check{p}_h(\bm{\omega}) - \mathbb{E}[\check{p}_h(\bm{\omega})]|.
\]
\paragraph{Proof of Corollaries \ref{thm:pers-betti-kernel-bias} and \ref{thm:pers-betti-kernel-var}(a).} For every $\bm{x} \in \Omega_{2h}$, define
\begin{align*}
f_{\bm{x}}(\bm{\omega}) = \mathds{1}\left\{\bm{\omega} \in B_{\bm{x}}\right\},  
\end{align*}
and let
\begin{align*}
\mathscr{F}_{2h,R}= \left\{\Psi = \int_{\Omega_{2h}} f_{\bm{x}}(\bm{\omega}) \mathrm{d} \mathbb{E}[\mu] \bigg | \bm{x} \in \Omega_{2h} \right\}.
\end{align*}
Corollary~\ref{thm:pers-betti-kernel-bias} follows from Theorem \ref{thm:F-bias} and the fact that
\begin{align*}
\int_{\bm{\omega} \in \Omega_{2h}} f_{\bm{x}}(\bm{\omega}) \mathrm{d} \bm{\omega} \leq \frac{L^2}{4}
\end{align*}
for every $\bm{x} \in \Omega_{2h}$. 
Similarly, Corollary~\ref{thm:pers-betti-kernel-var} follows from Theorem \ref{thm:F-var} and the fact that
\begin{align*}
\int_{\bm{\omega} \in \Omega_{2h}} \ell_{\bm{\omega}}^{-q} f_{\bm{x}}(\bm{\omega}) \mathrm{d} \bm{\omega} \leq C\ell_{\bm{x}}^{2-q},
\end{align*}
for a constant $C$.


\paragraph{Proof of Corollary~\ref{thm:pers-betti-kernel-var}(b).} 
For every $\bm{x} \in \Omega$, we define
\begin{align*}
f_{\bm{x}}(\bm{\omega}) = \mathds{1}\left\{\bm{\omega} \in B_{\bm{x}}\right\}, 
\end{align*}
and let
\begin{align*}
\widetilde{\mathscr{F}}_R = \left\{\widetilde{\Psi} = \int_{\Omega} f_{\bm{x}} (\bm{\omega}) \mathrm{d} \mathbb{E}[\tilde{\mu}] \bigg | \bm{x} \in \Omega \right\}.
\end{align*}
Corollary~\ref{thm:pers-betti-kernel-bias}(b) follows directly from Theorem \ref{thm:F-var} and the fact that for every $\bm{x} \in \Omega$,
\begin{align*}
\int_{\bm{\omega} \in \Omega} f_{\bm{x}}(\bm{\omega}) \mathrm{d} \bm{\omega} \leq \frac{L^2}{4}.
\end{align*}

\subsection{Proof of Theorem \ref{thm:surface-var}}
This proof again involves the Talagrand's inequality, and therefore takes a similar shape to the proof of Theorem \ref{thm:intensity-var}. We begin by defining an auxiliary function class.
\paragraph{Defining the auxiliary function class $\mathcal{G}$.} Recall that we choose the weight function as $f(\bm{\omega}) = \|\bm{\omega} - \partial \Omega\|_2^q$. Therefore, for any persistence measure $\mu \in \mathcal{Z}_{L,M}^q$, its corresponding persistence surface is characterized by
\begin{align*}
\rho_h(\mu)(\bm{u}) = \int_{\Omega} \|\bm{\omega} - \partial \Omega\|_2^q \frac{1}{h^2} K\left(\frac{\bm{u} - \bm{\omega}}{h}\right) \mathrm{d} \mu(\bm{\omega});
\end{align*}
hence, by defining
\begin{align*}
g_{\bm{u}}(\mu) = \int_{\Omega} \|\bm{\omega} - \partial \Omega\|_2^q \frac{1}{h^2} K\left(\frac{\bm{u} - \bm{\omega}}{h}\right) \mathrm{d} \left(\mu - \mathbb{E}[\mu]\right)(\bm{\omega})
\end{align*}
and letting $\mathcal{G} = \{g_{\bm{u}}(\bm{\mu}): \bm{u} \in \Omega\}$, we observe that $\mathbb{E}[g] = 0$ for all $g \in \mathcal{G}$ and 
\begin{align*}
\left\|\rho_h(\bm{\mu}_n) - \mathbb{E}[\rho_h(\bm{\mu})]\right\|_{\infty} = \sup_{g \in \mathcal{G}} \left\|\frac{1}{n} \sum_{i=1}^n g(\mu_i)\right\|.
\end{align*}
\paragraph{Bounding $\|g\|_{\infty}$ and $\mathbb{E}[g^2]$.} Assumptions \ref{as:total-persistence} and \ref{assumption:K} directly implies that for any $g \in \mathcal{G}$ and any $\bm{u} \in \Omega$,
\begin{align*}
\left|g_{\bm{u}}(\mu)\right| &\leq \frac{\|K\|_{\infty}}{h^2} \max \left\{\int_{\Omega} \|\bm{\omega} - \partial \Omega\|_2^q \mathrm{d} \mu, \int_{\Omega} \|\bm{\omega} - \partial \Omega\|_2^q \mathrm{d} \mathbb{E}[\mu]\right\} \\ 
&= \frac{\|K\|_{\infty}}{h^2} \max\left\{\mathsf{Pers}_q(\mu),\mathsf{Pers}_q(\mathbb{E}[\mu])\right\}\leq \frac{M\|K\|_{\infty}}{h^2}. 
\end{align*}
Regarding the variance of $g$, Assumption \ref{as:bound} implies that
\begin{align*}
\mathbb{E}[g_{\bm{u}}(\mu)^{2}] & \leq\|g\|_{\infty}\cdot\int_{\Omega}\|\bm{\omega}-\partial\Omega\|_{2}^{q}\frac{1}{h^{2}}\left|K\left(\frac{\bm{u}-\bm{\omega}}{h}\right)\right|\mathrm{d}\mathbb{E}[\mu]\\
 & \leq\frac{M\|K\|_{\infty}}{h^{2}}\int_{\Omega}\frac{1}{h^{2}}\left|K\left(\frac{\bm{u}-\bm{\omega}}{h}\right)\right|\|\bm{\omega}-\partial\Omega\|_{2}^{q}p(\bm{\omega})\mathrm{d}\bm{\omega}\\
 & \leq\frac{M\|K\|_{\infty}}{h^{2}}\int_{\|\bm{v}\|_{2}\leq1}\left|K(\bm{v})\right|\mathrm{d}\bm{v}\cdot\sup_{\bm{\omega}\in\Omega}\|\bm{\omega}-\partial\Omega\|_{2}^{q}p(\bm{\omega})\\
 & \leq\frac{M\|K\|_{1}\|K\|_{\infty}\|\bar{p}\|_{\infty}}{h^{2}},
\end{align*}
where in the third line we applied the change of variable $\bm{v} = (\bm{u} - \bm{\omega})/h$, and let
\begin{align*}
\|K\|_1 := \int_{\|\bm{v}\|_{2}\leq1}\left|K(\bm{v})\right|\mathrm{d}\bm{v}.
\end{align*}
\paragraph{Covering number of $\mathcal{G}$.} Similar to the proof of Theorem \ref{thm:intensity-var}, we bound the covering number of $\mathcal{G}$ by the Lipchitz property of the kernel function $K$. For any two points $\bm{u,u}' \in \Omega$, Assumption \ref{assumption:K} guarantees that
\begin{align*}
\left|K\left(\frac{\bm{u}-\bm{\omega}}{h}\right) - K\left(\frac{\bm{u}'-\bm{\omega}}{h}\right)\right| \leq \frac{L_K \|\bm{u} - \bm{u}'\|_2}{h}.
\end{align*}
Therefore, it is easy to verify that
\begin{align*}
|g_{\bm{u}}(\mu) - g_{\bm{u}'}(\mu)| \leq \frac{ML_K \|\bm{u}-\bm{u}'\|_2}{h^3}.
\end{align*}
A similar reasoning to the proof of Theorem \ref{thm:intensity-var} yields that the covering number of $\mathcal{G}$ is upper bounded by
\begin{align*}
\mathcal{N}(\mathcal{G},L^2(Q),\eta) \leq \mathcal{N}\left(\Omega, \|\cdot\|_2 , \frac{\eta h^3}{ML_K}\right) \leq 2 \left(\frac{LML_K}{\eta h^3}\right)^2.
\end{align*}
\paragraph{Completing the proof.} Theorem \ref{thm:surface-var} is a direct application of Theorems \ref{thm:Talagrand} and \ref{thm:expectation} with the following choice of parameters:
\begin{align*}
\begin{cases}
&AB = \frac{2LML_k}{h^3};\\ 
&B = \frac{M\|K\|_{\infty}}{h^2}; \\
&\sigma^2 = \frac{M\|K\|_1\|K\|_{\infty} \|\bar{p}\|_{\infty}}{h^2};\\
&\nu = 2.
\end{cases}
\end{align*}

\subsection{Proof of Theorems~\ref{thm:bound-intensity} and \ref{thm:bound-density}}\label{app:proof-bound}
Observe that the persistence diagram of the Vietoris-Rips filtration of $\bm{X} = (\bm{X}_1,\bm{X}_2,...,\bm{X}_N)$ is decided purely by $\{\varphi[J](\bm{X})\}_{J \subset [N],|J| = 2}$, in which 
\begin{align*}
\varphi[J](\bm{X}) = \|\bm{X}_i-\bm{X}_j\|_2,
\end{align*}
for $J = \{i,j\}$. In what follows, we firstly focus on the proof of Theorem~\ref{thm:bound-intensity}, and apply the techniques to that of Theorem~\ref{thm:bound-density} in a similar manner.

\paragraph{Proof of Theorem~\ref{thm:bound-intensity}.} Propositions \ref{prop:partition-V} and \ref{prop:partition-W} imply that for any Borel set $B \subseteq \Omega$,
\begin{align*}
\mathbb{E}[\mu](B) &= \sum_{r=1}^R \sum_{i=1}^{N_r} \sum_{s \in S} \int_{V_r \cap W_{J_{ir}^1,J_{ir}^2}^s \cap \Phi[J_{ir}^1,J_{ir}^2]^{-1}(B)} \kappa(\bm{X}) \mathrm{d}\bm{X}\\
&=\sum_{r=1}^R \sum_{i=1}^{N_r} \sum_{s \in S} \\ &\int_{\Psi_{J_{ir}^1,J_{ir}^2}^s(V_r \cap W_{J_{ir}^1,J_{ir}^2}^s \cap \Phi[J_{ir}^1,J_{ir}^2]^{-1}(B))} \kappa((\Psi_{J_{ir}^1,J_{ir}^2}^s)^{-1}(u,y)) J [\Psi_{J_{ir}^1,J_{ir}^2}^s]^{-1}(\bm{u},\bm{Y})\mathrm{d}\bm{Y} \mathrm{d}u,
\end{align*}
where in the second line we change the variable from $\bm{X} \in [0,1]^{d \times n}$ to $(\bm{Y},\bm{u})$ with $\bm{Y} \in [0,1]^{nd-2}$ and $\bm{u} \in \Omega$. Now, a change of order of summation gives
\begin{align}\label{eq:E-mu-B-1}
\mathbb{E}[\mu](B) & =\sum_{s\in S}\sum_{\substack{J_{1},J_{2}\subset[N]\\
|J_{1}|=|J_{2}|=2\\
J_{1}\neq J_{2}
}
}\sum_{r=1}^{R}\sum_{i=1}^{N_{r}}I(J_{ir}^{1}=J_{1},J_{ir}^{2}=J_{2})\nonumber \\
 & \times\int_{\Psi_{J_{1},J_{2}}^{s}(V_{r}\cap W_{J_{1},J_{2}}^{s}\cap\Phi[J_{1},J_{2}]^{-1}(B))}\kappa((\Psi_{J_{1},J_{2}}^{s})^{-1}(\bm{u},\bm{Y}))J[\Psi_{J_{ir}^{1},J_{ir}^{2}}^{s}]^{-1}(u,y)\mathrm{d}\bm{Y}\mathrm{d}\bm{u}\nonumber \\
 & \leq\sum_{s\in S}\sum_{\substack{J_{1},J_{2}\subset[N]\\
|J_{1}|=|J_{2}|=2\\
J_{1}\neq J_{2}
}
}\sum_{r=1}^{R}\sum_{i=1}^{N_{r}}I(J_{ir}^{1}=J_{1},J_{ir}^{2}=J_{2})\nonumber \\
 & \times\int_{\Psi_{J_{1},J_{2}}^{s}(V_{r}\cap W_{J_{1},J_{2}}^{s}\cap\Phi[J_{1},J_{2}]^{-1}(B))}d\sup\kappa\mathrm{d}\bm{Y}\mathrm{d}\bm{u}\nonumber \\
 & \leq\sum_{s\in S}\sum_{\substack{J_{1},J_{2}\subset[N]\\
|J_{1}|=|J_{2}|=2\\
J_{1}\neq J_{2}
}
}N(B)\int_{\Psi_{J_{1},J_{2}}^{s}(W_{J_{1},J_{2}}^{s}\cap\Phi[J_{1},J_{2}]^{-1}(B))}d\sup\kappa\mathrm{d}\bm{Y}\mathrm{d}\bm{u},
\end{align}
where $N(B)$ is the number of persistent homology points in $B$,
and in the second line we use the facts that $\{V_r\}_{r=1}^R$ are disjoint, $\kappa \leq \sup \kappa$ and $J [\Psi_{J_{ir}^1,J_{ir}^2}^s]^{-1} \leq d$. Hence, bounding $\mathbb{E}[\mu](B)$ boils down to characterizing the domain of integration on the right hand side of \eqref{eq:E-mu-B-1}. For this, notice that by definition, 
\begin{align*}
&(\bm{Y},\bm{u}) \in \Psi_{J_1,J_2}^s(W_{J_1,J_2}^s \cap  \Phi[J_1,J_2]^{-1}(B))\\
&\leftrightarrow \exists \bm{X} \in W_{J_1,J_2}^s, \text{ such that } \Phi[J_1,J_2](\bm{X}) \in B,  \Psi_{J_1,J_2}^s(\bm{X}) = (\bm{Y},\bm{u})\\
&\rightarrow \exists \bm{X} \in W_{J_1,J_2}^s, \text{ such that } \Phi[J_1,J_2](\bm{X}) \in B,  \Phi[J_1,J_2](\bm{X}) = \bm{u}, \text{ and } \bm{Y} \in [0,1]^{Nd-2}\\
&\rightarrow  \bm{u} \in B, \text{ and } \bm{Y} \in [0,1]^{Nd-2}.
\end{align*}
Hence, $\mathbb{E}[\mu](B)$ is upper bounded by
\begin{align*}
\mathbb{E}[\mu](B) & \leq N(B) \sum_{s \in S} \sum_{\substack{J_1,J_2 \subset [N] \\ |J_1| = |J_2| = 2 \\ J_1 \neq J_2}} \int_{\bm{u} \in B, \bm{Y} \in [0,1]^{Nd-2}} d \sup \kappa  \mathrm{d}\bm{Y}\mathrm{d}\bm{u}\\
&=  d \sup \kappa  N(B) \sum_{s \in S} \sum_{\substack{J_1,J_2 \subset [N] \\ |J_1| = |J_2| = 2 \\ J_1 \neq J_2}} \int_{[0,1]^{Nd-2}} \mathrm{d}\bm{Y} \int_{B} \mathrm{d}\bm{u}\\
&= d \sup \kappa N(B) \sum_{s \in S} \sum_{\substack{J_1,J_2 \subset [N] \\ |J_1| = |J_2| = 2 \\ J_1 \neq J_2}} \int_{B} \mathrm{d}\bm{u}.
\end{align*}
This effectively means that the intensity function $p(\bm{u})$ is upper bounded by
\begin{align*}
p(\bm{u}) & \leq\mathbb{E}\left[N(\{\bm{u}\})\right]d\sup\kappa\sum_{s\in S}\sum_{\substack{J_{1},J_{2}\subset[N]\\
|J_{1}|=|J_{2}|=2\\
J_{1}\neq J_{2}
}
}1\\
 & <\mathbb{E}\left[N(\{\bm{u}\})\right]\text{card}(S)|\{(J_{1},J_{2}):|J_{1}|=|J_{2}|=2,J_{1}\neq J_{2},J_{1}\subset[N],J_{2}\subset[N]\}|d\sup\kappa.
\end{align*}
Now, $N(\{\bm{u}\})\leq N_{\ell}$, so Lemma~\ref{lem:persistence-diagram-number-bound}
implies $\mathbb{E}\left[N(\{\bm{u}\})\right]\leq CN$. And $\text{card}(S)\leq4d^{2}$
and $|\{(J_{1},J_{2}):|J_{1}|=|J_{2}|=2,J_{1}\neq J_{2},J_{1}\subset[N],J_{2}\subset[N]\}|\leq\frac{N^{4}}{4}$,
so 
\begin{align*}
p(\bm{u}) & \leq(CN)\cdot(4d^{2})\cdot\left(\frac{N^{4}}{4}\right)\cdot d\sup\kappa\\
 & =C'N^{5}d^{3}\sup\kappa.
\end{align*}
Theorem~\ref{thm:bound-intensity} follows with the choice of 
\begin{align*}
\mathsf{poly}(N,d) = N^{5}d^3.
\end{align*}

\paragraph{Proof of Theorem~\ref{thm:bound-density}.} Propositions \ref{prop:partition-V} and \ref{prop:partition-W} imply that for any Borel set $B \subseteq \Omega$, the normalized persistence measure of $B$ is expressed by
\begin{align*}
\mathbb{E}[\tilde{\mu}](B) &= \sum_{r=1}^R \frac{1}{N_r}\sum_{i=1}^{N_r} \sum_{s \in S} \int_{V_r \cap W_{J_{ir}^1,J_{ir}^2}^s \cap \Phi[J_{ir}^1,J_{ir}^2]^{-1}(B)} \kappa(\bm{X}) \mathrm{d}\bm{X}\\ 
&\leq \sum_{r=1}^R \max_{1 \leq i \leq N_r} \sum_{s \in S} \int_{V_r \cap W_{J_{ir}^1,J_{ir}^2}^s \cap \Phi[J_{ir}^1,J_{ir}^2]^{-1}(B)} \kappa(\bm{X}) \mathrm{d}\bm{X}.
\end{align*}
Hence, same techniques can be applied to show that the persistence density function is upper bounded by
\begin{align*}
\tilde{p}(\bm{u}) & \leq d\sup\kappa\mathbb{E}\left[\frac{N(\{\bm{u}\})}{N(\{\bm{u}\})}\right]\sum_{s\in S}\sum_{\substack{J_{1},J_{2}\subset[N]\\
|J_{1}|=|J_{2}|=2\\
J_{1}\neq J_{2}
}
}1\\
 & \leq d\sup\kappa\max_{1\leq i\leq N({\bm{u}})}\sum_{s\in S}\sum_{\substack{J_{1},J_{2}\subset[N]\\
|J_{1}|=|J_{2}|=2\\
J_{1}\neq J_{2}
}
}1\\
 & \leq\text{card}(S)|\{(J_{1},J_{2}):|J_{1}|=|J_{2}|=2,J_{1}\neq J_{2},J_{1}\subset[N],J_{2}\subset[N]\}|d\sup\kappa\\
 & \leq(4d^{2})\cdot\left(\frac{N^{4}}{4}\right)\cdot d\sup\kappa.
\end{align*}
Theorem~\ref{thm:bound-density} follows from choosing
\begin{align*}
\mathsf{poly}(N,d) = N^4 d^3. 
\end{align*}

\subsection{Proof of Theorem \ref{thm:pers-betti-var}}
In this proof, we firstly define an auxiliary family of functions, and then verify the conditions in Theorem \ref{thm:discrimination}.
\paragraph{Defining the auxiliary function class.} For every $\bm{x} \in \Omega_{\ell}$ and $\mu \in \mathcal{Z}_{L,M}^q$, define 
\begin{align}
g_{\bm{x}}(\mu) = \mu(B_{\bm{x}})-\mathbb{E}[\mu](B_{\bm{x}}),
\end{align}
and let $\mathcal{G} = \{g_{\bm{x}}: \bm{x} \in \Omega_{\ell}\}$. It is easy to verify that $\mathbb{E}[g_{\bm{x}}(\mu)]=0$ for all $\bm{x} \in \Omega_{\ell}$, and that 
\begin{align*}
\sup_{\bm{x}\in \Omega_{\ell}} \left| \hat{\beta}_{\bm{x}} - \mathbb{E}[\hat{\beta}_{\bm{x}}]\right| = \left| \sup_{g \in \mathcal{G}} \frac{1}{n} \sum_{i=1}^n g(\mu_i)\right|.
\end{align*}
\paragraph{Bounding $||g_x||_{\infty}$ and $\mathbb{E}[g_x(\mu)^2]$.}
For any $\bm{x} \in \Omega_{\ell}$ , the set $B_{\bm{x}}$ is contained in $ \Omega_{\ell}$. Hence for any $\mu\in\mathcal{Z}_{L,M}^{q}$,
$\mu(B_{\bm{x}})$ and $\mathbb{E}[\mu](B_{\bm{x}})$ can be bounded
as 
\begin{align}
\mu(B_{\bm{x}}) & \leq\mu(\Omega_{\ell})\leq\ell^{-q}\mathsf{Pers}_{q}(\mu)\leq M\ell^{-q},\nonumber \\
\mathbb{E}[\mu](B_{\bm{x}}) & \leq\mathbb{E}[\mu](\Omega_{\ell})\leq\ell^{-q}\mathsf{Pers}_{q}(\mathbb{E}[\mu])\leq M\ell^{-q}.\label{eq:pers-betti-gx-var-muBx-bound}
\end{align}
Hence $\left\Vert g_{\bm{x}}\right\Vert _{\infty}$ can be bounded as
\begin{align}\label{eq:pers-betti-gx-infty}
&\left\Vert g_{\bm{x}}\right\Vert _{\infty}\leq\sup_{\mu\in\mathcal{Z}_{L,M}^{q}}\max\left\{ \mu(B_{\bm{x}}),\mathbb{E}[\mu](B_{\bm{x}})\right\} \leq M\ell^{-q}.
\end{align}
As for the variance of $g_x(\mu)$, we firstly observe that
\begin{align}\label{eq:pers-betti-gx-var}
& \mathbb{E}[g_{\bm{x}}(\mu)^2] \leq ||g_{\bm{x}}||_{\infty} \mathbb{E}[\mu](B_{\bm{x}})
\end{align}
Now, apart from using the bound $\mathbb{E}[\mu](B_{\bm{x}})\leq M\ell^{-q}$ from \eqref{eq:pers-betti-gx-var-muBx-bound}, we can also have tighter bound with respect to $\ell$ when $q>1$. To do this, we again take the coordinate transformation  
\begin{align*}
\begin{cases}
y_{1}=\frac{x_{2}-x_{1}}{\sqrt{2}}=\|\bm{x}-\partial\Omega\|_{2},\\
y_{2}=\frac{x_{2}+x_{1}}{\sqrt{2}}.
\end{cases}
\end{align*}
It can be easily verified that the determinant of the Jacobian
matrix between $\bm{x}$ and $\bm{y}$ coordinates is 1, and that
the $\Omega_{\ell}$ can be represented using $\bm{y}$ coordinates
by 
\[
\Omega_{\ell}=\left\{ (y_{1},y_{2}):\ell<y_{1}\leq\frac{L}{\sqrt{2}},y_{1}\leq y_{2}\leq\sqrt{2}L-y_{1}\right\} .
\]
Then, we have a tighter bound with respect to $\ell$ of $\mathbb{E}[\mu](B_{\bm{x}})$
when $q>1$ as 
\begin{align*}
\mathbb{E}[\mu](B_{\bm{x}}) & \leq\mathbb{E}[\mu](\Omega_{\ell})=\int_{\Omega_{\ell}}p(\bm{u})\mathrm{d}\bm{u}\\
 & =\int_{\Omega_{\ell}}\left\Vert \bm{u}-\partial\Omega\right\Vert _{2}^{-q}\bar{p}(\bm{u})\mathrm{d}\bm{u}\\
 & \leq\left\Vert \bar{p}\right\Vert _{\infty}\int_{\ell}^{\frac{L}{\sqrt{2}}}\left(\int_{y_{1}}^{\sqrt{2}L-y_{1}}\mathrm{d}y_{2}\right)y_{1}^{-q}\mathrm{d}y_{1}\\
 & \leq\left\Vert \bar{p}\right\Vert _{\infty}\int_{\ell}^{\frac{L}{\sqrt{2}}}\sqrt{2}Ly_{1}^{-q}\mathrm{d}y_{1}\\
 & \leq\frac{\sqrt{2}L\ell^{1-q}\left\Vert \bar{p}\right\Vert _{\infty}}{q-1}.
\end{align*}
Hence when we let $(q-1)_{+}=\max\{q-1,0\}$, 
\begin{equation}
\mathbb{E}[\mu](B_{\bm{x}})\leq\min\left\{ M\ell^{-q},\frac{\sqrt{2}L\ell^{1-q}\left\Vert \bar{p}\right\Vert _{\infty}}{(q-1)_{+}}\right\} .\label{eq:pers-betti-gx-var-EmuBx-bound}
\end{equation}
And hence by applying \eqref{eq:pers-betti-gx-var-EmuBx-bound} to \eqref{eq:pers-betti-gx-var}, the variance of $g_{x}(\mu)$ can be upper bounded as 
\begin{align}
\mathbb{E}[g_{\bm{x}}(\mu)^{2}] & \leq\left\Vert g_{\bm{x}}\right\Vert _{\infty}\mathbb{E}[\mu](B_{\bm{x}})\nonumber \\
 & \leq\min\left\{ M^{2}\ell^{-2q},\frac{\sqrt{2}ML\ell^{1-2q}\left\Vert \bar{p}\right\Vert _{\infty}}{(q-1)_{+}}\right\} 
\end{align}

\paragraph{Polynomial discrimination of $\mathcal{G}$.} By definition, the empirical persistent measure $\mu_i$ can be represented as
\begin{align*}
\mu_i = \sum_{j} \delta_{\bm{r}_{ij}},
\end{align*}
in which $\bm{r}_{ij} = (b_{ij},d_{ij})$ represents the $j$-th point in the corresponding persistent diagram, with $b_{ij}$ and $d_{ij}$ being its birth and death weight respectively . Without loss of generality, we can sort the points in descending order of their distance to the diagonal $\partial \Omega$. Let $N_i = \mu_i(\Omega_{\ell})$, then we have $N_i\leq M\ell^{-q}$.  Hence, for every $\bm{x}$ with $||\bm{x}-\partial \Omega||_2 = \ell$,  $\mu_i(B_{\bm{x}})$ can be represented as
\begin{align}
\mu_i(B_{\bm{x}}) = \sum_{j=1}^{N_i} \mathds{1}(b_{ij} < x_1) \mathds{1}(d_{ij} > x_2).
\end{align}
With this expression, we are ready to bound the cardinality of $\mathcal{G}(\bm{\mu}_1^n)$. Notice that for any fixed $\bm{x}$, the value of the tuple $(g_x(\mu_1),...,g_x(\mu_n))$ is completely decided by the Cartesian product of indicator functions
\begin{align*}
\{\mathds{1}(b_{ij}<x_1)\}_{i \in [n], j \in [N_i]} \times \{\mathds{1}(d_{ij}>x_2)\}_{i \in [n], j \in [N_i]} := S_b \times S_d.
\end{align*}
It is easy to see that with the variation of $\bm{x} = (x_1,x_2)$, the number of different values taken by $S_b$
 and $S_d$ can be bounded by 
 \begin{align*}
1 + \sum_{i=1}^n N_i \leq 1 + n \cdot M \ell^{-q}.
 \end{align*}
 Hence, the cardinality of $\mathcal{G}(\bm{\mu}_1^n)$ is bounded by
\begin{align}
\text{Card}(\mathcal{G}(\bm{\mu})) \leq \left(M\ell^{-q}n+1\right)^2.
\end{align}
\paragraph{Completing the proof.} The theorem is a direct result for applying Theorem \ref{thm:discrimination} with the following parameters:
\begin{align*}
\begin{cases}
A=M\ell^{-q};\\
B=M;\\
\sigma^{2}=\min\left\{ M^{2}\ell^{-2q},\frac{\sqrt{2}ML\ell^{1-2q}\left\Vert \bar{p}\right\Vert _{\infty}}{(q-1)_{+}}\right\};\\
\nu=2.
\end{cases}
\end{align*}

%% file: plots.tex
\section{Experimental details}\label{sec:experiments}

\begin{figure}[!t]
\centering
\begin{subfigure}[b]{0.48\textwidth}
         \centering
         \includegraphics[width=\textwidth]{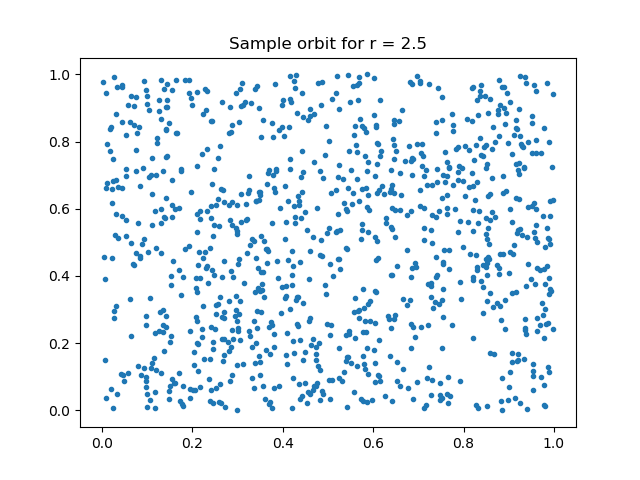}
     \end{subfigure}
     \begin{subfigure}[b]{0.48\textwidth}
         \centering
         \includegraphics[width=\textwidth]{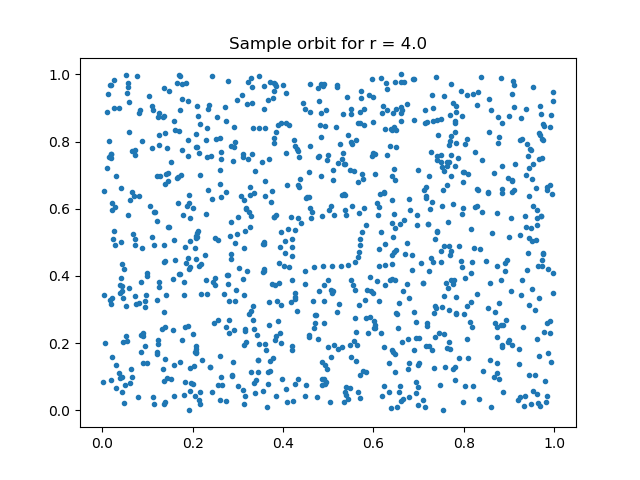}
     \end{subfigure}

     \begin{subfigure}[b]{0.48\textwidth}
\centering 
\includegraphics[width = \textwidth]{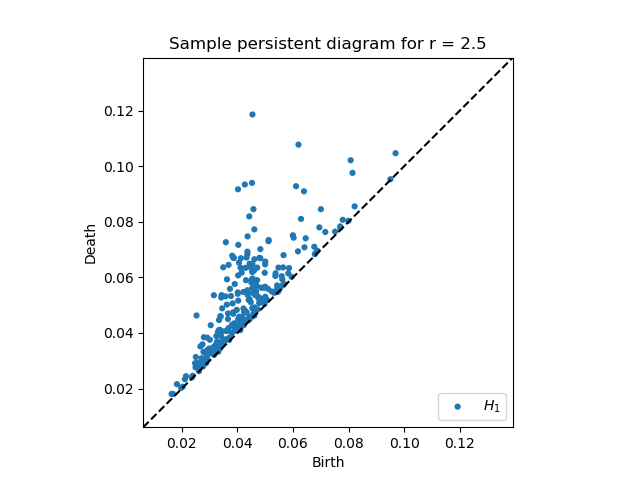}
\end{subfigure}
\begin{subfigure}[b]{0.48\textwidth}
\centering 
\includegraphics[width = \textwidth]{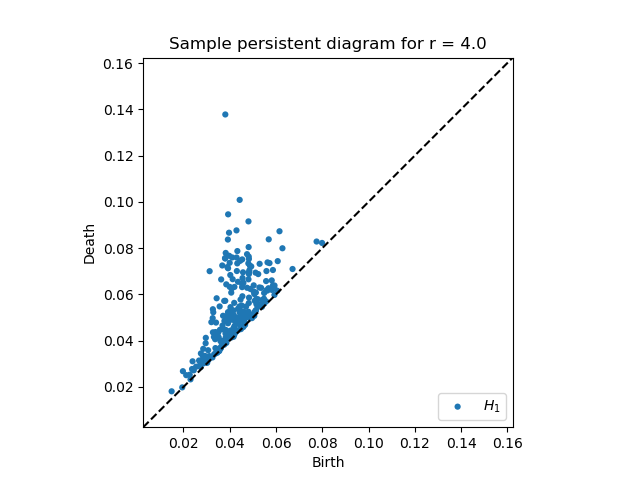}
\end{subfigure}
        \caption{Top row: sample orbits from the ORBIT5K data set with $r=2.5$ (left) and $r = 4.0$ (right). Bottom row: sample persistent diagrams.}
        \label{fig:ORBIT5K-orbits}
\end{figure}

\begin{figure}
\centering 
\begin{subfigure}[b]{0.48\textwidth}
\centering
\includegraphics[width = \textwidth]{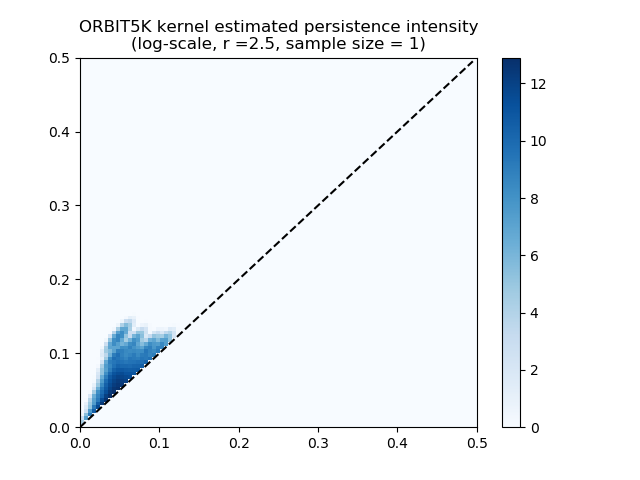}
\end{subfigure}
\begin{subfigure}[b]{0.48\textwidth}
\centering
\includegraphics[width = \textwidth]{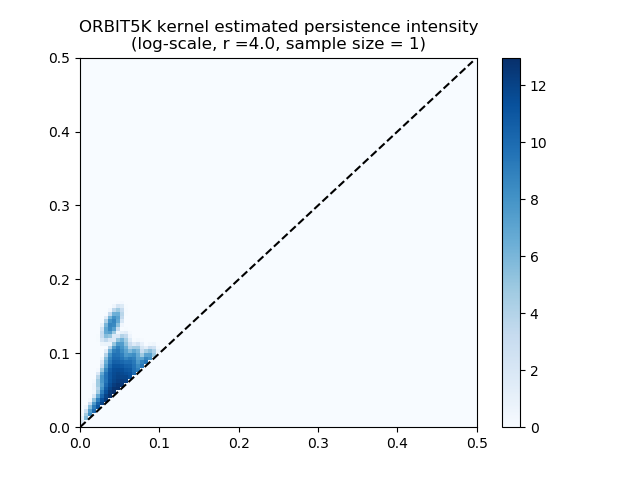}
\end{subfigure}

\begin{subfigure}[b]{0.48\textwidth}
\centering
\includegraphics[width = \textwidth]{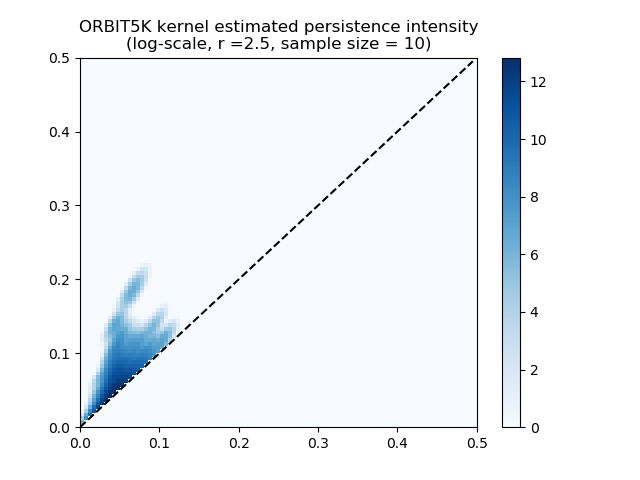}
\end{subfigure}
\begin{subfigure}[b]{0.48\textwidth}
\centering
\includegraphics[width = \textwidth]{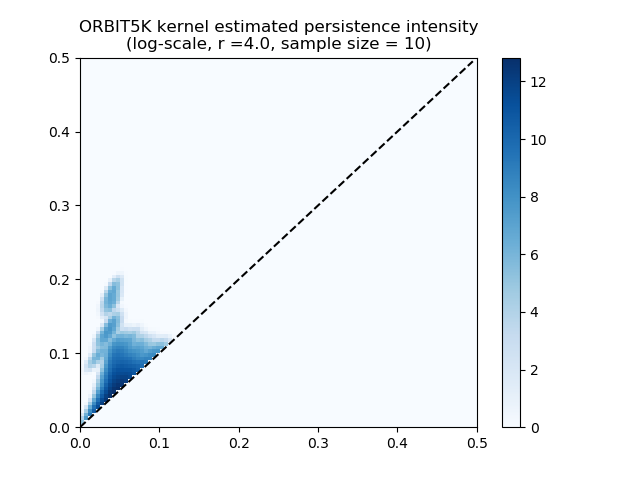}
\end{subfigure}

\begin{subfigure}[b]{0.48\textwidth}
\centering
\includegraphics[width = \textwidth]{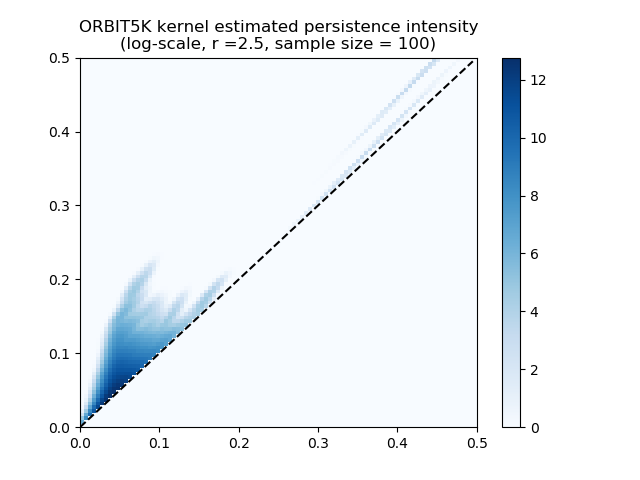}
\end{subfigure}
\begin{subfigure}[b]{0.48\textwidth}
\centering
\includegraphics[width = \textwidth]{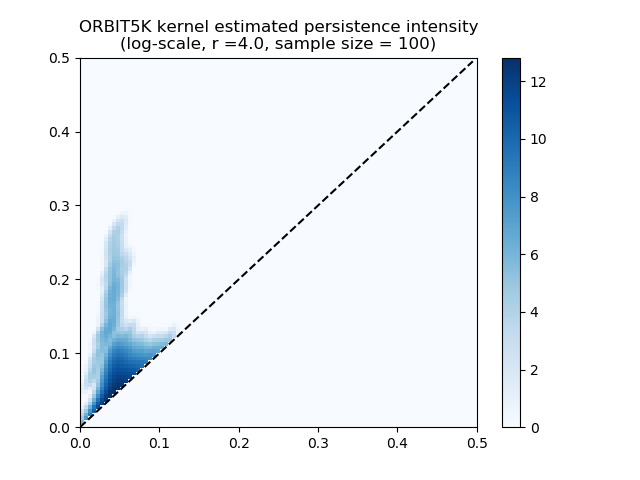}
\end{subfigure}

\begin{subfigure}[b]{0.48\textwidth}
\centering
\includegraphics[width = \textwidth]{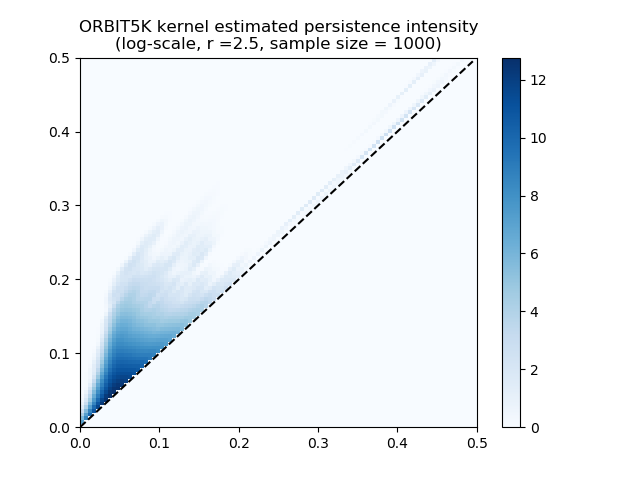}
\end{subfigure}
\begin{subfigure}[b]{0.48\textwidth}
\centering
\includegraphics[width = \textwidth]{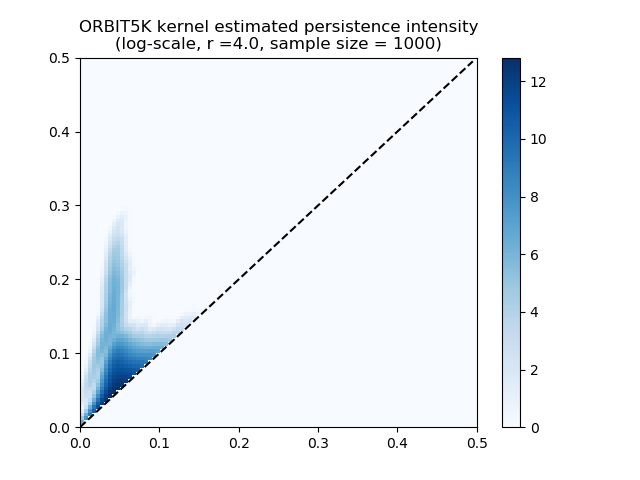}
\end{subfigure}
\caption{Kernel estimators for the persistence intensity function from the ORBIT5K data set with $r=2.5$ (left) and $r=4.0$ (right) and sample sizes $1$, $10$, $100$ and $1000$ (top to bottom).}
\label{fig:ORBIT5K-KIE}
\end{figure}

\begin{figure}
\begin{subfigure}[b]{0.48\textwidth}
\centering
\includegraphics[width = \textwidth]{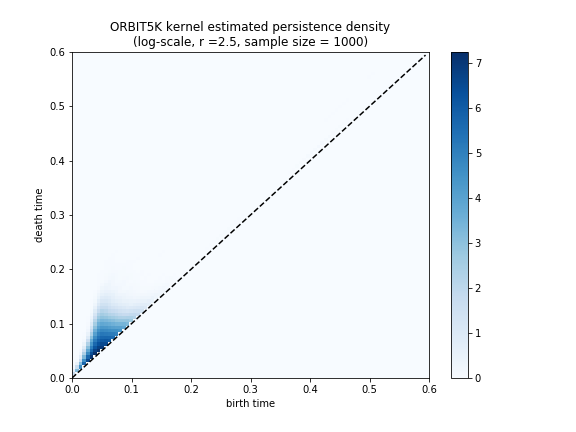}
\end{subfigure}
\begin{subfigure}[b]{0.48\textwidth}
\centering
\includegraphics[width = \textwidth]{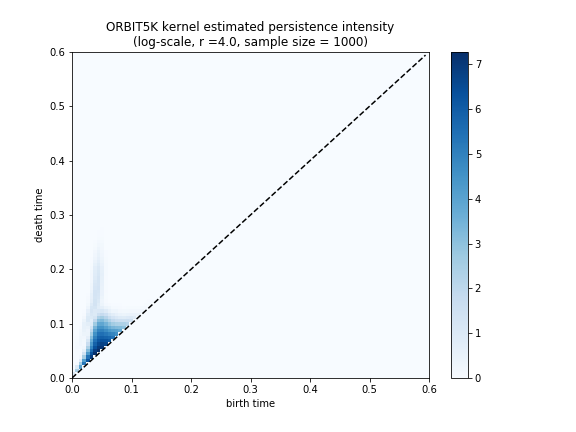}
\end{subfigure}
\caption{Kernel estimators for the persistence density function from the ORBIT5K data set with $r=2.5$ (left) and $r=4.0$ (right) and sample size $n = 1000$.}
\label{fig:ORBIT5K-KDE}
\end{figure}

\begin{figure}
\begin{subfigure}[b]{0.48\textwidth}
\centering
\includegraphics[width = \textwidth]{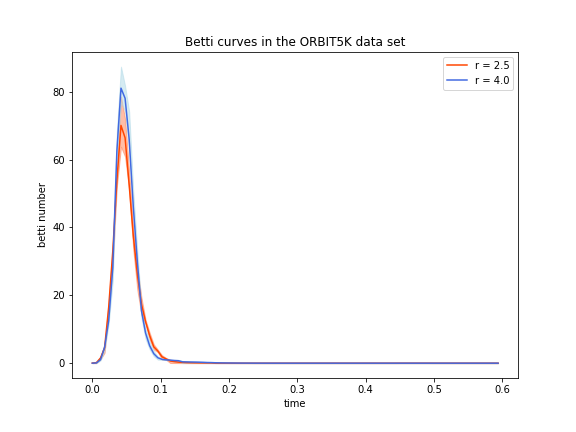}
\end{subfigure}
\begin{subfigure}[b]{0.48\textwidth}
\centering
\includegraphics[width = \textwidth]{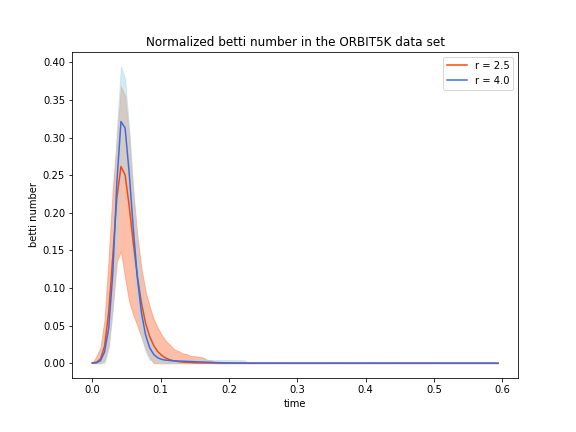}
\end{subfigure}
\caption{Empirical betti curves (left) and normalized betti curves (right) from the ORBIT5K data set with $r=2.5$ and $r=4.0$. Solid lines show sample average and the shades depict the lower and upper 2.5 percentiles.}
\label{fig:empirical-Betti-ORBIT5K}
\end{figure}

\begin{figure}
\begin{subfigure}[b]{0.48\textwidth}
\centering
\includegraphics[width = \textwidth]{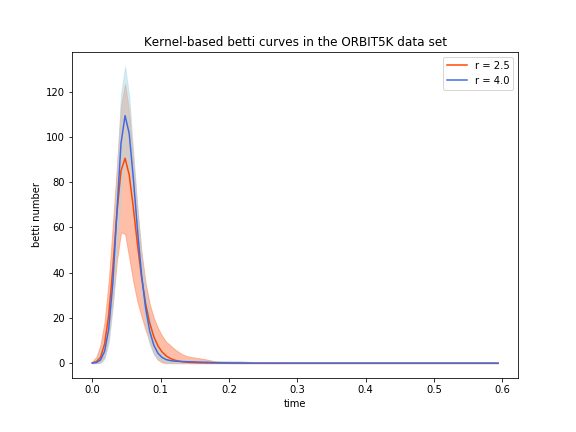}
\end{subfigure}
\begin{subfigure}[b]{0.48\textwidth}
\centering
\includegraphics[width = \textwidth]{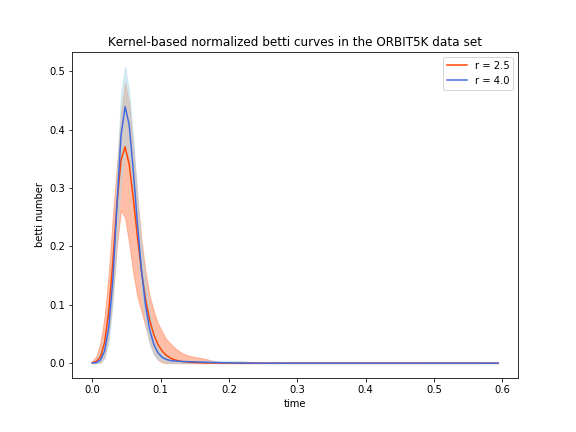}
\end{subfigure}
\caption{Kernel-based betti curves (left) and normalized betti curves (right) from the ORBIT5K data set with $r=2.5$ and $r=4.0$. Solid lines show sample average and the shades depict the lower and upper 2.5 percentiles.}
\label{fig:Betti-ORBIT5K}
\end{figure}

Figure \ref{fig:ORBIT5K-orbits} shows two ORBIT5K simulations with different values of $r$ ($2.5$ and $4$) and the corresponding persistent diagrams.
 Figure \ref{fig:ORBIT5K-KIE} displays the kernel intensity functions for the ORBIT5K simulations
set with $r=2.5$ and $r =4 $ for varying sample sizes, while Figure \ref{fig:ORBIT5K-KDE} shows persistence density functions. Figures \ref{fig:empirical-Betti-ORBIT5K} and \ref{fig:Betti-ORBIT5K} show the Betti curves and estimated Betti curves using the kernel density function for the ORBIT5K simulations for $r=2.5$ and $r =4 $.

\begin{figure}
\centering 
\begin{subfigure}[b]{0.48\textwidth}
\centering
\includegraphics[width = \textwidth]{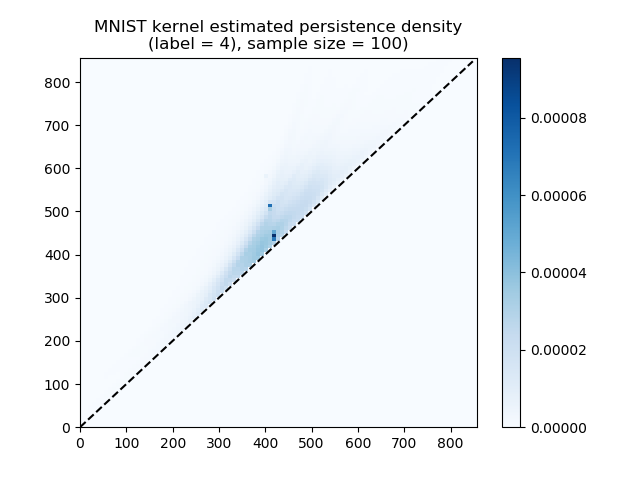}
\end{subfigure}
\begin{subfigure}[b]{0.48\textwidth}
\centering
\includegraphics[width = \textwidth]{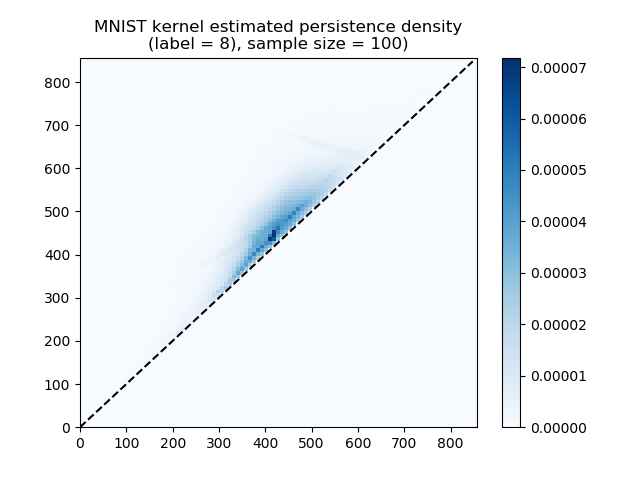}
\end{subfigure}

\begin{subfigure}[b]{0.48\textwidth}
\centering
\includegraphics[width = \textwidth]{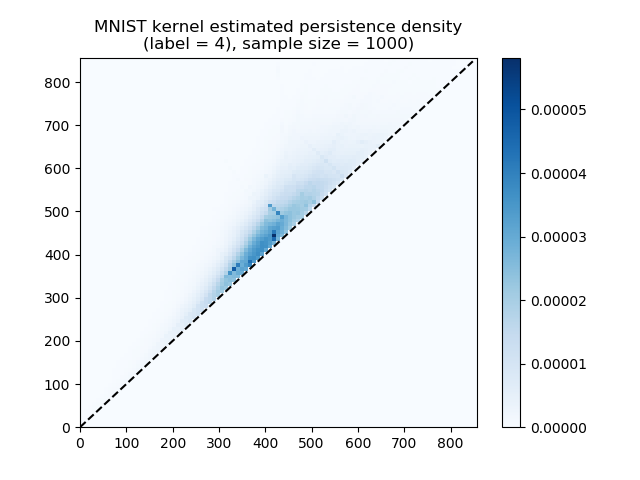}
\end{subfigure}
\begin{subfigure}[b]{0.48\textwidth}
\centering
\includegraphics[width = \textwidth]{imgs/MNIST/KDE-label4-n1000.png}
\end{subfigure}

\begin{subfigure}[b]{0.48\textwidth}
\centering
\includegraphics[width = \textwidth]{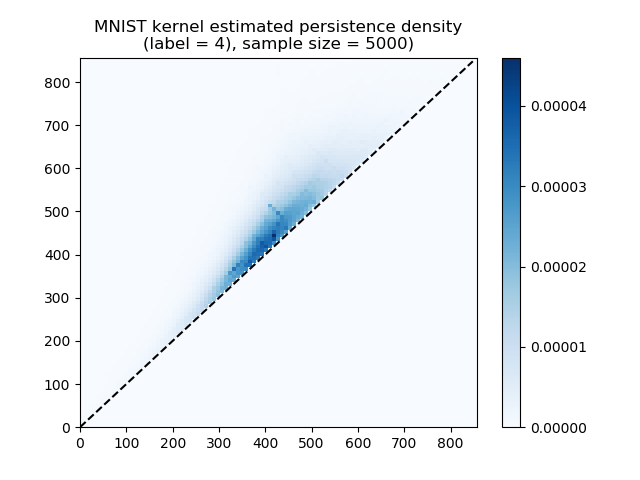}
\end{subfigure}
\begin{subfigure}[b]{0.48\textwidth}
\centering
\includegraphics[width = \textwidth]{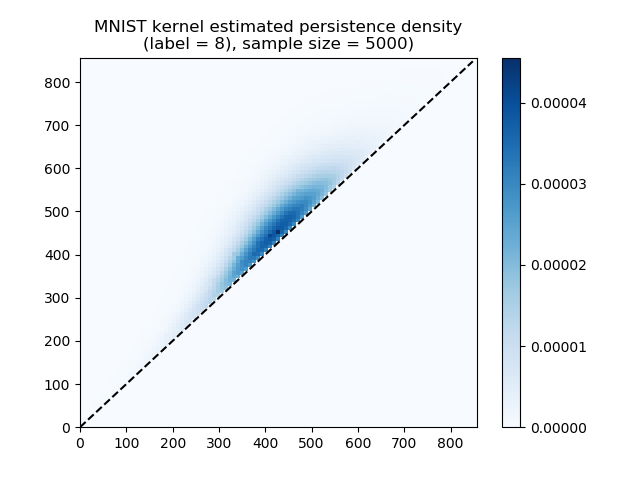}
\end{subfigure}
\caption{Kernel estimators for the persistence density function from the MNIST data set for the digits 4 (left column) and 8 (right column) based on random draws of sample sizes $100$, $1000$ and $5000$ (top to bottom).}\label{fig:MNIST-KDE}
\end{figure}

Figure \ref{fig:MNIST-KDE} displays the estimated persistence density functions computed over  random draws of varying size of the digits ``4" and ``8" from the MNIST dataset.

\begin{figure}[h!]
\centering
\begin{subfigure}{0.43\textwidth}
\centering
\includegraphics[width=\linewidth]{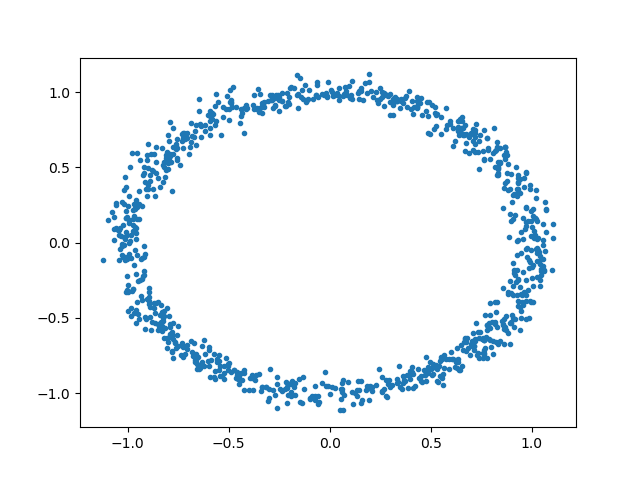}
\end{subfigure}
\hfill
\begin{subfigure}{0.43\textwidth}
\centering
\includegraphics[width=\linewidth]{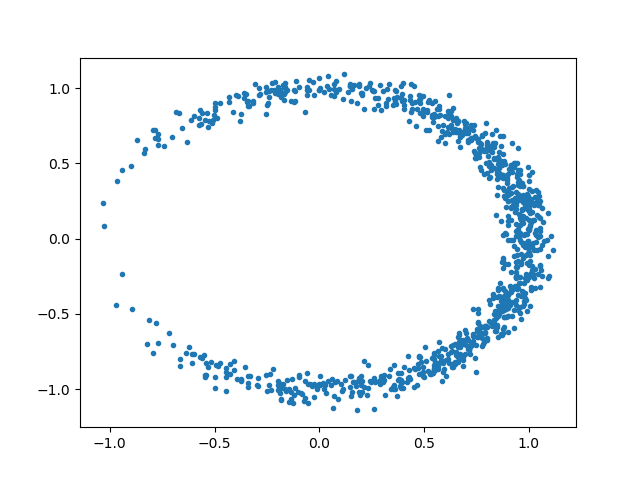}
\end{subfigure}

\begin{subfigure}{0.43\textwidth}
\centering
\includegraphics[width = \linewidth]{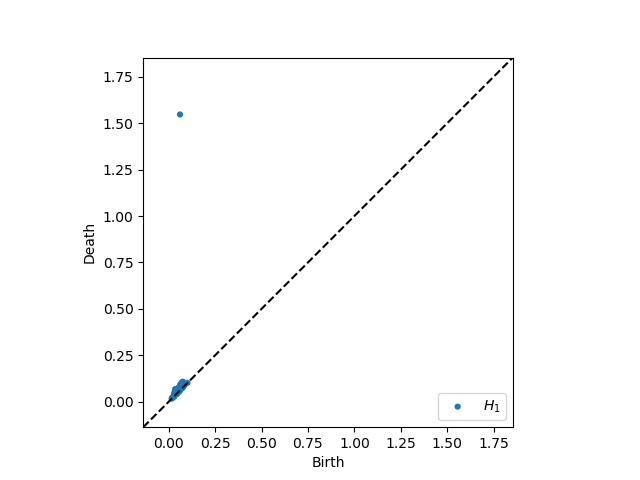}
\end{subfigure}
\hfill 
\begin{subfigure}{0.43\textwidth}
\centering
\includegraphics[width = \linewidth]{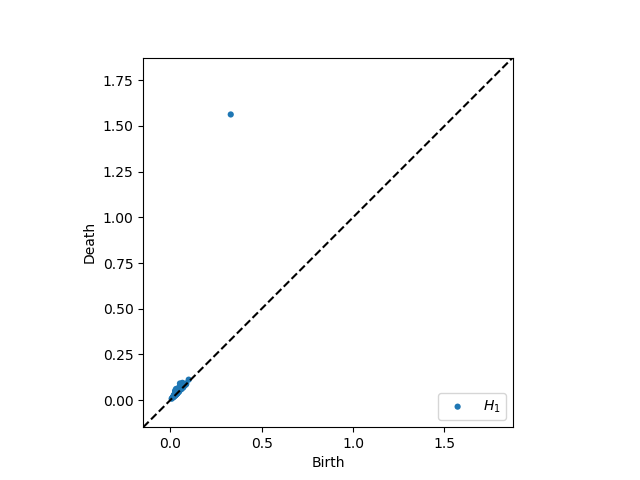}
\end{subfigure}

\begin{subfigure}{0.43\textwidth}
\centering
\includegraphics[width = \linewidth]{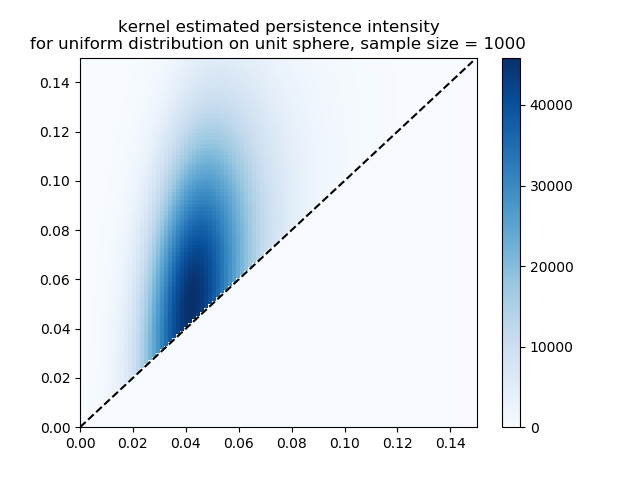}
\end{subfigure}
\hfill 
\begin{subfigure}{0.43\textwidth}
\centering
\includegraphics[width = \linewidth]{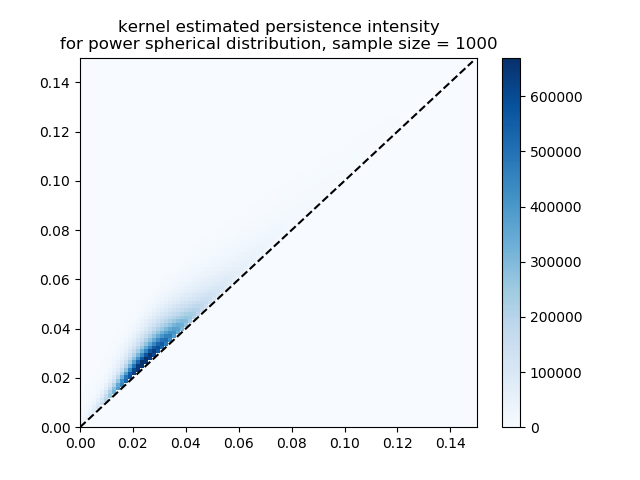}
\end{subfigure}

\begin{subfigure}{0.43\textwidth}
\centering
\includegraphics[width = \linewidth]{imgs/sphere/unif1000_density.png}
\end{subfigure}
\hfill 
\begin{subfigure}{0.43\textwidth}
\centering
\includegraphics[width = \linewidth]{imgs/sphere/power1000_density.png}
\end{subfigure}

\caption{Sample plot, sample persistence diagram, kernel-estimated persistence intensity and density functions (top to bottom) for the uniform distribution (left) and the power spherical distribution \cite{decao2020power} (right) on the unit circle $\mathbb{S}^1$. The parameters for the power spherical distribution are set to $\mu = \frac{\pi}{2}$ and $\kappa = 1$. Each sample contains $1000$ points generated in an i.i.d. manner from the corresponding distributions on the unit circle and each perturbed by an additive noise term, sampled i.i.d. from the  $N(0,0.05^2\bm{I}_2)$ distribution. The intensity and density functions are estimated through $1000$ independent samples.} 
\label{fig:sphere_app}
\end{figure}

Finally, Figure \ref{fig:sphere_app} shows the sample plots, sample persistent diagram, and the kernel-based estimators for the persistence intensity and density functions calculated from 1000 samples of 1000 data points drawn from the uniform distirbution and the power sphere distribution on the unit circle.

%% file: arxiv_version.bbl
\begin{thebibliography}{18}
\providecommand{\natexlab}[1]{#1}
\providecommand{\url}[1]{\texttt{#1}}
\expandafter\ifx\csname urlstyle\endcsname\relax
  \providecommand{\doi}[1]{doi: #1}\else
  \providecommand{\doi}{doi: \begingroup \urlstyle{rm}\Url}\fi

\bibitem[Adams et~al.(2017)Adams, Emerson, Kirby, Neville, Peterson, Shipman,
  Chepushtanova, Hanson, Motta, and Ziegelmeier]{adams2017persistence}
Henry Adams, Tegan Emerson, Michael Kirby, Rachel Neville, Chris Peterson,
  Patrick Shipman, Sofya Chepushtanova, Eric Hanson, Francis Motta, and Lori
  Ziegelmeier.
\newblock Persistence images: A stable vector representation of persistent
  homology.
\newblock \emph{Journal of Machine Learning Research}, 18, 2017.

\bibitem[Cao and Aziz(2020)]{decao2020power}
Nicola~De Cao and Wilker Aziz.
\newblock The power spherical distribution, 2020.

\bibitem[Chazal and Divol(2019)]{chazal2018density}
Fr{\'e}d{\'e}ric Chazal and Vincent Divol.
\newblock The density of expected persistence diagrams and its kernel based
  estimation.
\newblock \emph{Journal of Computational Geometry}, 10\penalty0 (2), 2019.

\bibitem[Chazal and Michel(2021{\natexlab{a}})]{DBLP:journals/frai/ChazalM21}
Fr{\'{e}}d{\'{e}}ric Chazal and Bertrand Michel.
\newblock An introduction to topological data analysis: Fundamental and
  practical aspects for data scientists.
\newblock \emph{Frontiers Artif. Intell.}, 4:\penalty0 667963,
  2021{\natexlab{a}}.

\bibitem[Chazal et~al.(2013)Chazal, Fasy, Lecci, Rinaldo, and
  Wasserman]{Chazal2013StochasticCO}
Fr{\'e}d{\'e}ric Chazal, Brittany~Terese Fasy, Fabrizio Lecci, Alessandro
  Rinaldo, and Larry~A. Wasserman.
\newblock Stochastic convergence of persistence landscapes and silhouettes.
\newblock \emph{Proceedings of the thirtieth annual symposium on Computational
  geometry}, 2013.

\bibitem[Chazal and Michel(2021{\natexlab{b}})]{fredmichelreview}
Frédéric Chazal and Bertrand Michel.
\newblock An introduction to topological data analysis: Fundamental and
  practical aspects for data scientists.
\newblock \emph{Frontiers in Artificial Intelligence}, 4, 2021{\natexlab{b}}.

\bibitem[Chen et~al.(2015)Chen, Wang, Rinaldo, and
  Wasserman]{original.intensity}
Yen-Chi Chen, Daren Wang, Alessandro Rinaldo, and Larry Wasserman.
\newblock Astatistical analysis of persistence intensity functions.
\newblock \emph{https://arxiv.org/pdf/1510.02502.pdf}, 2015.

\bibitem[Cohen-Steiner et~al.(2010)Cohen-Steiner, Edelsbrunner, Harer, and
  Mileyko]{cohen2010lipschitz}
David Cohen-Steiner, Herbert Edelsbrunner, John Harer, and Yuriy Mileyko.
\newblock Lipschitz functions have l p-stable persistence.
\newblock \emph{Foundations of computational mathematics}, 10\penalty0
  (2):\penalty0 127--139, 2010.

\bibitem[Divol and Lacombe(2021)]{divol2021estimation}
Vincent Divol and Th{\'e}o Lacombe.
\newblock Estimation and quantization of expected persistence diagrams.
\newblock In \emph{International Conference on Machine Learning}, pages
  2760--2770. PMLR, 2021.

\bibitem[Divol and Polonik(2019)]{divol2019choice}
Vincent Divol and Wolfgang Polonik.
\newblock On the choice of weight functions for linear representations of
  persistence diagrams.
\newblock \emph{Journal of Applied and Computational Topology}, 3\penalty0
  (3):\penalty0 249--283, 2019.

\bibitem[Gin{\'e} and Nickl(2021)]{gine2021mathematical}
Evarist Gin{\'e} and Richard Nickl.
\newblock \emph{Mathematical foundations of infinite-dimensional statistical
  models}.
\newblock Cambridge university press, 2021.

\bibitem[Kim et~al.(2020)Kim, Kim, Zaheer, Kim, Chazal, and
  Wasserman]{NEURIPS2020_b803a925}
Kwangho Kim, Jisu Kim, Manzil Zaheer, Joon Kim, Frederic Chazal, and Larry
  Wasserman.
\newblock Pllay: Efficient topological layer based on persistent landscapes.
\newblock In H.~Larochelle, M.~Ranzato, R.~Hadsell, M.F. Balcan, and H.~Lin,
  editors, \emph{Advances in Neural Information Processing Systems}, volume~33,
  pages 15965--15977. Curran Associates, Inc., 2020.

\bibitem[Kusano et~al.(2016)Kusano, Hiraoka, and Fukumizu]{pmlr-v48-kusano16}
Genki Kusano, Yasuaki Hiraoka, and Kenji Fukumizu.
\newblock Persistence weighted gaussian kernel for topological data analysis.
\newblock In Maria~Florina Balcan and Kilian~Q. Weinberger, editors,
  \emph{Proceedings of The 33rd International Conference on Machine Learning},
  volume~48, pages 2004--2013, 20--22 Jun 2016.

\bibitem[Morgan(2016)]{morgan2016geometric}
Frank Morgan.
\newblock \emph{Geometric measure theory: a beginner's guide}.
\newblock Academic press, 2016.

\bibitem[Nietert et~al.(2021)Nietert, Goldfeld, and Kato]{nietert2021smooth}
Sloan Nietert, Ziv Goldfeld, and Kengo Kato.
\newblock Smooth $ p $-wasserstein distance: Structure, empirical
  approximation, and statistical applications.
\newblock In \emph{International Conference on Machine Learning}, pages
  8172--8183. PMLR, 2021.

\bibitem[Peyre(2018)]{peyre2018comparison}
R{\'e}mi Peyre.
\newblock Comparison between $w_2$ distance and $\dot{H}^{-1}$ norm, and
  localization of wasserstein distance.
\newblock \emph{ESAIM: Control, Optimisation and Calculus of Variations},
  24\penalty0 (4):\penalty0 1489--1501, 2018.

\bibitem[Steinwart and Christmann(2008)]{steinwart2008support}
Ingo Steinwart and Andreas Christmann.
\newblock \emph{Support vector machines}.
\newblock Springer Science \& Business Media, 2008.

\bibitem[Wilding et~al.(2021)Wilding, Nevenzeel, van de Weygaert, Vegter,
  Pranav, Jones, Efstathiou, and Feldbrugge]{10.1093/mnras/stab2326}
Georg Wilding, Keimpe Nevenzeel, Rien van de Weygaert, Gert Vegter, Pratyush
  Pranav, Bernard J~T Jones, Konstantinos Efstathiou, and Job Feldbrugge.
\newblock {Persistent homology of the cosmic web – I. Hierarchical topology
  in $\Lambda$CDM cosmologies}.
\newblock \emph{Monthly Notices of the Royal Astronomical Society},
  507\penalty0 (2):\penalty0 2968--2990, 2021.

\end{thebibliography}
